\documentclass[10pt]{amsart} 
 \normalfont

\usepackage[psamsfonts]{amssymb}
\usepackage{pifont}
\theoremstyle{plain}

\newtheorem{theorem}{Theorem}[section]
\newtheorem{lemma}[theorem]{Lemma}
\newtheorem{proposition}[theorem]{Proposition}
\newtheorem{corollary}[theorem]{Corollary}

\newtheorem*{cor}{Corollary}

\newtheorem*{lem}{Lemma}

\def\bc{\mathbb{C}}

\def\bq{\mathbb{Q}}
\def\bz{\mathbb{Z}}
\def\br{\mathbb{R}}
\def\bh{\mathbb{H}}

\def\fF{\mathfrak{F}}
\def\fB{\mathfrak{B}}
\def\fC{\mathfrak{C}}

\def\fU{\mathfrak{U}}
\def\fG{\mathfrak{G}}
\def\fH{\mathfrak{H}}

\def\fL{\mathfrak{L}}

\def\fN{\mathfrak{N}}
\def\fP{\mathfrak{P}}

\def\fU{\mathfrak{U}}

\def\fV{\mathfrak{V}}

\def\fZ{\mathfrak{Z}}
\def\fs{\mathfrak{s}}
\def\fl{\mathfrak{l}}

\def\iff{\Leftrightarrow}

\def\f2q{($\fF_2$(q)}
\def\imp{$\Rightarrow$}
\def\ip{\langle ~, ~\rangle}
\def\nvg{$\fN$ = V$\oplus \fG$}
\def\nug{$\fN$ = U$\oplus \fG_0$}

\def\ws{V = V$_0$ + \(\sum_{\lambda \epsilon \Lambda}\) V$_\lambda$}
\def\rs{$\fG$ = $\fH$ + \(\sum_{\alpha \epsilon \Phi}$ $\fG_\alpha$}

\def\submu{$_{\mu}$}

\def\12{{\scriptstyle{\frac12}}}

\def\<{\langle} \def\>{\rangle}

\newcounter{commentlabel}
\begin{document}
\title{2-step nilpotent Lie groups arising from semisimple modules}
\author{Patrick Eberlein \\ University of North Carolina at Chapel Hill}
\date{\today}
\maketitle

$\mathbf{Abstract}$  Let $\fG_{0}$ denote a compact semisimple Lie algebra and U a finite dimensional real $\fG_{0}$ - module.  The vector space $\fN_{0} = U \oplus \fG_{0}$ admits a canonical 2-step nilpotent Lie algebra structure with $[\fN_{0}, \fN_{0}] = \fG_{0}$ and an inner product $\langle , \rangle$, unique up to scaling, for which the elements of $\fG_{0}$ are skew symmetric derivations of $\fN_{0}$.  Let N$_{0}$ denote the corresponding simply connected 2-step nilpotent Lie group with Lie algebra $\fN_{0}$, and let $\langle , \rangle$ also denote the left invariant metric on N$_{0}$ determined by the inner product $\langle , \rangle$ on $\fN_{0}$.   In this article we investigate the basic differential geometric properties of N$_{0}$ by using elementary representation theory to study the complexification $\fN = \fN_{0}^{\bc} = V \oplus \fG$, where V $= U^{\bc}$ and $\fG = \fG_{0}^{\bc}$.  The weight space decomposition for V and its real analogue for U describe the bracket structures for $\fN = \fN_{0}^{\bc}$ and $\fN_{0}$.  The Weyl group W of $\fG$ acts on the real Lie algebra $\fN_{0}$ by automorphisms and isometries.  The Lie algebra $\fN_{0}$ admits a Chevalley rational structure for which the the weight spaces of U are rational.  We use the roots of $\fG$ and the weights of V to construct totally geodesic, rational subalgebras of $\fN_{0} = U \oplus \fG_{0}$.
\newline

$\mathbf{Mathematics~2000~ Subject~Classification}$

Primary : 53C30 Secondary : 22E25, 22E46
\newline

$\mathbf{Key~words~and~phrases}$
2-step nilpotent, compact semisimple Lie algebra, Chevalley rational structure, rationality of weight spaces, totally geodesic and rational subalgebras, admissible abstract weights

\section*{Table~ of~ Contents}

Introduction

Section 1 \hspace{.2in}  Notation and preliminaries

Section 2  \hspace{.2in}  Some 2-step nilpotent Lie algebras arising from representations

\hspace{.75in}  Complex 2-step nilpotent Lie algebras

\hspace{.75in}  Uniqueness of the complex Lie algebra structure

\hspace{.75in}  Real 2-step nilpotent Lie algebras

\hspace{.75in}  Uniqueness of the real Lie algebra structure

\hspace{.75in}  Commutativity of irreducible submodules

Section 3 \hspace{.2in}  Chevalley bases and the universal enveloping algebra

\hspace{.75in}  Chevalley bases and $\fG$ - automorphisms of order 2.

\hspace{.75in}  Universal enveloping algebra and the subring $\fU(\fG)_{\bz}$

\hspace{.75in}  Chevalley bases of compact real forms $\fG_{0}$

Section 4 \hspace{.2in}  Weight space decomposition of a complex $\fG$ - module

\hspace{.75in}  Some relations between roots and weights

\hspace{.75in}  $\fG_{0}$ - invariant Hermitian inner products

\hspace{.75in}  $\bz$ - bases adapted to weight spaces

\hspace{.75in}  Orthogonality of weight spaces

\hspace{.75in}  Root space decomposition of $\fG$

\hspace{.75in}  Bracket relations in $\fN = V \oplus \fG$

\hspace{.75in}  The range of ad v : v $\in$ V

\hspace{.75in}  Surjectivity of ad v

\hspace{.75in}  The range of ad v$_{\lambda} : v_{\lambda} \in V_{\lambda}$

\hspace{.75in}  Automorphisms and derivations of $\fN = V \oplus \fG$

Section 5 \hspace{.2in}  Action of the Weyl group by unitary automorphisms

\hspace{.75in}  The group G$_{\bz}$

\hspace{.75in}  The Weyl group and complex $\fG$ - modules

\hspace{.75in}  The Weyl group and real $\fG_{0}$ - modules

\hspace{.75in}  The Weyl group and Aut($\fN$)

Section 6 \hspace{.2in}  Weight space decomposition of a real $\fG_{0}$ - module

\hspace{.75in}  Conjugations

\hspace{.75in}  Weight space decomposition of U

\hspace{.75in}  Root space decomposition of $\fG_{0}$

\hspace{.75in}  Bracket relations in $\fN_{0} = U_{0} \oplus \fG_{0}$

\hspace{.75in}  Abstract weights and real weight vectors

\hspace{.75in}  Relationship to complex weight vectors

\hspace{.75in}  Rationality of weight vectors

\hspace{.75in}  The range of ad u : u $\in$ U.

\hspace{.75in}  The range of ad u$_{\lambda} : u_{\lambda} \in U_{\lambda}$

\hspace{.75in}  Automorphisms and derivations of $\fN_{0} = U \oplus \fG_{0}$

\hspace{.75in}  The Weyl group and Aut($\fN_{0}$)

Section 7 \hspace{.2in} Adapted bases of $\fG_{0}$

\hspace{.75in}  The elements $A_{\beta},   B_{\beta} : \beta \in \Phi$

\hspace{.75in}  Bracket relations between $A_{\beta},   B_{\beta}$

\hspace{.75in}  Kernel of $A_{\beta},   B_{\beta}$

\hspace{.75in}  Range of $A_{\beta},   B_{\beta}$

\hspace{.75in}  Nonsingular subspaces for $A_{\beta},   B_{\beta}$

\hspace{.75in}  $\fG_{0}$ - modules with nontrivial zero weight space

Section 8 \hspace{.2in}  Rational structures on $\fN_{0} = U \oplus \fG_{0}$

\hspace{.75in}  Existence of rational structures

\hspace{.75in}  Rationality of weight spaces

Section 9 \hspace{.2in}  Admissible abstract weights

\hspace{.75in}  Simple Lie algebras

\hspace{.75in}  Semisimple Lie algebras

Section 10 \hspace{.15in} Totally geodesic subalgebras

\hspace{.75in}  A criterion for a subalgebra to be totally geodesic

\hspace{.75in}  Rational totally geodesic subalgebras with 1-dimensional center

\hspace{.75in}  Rational totally geodesic subalgebras of quaternionic type

Section 11 \hspace{.15in} Appendix

\hspace{.75in}  Proof of (9.1)

\hspace{.75in}  Proof of (9.2)

\section*{Introduction}

    Let N$_{0}$ denote a simply connected, 2-step nilpotent Lie group with a left invariant metric $\ip$.  and let $\fN_{0}$ denote the Lie algebra of N$_{0}$.  Write $\fN_{0} = \fV_{0} \oplus \fZ_{0}$, where $\fZ_{0}$ is the center of $\fN_{0}$ and $\fV_{0}$ its orthogonal complement.  If Ric : $\fN_{0}$ x $\fN_{0} \rightarrow \br$ denotes the Ricci tensor, then it is well known that Ric($\fV_{0}$, $\fZ_{0}$) = $\{0\}$, Ric is negative definite on $\fV_{0}$ and positive semidefinite on $\fZ_{0}$.  If N$_{0}$ has no Euclidean de Rham factor then Ric is positive definite on $\fZ_{0}$.  (See for example [Eb 3].)  This means that the comparison methods, which have been used so widely to study Riemannian manifolds with curvature of a fixed sign, do not apply in this situation.  In compensation, one has the usual algebraic advantages common to all homogeneous spaces and a few others that are specific to 2-step nilpotent Lie algebras.

    One of the most interesting classes of examples of 2-step nilpotent Lie groups with a left invariant metric arises from finite dimensional real representations of compact semisimple Lie algebras.  These examples are central to this article.

    Let $\fG_{0}$ be a compact, semisimple real Lie algebra, and let U be a finite dimensional real $\fG_{0}$ - module.  Let $\langle , \rangle_{U}$ and $\langle , \rangle_{\fG_{0}}$ be inner products on U and $\fG_{0}$ such that X : U $\rightarrow$ U and ad X : $\fG_{0} \rightarrow \fG_{0}$ are skew symmetric for all X $\in \fG_{0}$.  The orthogonal direct sum $\fN_{0} = U \oplus \fG_{0}$ admits a 2-step nilpotent Lie algebra structure such that [U, U] $= \fG_{0}$ and $\fG_{0}$ is contained in the center of $\fN_{0}$.  Given u$_{1}$, u$_{2} \in$ U we define [u$_{1}$ , u$_{2}$] $\in \fG_{0}$ by the condition that $\langle [u_{1} , u_{2}] , X \rangle_{\fG_{0}} = \langle X(u_{1}), u_{2} \rangle_{U}$ for all X $\in \fG_{0}$.  The isomorphism type of $\{\fN_{0}, \ip \}$ is independent of the inner products $\langle , \rangle_{U}$ and $\langle , \rangle_{\fG_{0}}$.   The metric Lie algebra $\{\fN_{0}, \ip\}$ determines a simply connected 2-step nilpotent Lie group N$_{0}$ with Lie algebra $\fN_{0}$ and a left invariant Riemannian metric $\ip$.  See section 2.

    The goal of this article is to develop the basic structure of the Lie algebra $\{\fN_{0}, \ip \}$ in order to study the geometry of $\{N_{0}, \ip \}$.  If $\fN, \fG$ and V denote the complexifications of  $\fN_{0}, \fG_{0}$ and U, then $\fG$ is a complex semisimple Lie algebra, V is a $\fG$ -module and $\fN = V \oplus \fG$ as a vector space.  Our approach is to use elementary representation theory to study the geometry and basic structure of $\fN = V \oplus \fG$ and then deduce the geometry and basic structure of the real form $\fN_{0} = U \oplus \fG_{0}$.  Some of the complex representation theory of $\fG$ and its interplay with the real representation theory of $\fG_{0}$ does not seem to be easily locatable in the literature.  Where needed, we develop this theory.

    If $\fH$ is a Cartan subalgebra of $\fG$, then $\fH$ determines a Chevalley basis $\fC = \{X_{\beta} : \beta \in \Phi ; \tau_{\alpha} : \alpha \in \Delta \}$ that is unique up to a certain type of rescaling.  Here $\Phi \subset \fH^{*}$ denotes the roots determined by $\fH, \Delta$ is a base for the positive roots and X$_{\beta}$ is a suitably chosen element of the 1-dimensional root space $\fG_{\beta}$.  The vectors
$\{ \tau_{\alpha} : \alpha \in \Delta\}$ are normalized root vectors and form a $\bc$ - basis of $\fH$.  Each Chevalley basis $\fC$ has structure constants in $\bz$ and determines in a canonical way a basis $\fB$ of each $\fG$ - module   V such that $\fC$ leaves invariant $\bz$ - span ($\fB$).  The set $\fL = \fC \cup \fB$ is a $\bc$ - basis for $\fN = V \oplus \fG$ whose structure constants are integers.

    All of this structure has an analogue for the real form $\fN_{0} = U \oplus \fG_{0}$ of $\fN$.  If $\fH$ is a Cartan subalgebra of $\fG = \fG_{0}^{\bc}$ that is the complexification of a maximal abelian subalgebra $\fH_{0}$ of $\fG_{0}$, then we say that $\fH$ is $\mathit{adapted~ to}~ \fG_{0}$.  If $\fH$ is a Cartan subalgebra  adapted to $\fG_{0}$, then the Chevalley basis of $\fG$ determines a compact Chevalley basis $\fC_{0}$ of $\fG_{0}$, and each $\fG_{0}$ - module U admits a basis $\fB_{0}$ such that $\fC_{0}$ leaves invariant $\bq$ - span ($\fB_{0}$)([R1]).  If $\fL_{0} = \fC_{0}~ \cup~ \fB_{0}$, then $\fL_{0}$ is an $\br$ - basis of $\fN_{0}$ with structure constants in $\bq$.  In fact, the basis $\fB_{0}$ of U may be chosen so that $\fC_{0}$ leaves invariant $\bz$ - span ($\fB_{0}$), and in this case $\fL_{0}$ has structure constants in $\bz$.  We call $\fN_{0, \bq} = \bq$ - span($\fL_{0}$) a $\mathit{Chevalley~ rational~ structure}$ for $\fN_{0}$.  See sections 4 and 8 for details.

    Let exp : $\fN_{0} \rightarrow N_{0}$ denote the Lie group exponential map.  It is well known that exp is a diffeomorphism.  If $\Gamma_{0}$ is the subgroup of N$_{0}$ generated by exp($\bz$ - span ($\fL_{0}$)), then a result of Mal'cev ([M], [R2]) says that $\Gamma_{0}$ is a lattice in N$_{0}$; that is, $\Gamma_{0}$ is a discrete subgroup of N$_{0}$ and $\Gamma_{0} \backslash N_{0}$ is compact.

    We are particularly interested in developing methods for studying the geometry of $\Gamma_{0} \backslash N_{0}$.  Consider for example a $\mathit{rational}$ subalgebra $\fN_{0}'$ of $\fN_{0}$; that is, a subalgebra $\fN_{0}'$ spanned by vectors in $\fN_{0,\bq}$.  If N$_{0}' =$ exp($\fN_{0}'$), then it is well known (cf. Theorem 5.1.11 of [CG]) that $\Gamma_{0}' = \Gamma_{0} \cap N_{0}'$ is a lattice in N$_{0}'$.  We obtain an isometric immersion of $\Gamma_{0}' \backslash N_{0}'$ into  $\Gamma_{0} \backslash N_{0}$.  These compact, totally geodesic immersed submanifolds of $\Gamma_{0} \backslash N_{0}$ may be regarded as higher dimensional analogues of closed geodesics of $\Gamma_{0} \backslash N_{0}$.

    In section 10, as we explain in more detail below, we construct totally geodesic, rational subalgebras  $\fN_{0}'$ from the roots and weights of $\fG$ and V.

    This article sets up a general framework that should be useful for studying the geometry of N$_{0}$ and its quotient manifolds.  For example, the bracket operation in $\fN_{0} = U \oplus \fG_{0}$ and its complexification $\fN = V \oplus \fG$ can be nicely described in terms of the roots of $\fG = \fG_{0}^{\bc}$ and the weights of the $\fG$ - module V.  A Cartan subalgebra $\fH$ of $\fG$ decomposes V into a direct sum of weight spaces $\{V_{\mu} , \mu \in \Lambda \}$.  If $\mu$ and $\lambda$ are weights, then $[V_{\mu} , V_{\lambda}] = \fG_{\mu + \lambda}$, where $\fG_{\mu + \lambda} = \{0\}$ if $\mu + \lambda$ is not a root, and $\fG_{\mu + \lambda}$ is a root space of $\fG$ if $\mu + \lambda$ is a root.  If $\fH$ is adapted to $\fG_{0}$, then this simple bracket relation for $\fN = V \oplus \fG$  has a somewhat more complicated analogue for the real form $\fN_{0} = U \oplus \fG_{0}$.  The details may be found in sections 4 and 6 where we derive a " weight space " decomposition of U from the weight space decomposition of V = U$^{\bc}$.  The weight spaces $\{U_{\mu}, \mu \in \Lambda \}$ and the zero weight space U$_{0}$ are $\mathit{rational}$ with respect to the Chevalley rational structure $\fN_{0, \bq}$ for $\fN_{0}$.  Equivalently, we can find a basis of U that is a union of bases contained in $\fN_{0, \bq}$ of the weight spaces $\{U_{0}, U_{\mu} ; \mu \in \Lambda \}$.

    It is also useful to know the range of individual transformations ad u : $\fN_{0} \rightarrow \fG_{0}$, where u is an arbitrary element of U.  We obtain partial results in section 4 for the complex case $\fN = V \oplus \fG$ and in section 6 for the real case $\fN_{0} = U \oplus \fG_{0}$.  The results in the real case are relevant to the problem of computing the lengths of closed geodesics of the compact nilmanifold $\Gamma_{0} \backslash N_{0}$ defined above.  (See Proposition 4.2 of [Eb 3]).  The results in the real case are also useful for computing the Lie algebra of almost inner derivations of $\fN_{0}$, which in turn is important for determining isometry classes among compact, Laplace isospectral nilmanifolds ([GW]).  Furthermore, for compact 2-step nilmanifolds a knowledge of the almost inner derivations is also important for resolving the question of when two compact nilmanifolds with conjugate geodesic flows are isometric.  See [GM] for further details.  For a 2-step nilpotent Lie algebra $\fN_{0}$ with center $\fZ_{0}$  a derivation D : $\fN_{0} \rightarrow \fN_{0}$ is said to be $\mathit{almost~ inner}$ if $D(\zeta) \in ad \zeta(\fN_{0})$ for all $\zeta \in \fN_{0}$.

    The Lie algebras $\fN_{0} = U \oplus \fG_{0}$ have large groups of isometries that are also automorphisms.  In particular, the Lie algebra $\fG_{0}$ and the Weyl group W of $\fG = \fG_{0}^{\bc}$ determine natural subgroups of Aut$(\fN_{0})~ \cap$ I$(\fN_{0})$, where I$(\fN_{0})$ denotes the linear isometry group of $\fN_{0}$.  The identity component of  Aut$(\fN_{0})~ \cap$ I$(\fN_{0})$ has previously been described by J. Lauret in Theorem 3.12 of [L].  The finite subgroups of Aut$(\fN_{0})~ \cap$ I$(\fN_{0})$ determined by the Weyl group W appear to be previously unknown.  Section 6 contains the details, which depend on complex analogues for $\fN = V \oplus \fG$ established in sections 4 and 5.

    Let $\fH$ be a Cartan subalgebra of $\fG = \fG_{0}^{\bc}$ that is adapted to $\fG_{0}$, and let $\fH_{0}$ be a maximal abelian subalgebra of $\fG_{0}$ such that $\fH = \fH_{0}^{\bc}$.   Let $\fC = \{X_{\beta} : \beta \in \Phi ; \tau_{\alpha} : \alpha \in \Delta \}$ be a Chevalley basis for $\fG$ determined by $\fH$.  For each positive root $\beta$ we define A$_{\beta} = X_{\beta} - X_{- \beta}$ and  B$_{\beta} = i X_{\beta} + i X_{- \beta}$.  Then $\fC_{0} = \{A_{\beta}, B_{\beta} : \beta \in \Phi^{+} ; i\tau_{\alpha} : \alpha \in \Delta \}$ is the compact Chevalley basis of $\fG_{0}$ that was mentioned earlier.  The elements $\{i\tau_{\alpha} : \alpha \in \Delta \}$ are a basis for $\fH_{0}$ and leave invariant the weight spaces $\{U_{0}, U_{\mu} : \mu \in \Lambda \}$.  The behavior of $\{A_{\beta} , B_{\beta} : \beta \in \Phi^{+} \}$ as elements of End(U) encodes the basic information of U as a $\fG_{0}$ - module.  In section 7 we discuss basic information about the transformations $\{A_{\beta}, B_{\beta}\}$.

    Let $\fN_{0} = U \oplus \fG_{0}$ with complexification $\fN = V \oplus \fG$.  Fix a Cartan subalgebra $\fH$ that is adapted to $\fG_{0}$.  A weight $\lambda$ for V is said to be $\mathit{admissible}$ if there are no solutions to the equation 2$\lambda = k\alpha + \beta$ for any integer k and any roots $\alpha, \beta \in \Phi$.   In section 9 we show that almost every weight for any $\fG$ - module V is admissible and we classify, up to action of the Weyl group W, those weights that are not admissible.  Proofs of the classifications stated in section 9 are given in section 11, the appendix.

    We discuss two applications of the notion of admissible weight.

    1)  If $\lambda$ is an admissible weight, then we may determine the range of ad u$_{\lambda} : \fN_{0} \rightarrow \fG_{0}$ for any element u$_{\lambda}$ in the real weight space U$_{\lambda}$ defined by $\lambda$.  See Proposition 6.12.

    2)  In section 10 we show that admissible weights may be used to construct various types of rational, totally geodesic subalgebras of $\fG_{0}$.  To every weight $\lambda$ we may associate a complex weight vector H$_{\lambda}$ in the Cartan subalgebra $\fH$.  If $\fH = \fH_{0}^{\bc}$ is adapted to $\fG_{0}$, then we also obtain a real weight vector $\tilde{H}_{\lambda} \in \fH_{0} \subset \fG_{0}$.  Fix $\fH_{0} \subset \fG_{0}$ and $\fH = \fH_{0}^{\bc} \subset \fG$.

\indent \indent a)  If $\lambda$ is an admissible weight, then $\fN_{0}(\lambda) = U_{\lambda} \oplus \br \tilde{H}_{\lambda}$ is a rational, totally geodesic subalgebra of $\fN_{0}$ with 1- dimensional center.  In fact, $\lambda$ need only satisfy the weaker condition that 2$\lambda$ not be a root of $\fH$.  As a special case we may choose an arbitrary nonzero element u$_{\lambda}$ of U$_{\lambda}$.  Then $\br$ - span $\{u_{\lambda}, \tilde{H}_{\lambda}(u_{\lambda}), \tilde{H}_{\lambda} \}$ is a  3-dimensional, totally geodesic subalgebra of $\fN_{0}$, which is rational whenever u$_{\lambda}$ is a rational element of U$_{\lambda}$.  Any such 2-step nilpotent Lie algebra must be isomorphic to the 3-dimensional Heisenberg algebra.

\indent \indent b)  Let $\lambda$ be an admissible weight, and let $\beta$ be any root of $\fH$.  Let A$_{\beta} = X_{\beta} - X_{- \beta}$ and  B$_{\beta} = i X_{\beta} + i X_{- \beta}$ as above.  Let U$'_{\lambda,\beta} =$ \(\sum_{k=-\infty}^{\infty}U_{\lambda+k\beta} \) and let $\fG_{\lambda, \beta} =
\br$ - span $\{A_{\beta}, B_{\beta}, \tilde{H}_{\beta}, \tilde{H}_{\lambda} \}$.  Let  $\fN_{0}(\lambda,\beta) = U'_{\lambda,\beta} \oplus \fG_{\lambda, \beta}$.  Then $\fN_{0}(\lambda,\beta)$ is a rational, totally geodesic subalgebra of $\fN_{0}$ whose center $\fG_{\lambda, \beta}$ is 4-dimensional and Lie algebra isomorphic to the quaternions.

    We note that by the discussion of rational subalgebras above we obtain totally geodesic, isometric immersions into $\Gamma_{0} \backslash N_{0}$ of the compact nilmanifolds $\Gamma_{0}(\lambda) \backslash N_{0}(\lambda)$ and $\Gamma_{0}(\lambda, \beta) \backslash N_{0}(\lambda, \beta)$, where $ N_{0}(\lambda) = exp(\fN_{0}(\lambda)), \Gamma_{0}(\lambda) = \Gamma_{0} \cap  N_{0}(\lambda), N_{0}(\lambda, \beta) = exp(\fN_{0}(\lambda,\beta))$ and $\Gamma_{0}(\lambda, \beta) = \Gamma_{0} \cap N_{0}(\lambda, \beta)$.

\section {Notation and preliminaries}

    The following notation will be standard in this paper :

$\fG_0$ = a compact, semisimple real Lie algebra, a compact real form for $\fG$

B$_0$ = the Killing form of $\fG_0$ (negative definite)

$\fH_0$ = a maximal abelian subalgebra of $\fG_0$

J$_0$ = conjugation in $\fG$ = $\fG_0^\bc$ determined by $\fG_0$

$\fG$ = a complex, semisimple Lie algebra, the complexification of $\fG_{0}$

G = the connected Lie subgroup of SL(V) with Lie algebra $\fG$, V a complex $\fG$-module

B = the Killing form of $\fG$ (nondegenerate)

$\fH$ = a Cartan subalgebra of $\fG$

H$_\alpha$ = the root vector in $\fH$ determined by $\alpha~\epsilon~\Phi$,
B(H, H$_\alpha$) = $\alpha$(H) for all H $\epsilon~\fH$

$\tau_\alpha$ = 2H$_\alpha$ / B(H$_\alpha$, H$_\alpha$) \hspace{.1in} for $\alpha \in \Phi$

$\tilde{\tau}_{\alpha}$ = i $\tau_\alpha$ \hspace{.1in} for $\alpha \in \Phi$

$\langle \lambda,\alpha \rangle := \lambda(\tau_{\alpha})$ for $\lambda \in$ Hom($\fH,\bc$) and $\alpha \in \Phi$

$\fH_\br$ = $\br$ - span \{H$_\alpha : \alpha ~\epsilon ~\Phi$\} = $\br$ - span \{$\tau_\alpha : \alpha~
\epsilon~\Delta$\}
= $\{H ~\epsilon ~\fH$ : $\alpha$(H) ~$\epsilon ~\br$ for all $\alpha ~\epsilon ~\Phi\}$

$\Phi$ = the finite collection of roots in Hom ($\fH$, $\bc$) determined by $\fH$

$\Delta = \{\alpha_{1},~ ...~,~\alpha_{n}\}$, a base of simple roots for $\Phi$

$\Phi^+$ = the positive roots determined by $\Delta$

W = the Weyl group in Hom ($\fH$, $\bc$) determined by $\Phi$

U = a finite dimensional real $\fG_0$ module

V = a finite dimensional complex $\fG$ module

J = conjugation in V = U$^\bc$ determined by U

$\Lambda = \Lambda$(V) $=$ the finite collection of weights in Hom($\fH$, $\bc$) determined by V and $\fH$

$\Lambda^{+} = \Lambda^{+}$(V), V $=$ U$^{\bc}, =$ a subset of $\Lambda$ that contains exactly one of $\{\lambda, - \lambda \}$ for each $\lambda \in \Lambda$

$\tilde{\Lambda} =$ \{$\lambda \in$ Hom($\fH,\bc$) : $\lambda(\tau_{\alpha}) \in \bz$ for all $\alpha \in \Phi$\}, the abstract
weights determined by $\fH$.

V$_\lambda$ = \{v ~$\epsilon$ ~V : Hv = $\lambda$(H)v  for all H ~$\epsilon ~\fH$\}, the weight space
determined by $\lambda~ \epsilon~ \Lambda$

$\fG_\alpha$ = \{X ~$\epsilon ~\fG$ : [H,X] = $\alpha$(H)X for all H ~$\epsilon ~\fH$\}, the root space
determined by $\alpha~ \epsilon~ \Phi$

V = V$_0$ + \(\sum_{\lambda \epsilon \Lambda}\) V$_\lambda$, the complex weight space
decomposition determined by $\fH$

U = U$_0$ + \(\sum_{\lambda \epsilon \Lambda^{+}}\) V$_\lambda$, the real weight space
decomposition determined by $\fH_{0}$

$\fG = \fH$ + \(\sum_{\alpha \epsilon \Phi}$ $\fG_\alpha$, the root space decomposition of $\fG$
 determined by $\fH$

 $\fG_{0} = \fH_{0}$ + \(\sum_{\alpha \epsilon \Phi^{+}}$ $\fG_{0, \alpha}$, the root space decomposition of $\fG_{0}$ determined by $\fH_{0}$

 $\fN_{0} = U \oplus \fG_{0}$, a real 2-step nilpotent Lie algebra

 $\fN = V \oplus \fG$, the complexification of $\fN_{0}$.

 $\fC_{0}$ = a compact Chevalley basis.

 $\fB_{0}$ = a basis of U such that $\fC_{0}$ leaves invariant $\bq$ - span ($\fB_{0}$)

 $\fN_{0, \bq} = \bq$ - span ($\fB_{0}~ \cup~ \fC_{0}$), a Chevalley rational structure for $\fN_{0}$
 \newline

    We recall that B is nondegenerate on $\fH$ and positive definite on $\fH_\br$
(cf. [He, p.170]).  It is well known that $\beta(\tau_\alpha) ~\epsilon ~\bz$ for all $\alpha,\beta ~\epsilon
~\Phi$, and dim $\fG_\alpha$ = 1 for all $\alpha ~\epsilon ~\Phi$.
Typically, we will consider the situation that $\fG$ = $\fG_0^{\bc}$, $\fH$ = $\fH_0^{\bc}$ and V = U$^{\bc}$.

\section {Some 2-step nilpotent Lie algebras arising from representations}

    In this section we describe some complex and real 2-step nilpotent Lie algebras that arise from finite dimensional
 representations of a semisimple Lie algebra, either complex or real.
 \newline

$\mathit{Complex~ 2-step~ nilpotent~ Lie~ algebras}$

    Let $\fG$ be a complex, semisimple real Lie algebra, and let V be a finite dimensional,
 complex $\fG$ -module.  As a complex vector space we define $\fN$ = V$\oplus
\fG$.  Let B* denote a nondegenerate, symmetric, bilinear form on
 $\fN$ such that

\indent 1)  V and $\fG$ are orthogonal relative to B*.

\indent 2)  B*( Xv, v$'$ ) + B*( v, Xv$'$ ) = 0 for all v, v$'$ $\epsilon$ ~V and all X $\epsilon ~\fG$.

\indent 3)  B*( ad X(Y), Z ) + B*( Y, ad X(Z) ) = 0 for all X,Y, Z $\epsilon ~\fG$.

    We now let [ , ] denote the complex, 2-step nilpotent  Lie algebra structure on $\fN$ such that $\fG$ lies in the
center of $\fN$, [$\fN, \fN$] $\subseteq$ ~$\fG$ and B*( [v, v$'$], Z ) = B*( Zv, v$'$ ) for all v,v$'$ $\epsilon ~$V and
all Z $\epsilon ~\fG$.
\newline

$\mathit{Remarks}$

\indent 1)  A nondegenerate, symmetric bilinear form B* on a $\fG$-module V will be called $\mathit{\fG-invariant}$
if it satisfies property 2) above.  Property 3) may be restated by saying that
B* is $\fG$-invariant on $\fG$ relative to the adjoint representation.  If B* = B, the Killing form of $\fG$, then
B* satisfies 3).

\indent 2)  To construct an example satisfying 2) let $\fG_0$ be a compact real form of $\fG$, and let U be a finite
dimensional real $\fG_0$-module.   Let G$_0$ be the simply connected Lie group with Lie algebra $\fG_0$.
The group G$_0$ is compact since the Killing form B$_{0}$ on $\fG_0$ is negative definite.  See for example
[He, pp.133-134].  Let $\langle$ , $\rangle$ be a positive definite G$_0$-invariant inner product on U.  Then the
elements of $\fG_0$ are skew symmetric relative to $\langle$ , $\rangle$.  Let V = U$^\bc$ , which becomes a
complex $\fG$-module, and let B* denote
the complex extension of $\langle$ , $\rangle$ from U to V.  It is easy to see that B* is a $\fG$-invariant nondegenerate,
symmetric complex bilinear form on V.

    The next result says that the example above is the only way to find a nondegenerate, symmetric, $\fG$-invariant
 bilinear form B* on a complex $\fG$-module V.
 \newline

 \begin{proposition} Let V be a complex $\fG$-module, and let $\fG_0$ be a compact real form of $\fG$.  Suppose that V admits a nondegenerate, $\fG$-invariant, symmetric, complex bilinear form B*.  Then

\indent 1)  There exists a real $\fG_0$-module  U such that V = U$^\bc$ and  the restriction of B* to U is a positive   definite, real valued bilinear form.

\indent 2)  If J : V $\rightarrow$ V denotes the conjugation determined by U, then

\indent B*(Jv,w) $= \overline{B^{*}(v,Jw)}$ for all v,w $\in$ V.

\indent 3)  B*(v,Jv) is a positive real number for all nonzero v $\in$ V.
\end{proposition}

\begin{proof} 1)  We follow the proof of Proposition 6.4 of [B-tD].  Let G$_0$ be the compact simply connected Lie group with Lie algebra $\fG_0$.  It follows that V is a complex G$_0$-module since G$_0$ is simply connected.  Let $\ip$ be a fixed G$_0$-invariant Hermitian inner product on V, and define a conjugate linear isomorphism f : V $\rightarrow$ V by B*(v,w) = $\langle v, fw \rangle$ for v,w $\epsilon$~V.  Then f commutes with G$_0$, and $\langle f^{2}v,w \rangle = \langle f^{*}w,f^{*}v \rangle$ for all v,w $\in$ V.  It follows that f$^2$ is self adjoint and positive definite with respect to $\ip$.  On the eigenspace V$_\lambda$ for f$^2$, which is invariant under G$_0$, we define h = $\sqrt{\lambda}$ Id.  Then h extends to a complex linear isomorphism of V such that h commutes with G$_0$ , h$^2$ = f$^2$ and hf = fh.  Moreover, h is self adjoint and positive definite with respect to $\ip$ since it has these properties on the orthogonal subspaces $\{V_\lambda\}$.

    If J = hf$^{-1}$, then J : V $\rightarrow$ V is a conjugate linear map such that J$^2$ = Id and J commutes with
G$_0$.  If U is the +1 eigenspace of J, then U is invariant under G$_0$, and U$^\bc$ = V since iU is the $-$1 eigenspace
of J.  Note that h = f on U, and h(U) $\subseteq$ U since h commutes with J.

    We observe that B*(u,u$'$) = $\langle u, fu' \rangle$ = $\langle u, hu' \rangle$  for all
u, u$'$ $\epsilon$~U.  The restriction of B* to U has real values since = $\overline{\langle u, hu' \rangle}$ =
$\langle hu', u \rangle$ =$\langle u' , hu \rangle$ = B*(u$'$ , u) =
B*(u , u$'$) = $\langle u , hu' \rangle$ for all u , u$'$ ~$\epsilon$ ~U.  Since h is positive definite on V it follows that B* is positive definite on U.

    The proofs of 2) and 3) are straightforward, and we omit the details.
\end{proof}
    For later use we relate B* to the existence of $\fG_0$ - invariant Hermitian inner products on V.
\newline
\begin{proposition} Let V be a complex $\fG$-module, and let $\fG_0$ be a compact real form of $\fG$.  Suppose that V admits a nondegenerate, symmetric, complex bilinear form B* such that the elements of $\fG$ are skew symmetric relative to B*.  Let U be a real $\fG_0$-module such that V = U$^\bc$ and  the restriction of B* to U is a positive   definite, real valued bilinear form ( , ).  Let J : V $\rightarrow$ V denote the conjugation determined by U.  Define $\< , \>$ on V by $\<$v,w$\> =$ B*(v,Jw) for all v,w $\in$ V.  Then

\indent \indent a)  $\< , \>$ is a Hermitian inner product such that X*$ = -$ X for all X $\in \fG_0$, where X* denotes the metric adjoint of X.

\indent \indent b) $\langle$ Jv, Jv$' \rangle  = \overline{\langle \mathrm{v,v}' \rangle}$ for all v,v$' \in$ V.
\end{proposition}

\begin{proof} a) It is easy to see that $\< ,\>$ is conjugate bilinear on V since J is conjugate linear on V.  It now follows from 2) and 3) of (2.1) that $\< , \>$ is a Hermitian inner product on V.  Note that $\< , \> =$ B* $=$ ( , ) on U.  Hence $\< , \>$ is the unique Hermitian extension to V of the inner product ( , ) on U.  By condition 2) in the definition of B* we see that the elements of $\fG_0$ are skew symmetric on U relative to B*.  It follows that the elements of $\fG_0$ are skew Hermitian on V relative to $\< , \>$.

    b) Use the definition of $\< , \>$ and 2) of (2.1).
\end{proof}

$\mathit{Uniqueness~ of~ the~ complex~ Lie~ algebra~ structure}$

    The Lie algebra structure of $\fN$ = V$\oplus \fG$ seems to depend on the choice of a nondegenerate
symmetric bilinear form B* that satisfies the three properties listed at the beginning of this section.  However, the next result shows that the isomorphism type is uniquely determined.
\begin{proposition} Let B$_1$* and B$_2$* be two nondegenerate symmetric bilinear forms on $\fN$ = V$\oplus \fG$ that satisfy
properties 1), 2) and 3) above, and let [ , ]$_1$ and [ , ]$_2$ be the corresponding  2-step nilpotent Lie algebra
structures on $\fN$.  Then $\{\fN, [ , ]_1\}$ is isomorphic to $\{\fN, [ , ]_2\}$.
\end{proposition}

    We shall need a preliminary result.  If B* : V x V $\rightarrow \bc$ is a bilinear form on V, then we say that
T ~$\epsilon$ ~End(V) is $\mathit{symmetric}$ relative to B* if B*(Tx, y) = B*(x, Ty) for all x,y ~$\epsilon$~V.
\begin{lem}Let V be a finite dimensional complex vector space.  For every A ~$\epsilon$~GL(V) there exists
Y ~$\epsilon$~End(V) such that

\indent 1)  exp(Y) = A

\indent 2)  If X~$\epsilon$~End(V) and X commutes with A, then X commutes with Y.

\indent 3)  If B* is a bilinear form on V and A is symmetric with respect to B*, then Y is symmetric with respect to B*.
\end{lem}
\begin{cor}  For every A~$\epsilon$~GL(V) there exists T~$\epsilon$~GL(V) such that

\indent 1)  T$^2$ = A

\indent 2)  If X ~$\epsilon$~End(V) and X commutes with A, then X commutes with T.

\indent 3)  If B* is a bilinear form on V and A is symmetric with respect to B*, then T is symmetric with respect to B*.
\end{cor}

$\mathit{Proof ~of ~the ~lemma}$

    1) and 2)  From a  Jordan canonical form for A we may construct S and N in End(V) such that  A = S + N, SN = NS, S is semisimple and N is nilpotent.  Moreover, the transformations S and N are polynomials in A ([HK, Theorem 8, p.217]).

    Since A is nonsingular the transformation S is also nonsingular, and we may write A = SU, where
U = I + S$^{-1}$N.  If Z $=$ S$^{-1}$N, then Z is nilpotent since S$^{-1}$commutes with N.

    Since S is semisimple we may write V as a direct sum V$_1\oplus~...~\oplus~$V$_k$, where S =
$\lambda_i$ Id on V$_i$ for some distinct nonzero complex numbers
$\{\lambda_1, ... , \lambda_k\}$. Note that Z leaves each V$_i$
invariant since Z commutes with S.  From the Taylor series for
log(1+x) we define  $\lambda(Z) = \sum_{n=1}^{\infty} (-1)^{n} Z^{n}
/  n $ .  This sum is finite since Z is nilpotent, and it is known
that $exp(\lambda(Z)) = I + Z = U$ (cf. [He, Lemma 4.5, p. 270]).
The transformation  $\lambda$(Z) leaves each V$_i$ invariant since Z
does.  Moreover, A = SU leaves each V$_{i}$ invariant since U = I +
Z does.

    Choose $\mu_i~\epsilon~\bc$ so that  e$^{\mu_i}$ = $\lambda_i$ and define  Y~$\epsilon$~End(V) by
Y = $\mu_i$ Id + $\varphi$(Z) on each V$_i$.  On V$_i$ , exp(Y) = $\lambda_i$ exp($\varphi$(Z)) = SU = A,
and hence exp(Y) = A on V.

    Now let X~$\epsilon~$End(V) be an element that commutes with A.  Then X commutes with both S and N,
which are polynomials in A, and it follows that X commutes with Z =
S$^{-1}$N and $\varphi$(Z).  Note that X leaves invariant the
eigenspaces $\{V_i\}$ of S, and X commutes with Y = $\mu_i$ Id +
$\varphi$(Z) on V$_i$. Hence X commutes with Y on V, which proves
2).

\indent 3)  Let B* : V x V $\rightarrow \bc$ be a symmetric bilinear form, and let Sym(B*) =
$\{$T$~\epsilon$~End(V)
: B*(Tx, y) = B*(x, Ty) for all x,y~$\epsilon~$V\}.  It is easy to see that Sym(B*) is a subspace of End(V)
with the following additional properties :

\indent a)  $\lambda$ Id $\epsilon$~Sym(B*) for all $\lambda~\epsilon~\bc$.

\indent b)  If T~$\epsilon$~Sym(B*) is invertible, then T$^{-1}~\epsilon$~Sym(B*).

\indent c)  If T$_1$,T$_2~ \epsilon$~Sym(B*) and T$_1$,T$_2$ commute , then T$_1$T$_2~\epsilon$~Sym(B*).

\indent d)  If T~$\epsilon$~Sym(B*), then Sym(B*) contains any convergent power series in T.

    We use the notation in the proof of 1) and 2).  If A~$\epsilon$~Sym(B*), then from d) it follows
that S and N are in Sym(B*).  Note that Z$ = $S$^{-1}$N
$~\epsilon$~Sym(B*) by b) and c), and
$\varphi$(Z)$~\epsilon$~Sym(B*) by d). The restriction of Y to each
V$_i$ is $\mu_i$ Id + $\varphi$(Z), which is symmetric with respect
to B* on V$_i$ by a).  We conclude that Y~$\epsilon$~Sym(B*) since
S~$\epsilon$~Sym(B*) implies that B*(V$_i$, V$_j$) = $\{0\}$ if
i$\neq$j.
\newline

$\mathit{Proof ~of ~the ~corollary}$  Let Y~$\epsilon$~End(V) be chosen as in the lemma, and let T = exp(Y/2).
Clearly T$^2$ = A.  If X commutes with A, then X commutes with Y by 2) of the lemma, and hence X commutes with T.
 Finally, suppose that A~$\epsilon$~Sym(B*).  Then Y~$\epsilon$~Sym(B*) by 3) of the lemma, and
 T~$\epsilon$~Sym(B*) by d) in the proof of 3).
 \newline

$\mathit{Proof ~of ~the ~Proposition}$

    Let S : $\fN$ = V$\oplus \fG \rightarrow \fN$ be the invertible linear transformation such that
B$_2$*($\xi$,$\eta$) =  B$_1$*(S($\xi$),$\eta$) for all elements $\xi$,$\eta$ in $\fN$.  Note that
S is symmetric with respect to both  B$_1$* and  B$_2$* , and S leaves invariant both V and $\fG$.  Since
B$_1$* and B$_2$* are nondegenerate and $\fG$-invariant it follows from the definitions that for all
X~$\epsilon$~$\fG$ , S commutes with ad X on $\fG$ and with X on V.

    Note that S is a composition of transformations in End($\fN$) of the following two types :

\indent 1) S = Id on V and S is arbitrary on $\fG$.

\indent 2) S = Id on $\fG$ and S is arbitrary on V.

Hence it suffices to prove the Proposition for S in these two special cases.

\indent 1)  Define $\varphi  :  \{\fN, [ , ]_{1}\} \rightarrow
\{\fN, [ , ]_{2}\}$ by $\varphi$  = Id on V  and $\varphi  = S^{-1}$
on $\fG$.  We assert that $\varphi$  is a Lie algebra isomorphism.
For X,Y~$\epsilon$~V and Z~$\epsilon~\fG$ we have
B$_2$*($\varphi$([X,Y]$_1$), Z) = B$_2$*(S$^{-1}$([X,Y]$_1$),Z) =
B$_1$*([X,Y]$_1$,Z) = B$_1$*(ZX,Y) = B$_2$*(ZX,Y) =
B$_2$*([X,Y]$_2$, Z) = B$_2$*([$\varphi$X,$\varphi$Y]$_2$,Z).  Hence
$\varphi$([X,Y]$_1$) = [$\varphi$X,$\varphi$Y]$_2$ by the
nondegeneracy of B$_2$*.  The proof is complete since $\fG$ lies in
the center of $\fN$.

\indent 2)  By the corollary to the lemma above we may choose
T~$\epsilon$~GL(V) such that T$^2$ = S ,T commutes with $\fG$ on V
and T is symmetric with respect to both B$_1$* and B$_2$*.  Now
define $\varphi : \{\fN, [ , ]_{2}\} \rightarrow \{\fN, [ , ]_{1}\}$
by $\varphi$  = T on V  and $\varphi$ = Id on $\fG$.  We assert that
$\varphi$  is a Lie algebra isomorphism.  For X,Y~$\epsilon$~V and
Z~$\epsilon~\fG$ we have B$_2$*($\varphi$([X,Y]$_2$), Z) =
B$_2$*([X,Y]$_2$, Z) = B$_2$*(ZX, Y) = B$_1$*(SZX, Y) = B$_1$*(TZX,
TY) = B$_1$*(ZTX, TY) = B$_1$*([TX,TY]$_1$, Z) = B$_2$*([TX,TY]$_1$,
Z) = B$_2$*([$\varphi$X, $\varphi$Y]$_1$, Z). Hence
$\varphi$([X,Y]$_2$) = [$\varphi$X, $\varphi$Y]$_1$ by the
nondegeneracy of B$_2$*, and the proof is complete.
\newline

$\mathit{Real~ 2-step~ nilpotent~ Lie~ algebras}$

    We now define the real analogue of the complex 2-step nilpotent Lie algebra discussed above, and we prove in similar fashion the uniqueness of the Lie algebra structure, up to isomorphism.

    Let $\fG_0$ be a compact, semisimple real Lie algebra, and let U be a finite dimensional, real
$\fG_0$-module.  As a real vector space we define $\fN_{0}$ = U$\oplus\fG_0$.  Let $\ip$ denote a positive
definite inner product on $\fN$ such that

\indent 1)  U and $\fG_0$ are orthogonal relative to $\ip$.

\indent 2)  $\langle Xu, u' \rangle + \langle u, Xu' \rangle = 0$ for all u, u$' \epsilon$~U and all
X $\epsilon~\fG_0$.

\indent 3)  $\langle$ ad X(Y), Z $\rangle$ + $\langle$ Y, ad X(Z) $\rangle$
 = 0 for all X,Y, Z~$\epsilon~\fG_o$.

    We now let [ , ] denote the real 2-step nilpotent Lie algebra structure on $\fN_{0}$ such that $\fG_0$
 lies in the center of $\fN_{0}$, [$\fN_{0}$, $\fN_{0}$] $\subseteq \fG_0$, and $\langle [u, u'], Z \rangle =
 \langle Z(u), u' \rangle$ for all u,u$'~\epsilon$~U and all Z~$\epsilon~\fG_0$.
 \newline

$\mathit{Remarks}$

\indent 1)  If G$_0$ is the compact simply connected Lie group with Lie algebra $\fG_0$, then 2) holds  for any G$_0$-invariant inner product $\ip$ on U.  If U is an irreducible G$_0$-module, then $\ip$ is uniquely determined up to a constant positive multiple.  As above we define an inner product on $\fG_{0}$ to be $\fG_{0}- \mathit{invariant}$ if property 2) holds.

    To define an inner product $\ip$ on $\fG_0$ such that 3) holds we could set  $\ip = -$ B$_0$, where B$_0$ is the Killing form of $\fG_0$, or we could define $\langle$ X, Y $\rangle = -$ trace XY for all
X,Y ~$\epsilon~\fG_0$, where we use the fact that the elements of $\fG_0$ are skew symmetric relative to $\ip$.  If $\fG_0$ is simple, then the choice of $\ip$ on $\fG_{0}$ satisfying 3) is unique up to a constant positive multiple.

    We call an inner product $\ip$ on $\fG_0$ such that 3) holds an $\mathit{ad-\fG_{0}}~ invariant$ inner product.

\indent 2)  If U is a finite dimensional real $\fG_0$-module, then the elements of $\fG_0 \subseteq$ End(U) have eigenvalues in i$\br$ since the elements of $\fG_0$ are skew symmetric relative to a G$_0$-invariant inner product $\ip$ on U.

\indent 3)  It is not difficult to show that [$\fN_{0}$, $\fN_{0}$] = $\fG_0$.  Moreover, $\fG_0$ = $\fZ_0$, the center of $\fG_0$, $\Leftrightarrow$ Ker ($\fG_0$) = $\{$u~$\epsilon$~U : Z(u) = 0 for all Z~$\epsilon~\fG_0$\}
= $\{0\}$.  We omit the details.
\newline

$\mathit{Uniqueness~ of~ the~ real~ Lie~ algebra~ structure}$

\begin{proposition} Let $\langle~,~ \rangle_{1}$ and $\langle~,~ \rangle_{2}$ be two positive definite inner products on $\fN_{0} = U \oplus \fG_0$ that satisfy properties 1), 2) and 3) above, and let [ , ]$_1$ and [ , ]$_2$ be the corresponding  2-step nilpotent Lie algebra structures on $\fN$.  Then
\{$\fN_{0}, [ , ]_{1}$\} is isomorphic to \{$\fN_{0}, [ , ]_{2}$\}.
\end{proposition}

\begin{proof}  We proceed as in (2.2) but the proof is slightly easier.  Let S : $\fN_{0} = U\oplus\fG_0 \rightarrow \fN_{0}$ be the invertible linear transformation such that $\langle \xi, \eta \rangle_2 = \langle S\xi, \eta \rangle_1$ for all $\xi , \eta$ in $\fN_{0}$.  Note that S is symmetric and positive definite with respect to both $\ip_1$ and $\ip_2$.  Since $\ip_1$ and $\ip_2$ are nondegenerate and $\fG_0$-invariant it follows that for all X~$\epsilon~\fG_0$, S commutes with ad X on $\fG_0$ and with X on U.

    Since S is symmetric we may write U as a direct sum U$_1\oplus~... ~\oplus~$U$_k$, where S =
$\lambda_i$ Id on V$_i$ for some distinct positive real numbers $\{\lambda_1,~...~,\lambda_k\}$.  Note that
$\fG_0$ leaves each U$_i$ invariant since $\fG_0$ commutes with S.  Define T on U by setting T = $\sqrt{\lambda_i}$ Id on each U$_i$.  It is clear that T$^2$ = S, T commutes with $\fG_0$ and T is positive definite and symmetric with respect to both $\ip_1$ and $\ip_2$.  The proof of the proposition is now completed exactly as in the proof of (2.3).
\end{proof}

$\mathit{Commutativity~ of~ irreducible~ submodules}$

\begin{proposition}  Let $\fN_{0} = U \oplus \fG_0$, and let U = U$_1 \oplus$...$\oplus$U$_N$ be an orthogonal direct sum decomposition of U into irreducible $\fG_0$-submodules $\{U_i\}$.  Then
[U$_i$, U$_j$] = $\{0\}$ if i$\neq$j.
\end{proposition}
\begin{proof} If Z is any element of $\fG_0$ and
u$_i$,u$_j$ are any elements of U$_i$ ,U$_j$ with i$\neq$j, then $\langle$ [u$_i$ , u$_j$], Z $\rangle$ =
$\langle$ Z(u$_i$),u$_j \rangle$ = 0 since U$_i$ is orthogonal to U$_j$.  It follows that
[U$_i$,U$_j$] = $\{0\}$.
\end{proof}

\section{Chevalley bases and the universal enveloping algebra}

    Let $\fG$ be a complex semisimple Lie algebra, and let $\fH$ be a Cartan subalgebra of $\fG$.  Let $\Phi \subset$ Hom($\fH, \bc$) denote the roots determined by $\fH$.  Let $\Delta$ be a base of simple
positive roots for $\Phi$, and let $\Phi^{+}$ denote the positive roots.  Let \rs ~denote the root space decomposition determined by $\fH$.  The root spaces $\fG_{\alpha}$ are 1-dimensional for all $\alpha \in \fH$.

    If B is the Killing form on $\fG$, then B is nondegenerate on $\fH$.  We define $\tau_{\alpha} =
2H_{\alpha}/B(H_{\alpha},H_{\alpha})$, where H$_{\alpha} \in \fH$ is the unique vector such that $\alpha(H) = B(H,H_{\alpha})$ for all H
$\in \fH$.

    Let X$_{\alpha} \in \fG_{\alpha}~,~\alpha \in \Phi$  be elements such that

\indent 1)  [X$_{\alpha}$, X$_{-\alpha}$] = $\tau_{\alpha}$ for all $\alpha \in \Phi$.

\indent 2)  If $\alpha~ , ~\beta~, ~ \alpha + \beta \in \Phi$, then [X$_{\alpha}$, X$_{\beta}$] = C$_{\alpha,\beta}$ X$_{\alpha+\beta}$,
where C$_{\alpha,\beta} = -$ C$_{-\alpha,-\beta}$.

A basis $\fC$ = \{X$_{\beta}$ : $\beta \in \Phi$ ; $\tau_{\alpha}$ : $\alpha \in \Delta$ \} is called a $\mathit{Chevalley~ basis}$
of $\fG$ (determined by $\fH$) if the elements X$_{\beta} \in \fG_{\beta}$, $\beta \in \Phi$, satisfy 1) and 2) above.  In this case
the structure constants of $\fC$ lie in $\bz$ .  See [Hu, pp.144-145] for a more precise statement about the structure constants and the nonuniqueness of Chevalley bases.
\newline

$\mathit{Chevalley~ bases~ and~ \fG-Automorphisms~ of~ order~ 2}$

    The existence of Chevalley bases is closely connected with the existence of automorphisms $\varphi$ of $\fG$ of order two such that
 $\varphi \equiv -$ Id on a Cartan subalgebra $\fH$.  This is explained in part by the next two results.

\begin{proposition} Let $\fG$ be a complex semisimple Lie algebra, and let $\fH$ be a Cartan subalgebra of $\fG$.  Let $\varphi$
be an automorphism of order 2 such that $\varphi \equiv -$ Id on $\fH$.  Let X$_{\alpha}$ be the unique element of $\fG_{\alpha}$ such that
[X$_{\alpha}$ , $\varphi$(X$_{\alpha}$)] $= -~ \tau_{\alpha}$.  Then
$\fC$ = \{X$_{\beta}$ : $\beta \in \Phi$ ; $\tau_{\alpha}$ : $\alpha \in \Delta$ \} is a Chevalley basis of $\fG$.
\end{proposition}

\begin{proof}  See [Hu : p.145].
\end{proof}

    We note that if $\varphi \equiv -$ Id on $\fH$, then  $\varphi$($\fG_{\alpha}$) = $\fG_{-\alpha}$ for all $\alpha \in \Phi$, and it
follows that [X$_{\alpha}$ , $\varphi$(X$_{\alpha}$)] is a multiple of $\tau_{\alpha}$ for any choice of X$_{\alpha}$ in $\fG_{\alpha}$.

\begin{proposition}  Let $\fG$ be a complex semisimple Lie algebra, and let $\fH$ be a Cartan subalgebra of $\fG$.  Let
$\fC$ = \{X$_{\beta}$ : $\beta \in \Phi$ ; $\tau_{\alpha}$ : $\alpha \in \Delta$ \} be a Chevalley basis of $\fG$.
For $\alpha \in \Delta$ let $\tilde{\tau}_{\alpha}$ = i $\tau_{\alpha}$, and for $\beta \in \Phi^{+}$ let A$_{\beta}$ =
X$_{\beta}~ -~ $X$_{-\beta}$ and B$_{\beta}$ = i X$_{\beta}$ + i X$_{-\beta}$.  Let
$\fC_0 =$ \{$\tilde{\tau}_{\alpha}$ : $\alpha \in \Delta$ ; A$_{\beta}$, B$_{\beta}$ : $\beta \in \Phi^{+}$ \}.  Then

\indent (1) $\fG_0 = \br$-span ($\fC_0$)  is a real, compact, semisimple Lie algebra such that $\fG_0^{\bc} = \fG$.

\indent (2) $\fH_0 = \br$-span \{$\tilde{\tau}_{\alpha}$: $\alpha \in \Delta$ \} is a maximal abelian subalgebra of $\fG_0$  and
$\fH_0^{\bc} = \fH$.

\indent (3)   J$_0$(X$_{\beta}$) = $-$ X$_{-\beta}$ for all $\beta \in \Phi$ , where J$_0$  denotes the conjugation in $\fG$ determined by $\fG_0$.

\indent (4)  If $\varphi$~ : $\fG \rightarrow \fG$ is the
$\bc$-linear map such that $\varphi \equiv -$ Id  on $\fH$ and
$\varphi$(X$_{\beta}$) = J$_0$(X$_{\beta}$) for all $\beta \in \Phi$
, then $\varphi$  is an automorphism of $\fG$ such that
$\varphi^{2}$ = Id.
\end{proposition}

$\mathit{Remarks}$

\indent 1)  Property (3) is equivalent to the condition (3$'$) [X$_{\beta}$, J$_0$(X$_{\beta}$)] = $-~ \tau_{\beta}$ for all $\beta \in \Phi$ since J$_0$(X$_{\beta}$) $\in \fG_{\beta}$ by (6.2) below and $\fG_{\beta}$ is 1-dimensional for all $\beta \in \Phi$.

\indent 2)  The automorphisms of $\fG$ act transitively on both the Cartan subalgebras $\fH$ and the compact real forms $\fG_0$ of $\fG$.  It follows that  every compact real form $\fG_0$ of $\fG$ can be obtained as in (1) of (3.2) for an appropriate choice of Cartan subalgebra $\fH$.
\newline

$\mathit{Proof~ of~ Proposition~ 3.2}$.  Assertion (1) is well known, and a proof may be found, for example, in [He, pp. 181-182].  We prove (2).
We note that by (1) $\fH_0$ is clearly an abelian subalgebra of $\fG_0$ such that $\fH = \fH_0^{\bc}$.  If $\fH_0$ were contained in some abelian subalgebra $\fH'_0$ of larger dimension, then $\fH'$ $= \fH'_{0}$$^{\bc}$ would be a Cartan subalgebra of $\fG$ that properly contains $\fH$.  This is impossible, which completes the proof of (2).

\indent (3)  By definition A$_{\beta}$ = X$_{\beta}~ -$ X$_{-\beta}$ and B$_{\beta}$ = i X$_{\beta}$ + i X$_{-\beta} \in \fG_{0}$
for all $\beta \in \Phi$.  It follows immediately that 2X$_{\beta}$ = A$_{\beta}~ -~$i B$_{\beta}$ and $-$2X$_{-\beta}$ =
A$_{\beta}~ +~$i B$_{\beta}$.  Hence J$_0$(X$_{\beta}$) = $-$X$_{-\beta}$.

\indent (4)  Let $\varphi$~ : $\fG \rightarrow \fG$ be the
$\bc$-linear map such that $\varphi \equiv -$ Id on $\fH$ and
$\varphi$(X$_{\beta}$)
 = J$_0$(X$_{\beta}$) for all $\beta \in \Phi$.  Clearly $\varphi^{2}$ = Id.  To show that $\varphi$  is an automorphism of $\fG$ it
 suffices to show that $\varphi$([H, X$_{\beta}$]) = [$\varphi$(H), $\varphi$(X$_{\beta}$)] for all H $\in \fH, \beta \in \Phi$ and
 $\varphi$([X$_{\alpha}$, X$_{\beta}$]) = [$\varphi$(X$_{\alpha}$), $\varphi$(X$_{\beta}$)] for all $\alpha~,~\beta \in \Phi$.
These equalities follow immediately from 3), the definition of
$\varphi$ , the definition of a Chevalley basis and the fact that
J$_{0}(\fG_{\beta}) = \fG_{- \beta}$.
\newline

$\mathit{Universal~ enveloping~ algebra~ \fU(\fG)~ and~ the~ subring~ \fU(\fG)_{\bz}}$

    Let $\fG$ be a complex semisimple Lie algebra.  The universal enveloping algebra $\fU(\fG)$ may be regarded as the set of finite $\bc$-linear combinations of formal expressions X$_{1}$ X$_{2}$ ... X$_{n}$, where n is any positive integer and the X$_{i}$ are arbitrary elements of $\fG$.  Multiplication and addition are defined in the obvious way.  The only relations in $\fU(\fG)$ are those that come from the bracket in $\fG$ : XY $-$ YX = [X,Y] for elements X,Y of $\fG$.  For a complex finite dimensional vector space V every Lie algebra homomorphism $\rho : \fG \rightarrow$ End(V) extends to an algebra homomorphism $\rho : \fU(\fG) \rightarrow$ End(V).  Moreover, $\rho(\fU(\fG))$ is the subalgebra of End(V) generated by $\rho(\fG)$.

    Let $\fH$ be a Cartan subalgebra of $\fG$ with roots $\Phi$ and base $\Delta$ for $\Phi$.  Let $\fC$ = \{X$_{\beta}$ : $\beta \in \Phi$ ; $\tau_{\alpha}$ : $\alpha \in \Delta$ \} be a Chevalley basis for $\fG$.  Let $\fU(\fG)_{\bz}$ denote the subring of $\fU(\fG)$ generated by 1 and the elements \{(1/n!)(X$_{\alpha})^{n} : \mathrm{n~ is~ any~ positive~ integer~ and~} \alpha \in \Phi$ \}.
\newline

$\mathit{Remarks}$

\indent 1)  We note that each element X$_{\alpha}$ from the Chevalley basis $\fC$ is nilpotent as an element of End(V) since ad  X$_{\alpha}$ is a nilpotent element of End($\fG$).  In particular the power series exp(X$_{\alpha}$) is a finite sum and belongs not only to GL(V) but to $\fU(\fG)_{\bz}$ regarded as a subring of $\fU(\fG) \subset$ End(V).  This observation will be used in section 5 in the discussion of the Weyl group.

\indent 2)  Note that $\fC \subset \fU(\fG)_{\bz}$.  To verify this, observe that the vectors \{X$_{\alpha} :     \alpha \in \Phi$ \} line in $\fU(\fG)_{\bz}$ by the definition of $\fU(\fG)_{\bz}$.  The vectors \{$\tau_{\beta} : \beta \in \Delta$ \} lie in $\fU(\fG)_{\bz}$ since $\tau_{\beta} = [X_{\beta}, X_{-\beta}] = X_{\beta} \ X_{-\beta} - X_{-\beta} \ X_{\beta}$.

\indent 3)  The ring $\fU(\fG)_{\bz}$ is an infinitely generated $\bz$-module with a particularly nice basis due to Kostant.  See [Hu, pp. 149-154] for details.

    For later convenience we introduce the following.
\newline

$\mathit{Chevalley~ bases~ of~ compact~ real~ forms~ \fG_0}$

    A $\mathit{compact~ Chevalley~ basis}$ of a compact, semisimple real Lie algebra $\fG_0$ is a basis $\fC_0 =$ \{$\tilde{\tau}_{\alpha}$: $\alpha \in \Delta$ ; A$_{\beta}$, B$_{\beta}$ : $\beta \in \Phi^{+}$\} for $\fG_0$ as defined in (3.2).

    Since the structure constants of a Chevalley basis $\fC$ of $\fG = \fG_0^{\bc}$ are integers we immediately obtain the following

\begin{proposition}  The structure constants of the basis $\fC_0$ of $\fG_0$ are integers.
\end{proposition}

    We recall that an inner product $\ip$ on $\fG_{0}$ is $ad-\fG_{0}~ invariant$ if ad X : $\fG_{0} \rightarrow \fG_{0}$ is skew symmetric relative to $\ip$ for every X $\in \fG_{0}$.

\begin{proposition}  Let $\fG_{0}$ be a compact, semisimple Lie algebra over $\br$, and let  $\fC_{0}$ be a compact Chevalley basis of $\fG_{0}$.  Then  $\fC_{0}$ is an orthogonal basis of $\fG_{0}$ for every $ad - \fG_{0}$ invariant inner product $\ip$ on $\fG_{0}$.
\end{proposition}

$\mathbf{Lemma~ 1}$  Let $\fG_{0} = \fG_{1} \oplus ~...~ \oplus \fG_{N}$ be a decomposition of $\fG_{0}$ into a direct sum of its simple ideal $\fG_{1}, ~...~, \fG_{N}$.  If $\ip$ is any $ad- \fG_{0}$ invariant inner product on $\fG_{0}$, then the ideals $\{\fG_{i}\}$ are orthogonal.

\begin{proof}  Fix distinct integers i,j.  Let $\fH_{i} = \{X \in \fG_{i} : \langle X , \fG_{j} \rangle = 0 \}$.  It follows that $\fH_{i}$ is an ideal of $\fG_{i}$ since $\ip$ is $ad-\fG_{0}$ invariant.  The $ad-\fG_{0}$ invariance of $\ip$ also shows that  $\fH_{i}$ contains $ad X(\fG_{i})$ for all X $\in \fG_{i}$.  The subspace  $ad X(\fG_{i})$ is nonzero for all nonzero X $\in \fG_{i}$ since $\fG_{i}$ has trivial center.  We conclude that $\fH_{i} = \fG_{i}$ since $\fG_{i}$ is simple.
\end{proof}

$\mathbf{Lemma~ 2}$ Let $\fG_{0}$ be a compact, simple Lie algebra over $\br$.  Then every $ad-\fG_{0}$ invariant inner product $\ip$ on $\fG_{0}$ is a negative multiple of the Killing form B$_{0}$.

\begin{proof} We recall that $-B_{0}$ is an $ad-\fG_{0}$ invariant inner product on $\fG_{0}$.  Now let $\ip$ be an arbitrary $ad-\fG_{0}$ invariant inner product on $\fG_{0}$.  Write $\langle X, Y \rangle = - B_{0}(SX,Y)$ for all X,Y $\in \fG_{0}$, where S : $\fG_{0} \rightarrow \fG_{0}$ is a linear transformation that is symmetric with respect to both $\ip$ and $-B_{0}$.  Since ad X : $\fG_{0} \rightarrow \fG_{0}$ is skew symmetric with respect to both $\ip$ and $-B_{0}$ it follows that ad X $\circ~ S = S~ \circ$ ad X for all X $\in \fG_{0}$.  Hence any eigenspace of S in $\fG_{0}$ is an ideal of $\fG_{0}$.  Since $\fG_{0}$ is simple and both $\ip$ and $-B_{0}$ are positive definite it follows that S is a positive multiple of the identity.
\end{proof}

$\mathit{Proof~ of~ 3.4 }$  Every compact Chevalley basis of $\fG_{0}$ is a union of compact Chevalley bases from its simple ideals $\{\fG_{i}\}$.  Hence by Lemma 1 it suffices to consider the case that $\fG_{0}$ is simple.  By Lemma 2 it suffices to consider the $ad-\fG$ invariant inner product $-B_{0}$ on $\fG_{0}$.

    We note that $\fG = \fG_{0}^{\bc}$ is a complex semisimple Lie algebra.  If B,B$_{0}$ denote the Killing forms for $\fG, \fG_{0}$ respectively, then it is routine to check that B(X,Y) = B$_{0}$(X,Y) for all X,Y $\in \fG_{0}$.  It now follows immediately from (4.5) below and the definition of a compact Chevalley basis $\fG_{0}$ that $\fG_{0}$ is orthogonal relative to B and hence to B$_{0}$. $\square$

\section{Weight space decomposition of a complex $\fG$-module}

    We collect some useful facts.  In this section we fix a Cartan subalgebra $\fH$ of $\fG$, and a complex $\fG$ - module V.  We let $\Lambda \subset$ Hom($\fH$,$\fG$) denote the set of weights determined by V, and we let \ws be the weight space decomposition determined by $\fH$.  We let $\Delta$ denote a base for the roots $\Phi$ of $\fH$, and let $\fC$ = \{X$_{\beta}$ : $\beta \in \Phi$ ; $\tau_{\alpha}$ : $\alpha \in \Delta$ \} be a Chevalley basis of $\fG$ as defined above in section 3.  Let $\fC_0 =$ \{$\tilde{\tau}_{\alpha}$: $\alpha \in \Delta$ ; A$_{\beta}$, B$_{\beta}$ : $\beta \in \Phi^{+}$\} denote the corresponding compact Chevalley basis for $\fG_0$ as defined in (3.2).
 \newline

 $\mathit{Some~ relations~ between~ roots~ and~ weights}$

 \begin{proposition}  Let $\lambda \in \Lambda$ and $\beta \in \Phi$. Then V$_{\lambda + \beta} \neq 0
\iff \fG_{\beta}$(V$_{\lambda}) \neq \{0\}$.
\end{proposition}

\begin{proof}  Since $\fG_{\beta}(V_{\lambda}) \subseteq V_{\lambda+\beta}$ it is obvious that  V$_{\lambda + \beta} \neq 0$  if $\fG_{\beta}(V_{\lambda}) \neq \{0\}$.  Suppose now that
V$_{\lambda + \beta} \neq 0$.  Since $\lambda$ is a weight we know that $\lambda(\tau_{\beta})$ is an integer.   Recall that $\beta$($\tau_{\beta}$) = 2 for all $\beta \in \Phi$.

    Choose elements X$_{\beta} \in \fG_{\beta}$ and X$_{-\beta} \in \fG_{-\beta}$ such that
$[X_{\beta}, X_{-\beta}] = ~\tau_{\beta}$.  Then $ [\tau_{\beta}, X_{\beta}] = 2X_{\beta}$ and $[\tau_{\beta}, X_{-\beta}] = -2X_{-\beta}$.  If we define $\fs \fl(2,\beta) = \br$-span $\{X_{\beta}, X_{-\beta}, \tau_{\beta}\}$, then clearly $\fs \fl(2,\beta)$ is a subalgebra of $\fG$ isomorphic to $\fs \fl(2,\bc)$ for each $\beta \in \Phi$.

    We need a preliminary result that is known but not so easy to find in the literature.  We include a proof for completeness.

\begin{lem}  For $\lambda \in \Lambda$ and $\beta \in \Phi$ we let k $\geq 0$ and j $\geq 0$ be the largest integers such that $V_{\lambda + k\beta} \neq \{0\}$ and  $V_{\lambda - j\beta} \neq \{0\}$.  Then

\indent a)  $\lambda(\tau_{\beta}) = j - k$.

\indent b)  V$_{\lambda+r\beta} \neq$ \{0\} if $ - j \leq r \leq k$.
\end{lem}

$\mathit{Proof~ of~ the~ Lemma}$  Let V$'_{\lambda,\beta} =$ \(\sum_{r=-\infty}^{\infty} V_{\lambda+r\beta} \).  Clearly V$'_{\lambda,\beta}$ is an
$\fs \fl(2,\beta)$-module.  The representation theory of $\fs \fl(2,\bc)$ states that there is a positive integer N such
that N is the largest eigenvalue and $-$ N the smallest eigenvalue for $\tau_{\beta}$ on V$'_{\lambda,\beta}$.  Moreover, the
eigenvalues of $\tau_{\beta}$ on V$'_{\lambda,\beta}$ contain \{$N-2p : 0 \leq p \leq N$ \}.  In particular $N
= (\lambda+k\beta)(\tau_{\beta}) = \lambda(\tau_{\beta}) + 2k$ and $-N = (\lambda-j\beta)(\tau_{\beta}) = \lambda(\tau_{\beta}) - 2j$.
The proof of a) follows immediately.

    If $-j \leq r \leq k$, then $(\lambda+r\beta)(\tau_{\beta}) = N - 2p$, where $p = k - r \geq 0$.  Note that
$p \leq N$ since $N - 2p = - N +2(r+j) \geq - N$.  By the discussion of the previous paragraph $(\lambda+r\beta)(\tau_{\beta})$
is an eigenvalue of $\tau_{\beta}$ on V$'_{\lambda,\beta}$ and this eigenvalue can only occur on V$_{\lambda+r\beta}$.  This proves b) and
completes the proof of the Lemma.

    We now complete the proof of the proposition.  Choose a largest integer k $\geq$ 1 such that V$_{\lambda + k\beta} \neq \{0\}$.  If N = $\lambda$($\tau_{\beta}$) + 2k, then $\tau_{\beta}$ = N Id on V$_{\lambda+k\beta}$.  We let v = v$_0$ be a nonzero element of V$_{\lambda+k\beta}$.  For
every integer r $\geq$ 1 we define inductively v$_r$ = (1/r)X$_{-\beta}$(v$_{r-1}$) or equivalently, v$_r$ =  (1/r!) (X$_{-\beta}$)$^r$v$_0$.  Observe that on V we have X$_{\beta}$ X$_{-\beta}$ = X$_{-\beta}$ X$_{\beta}$ + [X$_{\beta}$, X$_{-\beta}$] = X$_{-\beta}$ X$_{\beta}$ + $\tau_{\beta}$.  By induction on r we obtain in standard fashion
(cf. section 7.2 of [Hu])

\indent a)  For r $\geq 0$

\indent \indent     i)  v$_r \in V_{\lambda+(k-r)\beta}$
 \hspace {.5in} \imp $\tau_{\beta}(v_r) = (N - 2r)v_r$.

\indent \indent         ii)  X$_{-\beta}$(v$_r$) = ($r+1$)v$_{r+1}$.

\indent \indent         iii)  X$_{\beta}$(v$_r$) = (N$-r+1$)v$_{r-1}$.

    If $\{v_{0}, v_{1},~ ...~, v_{r}\}$ are all nonzero, then they are linearly independent since they belong to different
eigenspaces of $\tau_{\beta}$ by i).  Let m be the positive integer such that v$_{m} \neq 0$ but v$_{m+1} = 0$.  By iii) we have
$0 = X_{\beta}(v_{m+1}) = (N-m)v_{m}$, which implies that N = m. This proves

\indent b)  v$_N \neq$ 0 but v$_{N+1}$ = 0.  If W = $\bc$-span
$\{v_0, v_1, ... , v_N\}$, then W is a $\fs \fl(2,\bc)$-module
of dimension N+1.  The eigenvalues of $\tau_{\beta}$ on W are $\{N-2r : 0 \leq r \leq N\}$.

    If j $\geq 0$ is the largest integer such that V$_{\lambda - j\beta} \neq \{0\}$, then it follows from a) of the Lemma above that N $= \lambda(\tau_{\beta}) + 2k = k + j \geq k$.  Hence $ 0 \neq v_k \in (X_{-\beta})^k(V_{\lambda+k\beta}) \subseteq V_{\lambda}$ and $0 \neq(N-k+1) v_{k-1} = X_{\beta}(v_k) \in X_{\beta}(V_{\lambda})$.  This completes the proof of the Proposition.
\end{proof}

\begin{proposition}  Let V be a complex $\fG$ - module.  If $\Lambda \cap \Phi$ is nonempty, then the zero weight space V$_{0}$ is nonzero.  If V$_{0} \neq \{0\}$ and  V$_{0}$ is not a submodule of V, then  $\Lambda \cap \Phi$ is nonempty.
\end{proposition}

$\mathit{Remark}$  If V is irreducible and V$_{0} \neq \{0\}$, then with only a few exceptions $\Phi \subset \Lambda$.  Hence if $\Phi \cap \Lambda$ is nonempty, then $\Phi \subset \Lambda$ except in a few cases.  See [D] for a list of the exceptions.

\begin{proof}  Suppose first that V$_{0} \neq \{0\}$ and  V$_{0}$ is not a submodule of V.  Then $\fG_{\beta}(V_{0}) \neq \{0\}$ for some $\beta \in \Phi$.  This implies that $\beta \in \Lambda \cap \Phi$ since $\fG_{\beta}(V_{0}) \subseteq V_{\beta}$.

    Conversely, suppose that $\beta \in \Lambda~ \cap~ \Phi$.  Let k, j be the largest nonnegative integers such that V$_{\beta + k\beta} \neq \{0\}$ and V$_{\beta - j\beta} \neq \{0\}$.  By a) of the lemma in (4.1) we have j $-$ k $= 2$.  Since j $\geq 2$ we have $ - j \leq - 1 \leq k$.  Hence V$_{0} \neq \{0\}$ by b) of the Lemma in (4.1).
\end{proof}

$\mathit{\fG_{0}-invariant~ Hermitian~ inner~ products}$

    Let $\fG$ be a finite dimensional, complex, semisimple Lie algebra, and let $\fG_{0}$ be a compact real form for $\fG$.  Let V be a finite dimensional complex $\fG$-module.  We say that a Hermitian inner product $\< , \>$ on V is $\mathit{\fG_{0}-invariant}$ if X* $= -$ X for all X $\in \fG_{0}$, where X* denotes the metric transpose of X.  We saw already in (2.2) how such Hermitian inner products may arise when V = U$^{\bc}$.  In fact, any complex $\fG$-module V admits a $\fG_{0}$-invariant Hermitian inner product $\< , \>$.  Let G$_{0}$ denote the compact, simply connected Lie group with Lie algebra $\fG_{0}$.  If one starts with any Hermitian inner product on V, then by averaging over G$_{0}$ one obtains a Hermitian inner product $\< , \>$ preserved by the elements of G$_{0}$.  It follows immediately that $\< , \>$ is $\fG_{0}$-invariant.

\begin{proposition}  Let V be a finite dimensional, complex $\fG$-module, and let $\< , \>$ be a $\fG_{0}$-invariant Hermitian inner product on V.  Let G be the connected Lie subgroup of SL(V) with Lie algebra $\fG$.  Let $\fC$ = \{X$_{\beta}$ : $\beta \in \Phi$ ; $\tau_{\alpha}$ : $\alpha \in \Delta$ \} be a Chevalley basis of $\fG$.  Then

\indent 1)  G is self adjoint (invariant under the metric transpose *).   If G$_{0}$ is the connected Lie subgroup of G whose Lie algebra is $\fG_{0}$, then G$_{0}$ is the identity component of the subgroup of unitary elements of G.

\indent 2)  X$_{\beta}$* $=$ X$_{-\beta}$ for all $\beta \in \Phi$.

\indent 3)  H* = H for all H $\in \fH_{\br}$.

\indent 4)  The weight spaces \{V$_{0}$,V$_{\mu} : \mu \in \Lambda$ \} of V are orthogonal relative to $\< , \>$.
\end{proposition}

\begin{proof} 1)  If  $\< , \>$ is a  $\fG_{0}$-invariant Hermitian inner product on V, then $\fG$ is self adjoint since X* $= -$ X for all X $\in \fG_{0}$ and $\fG_{0}^{\bc} = \fG$.  It follows that the connected Lie group G $\subset$ SL(V) is self adjoint since G is generated by exp($\fG$).  Moreover, $\fG_{0}$ is the subalgebra of skew Hermitian elements of $\fG$.  This proves the remaining assertion in 1).

    2)  Recall that $\fG_{0} = \br$-span($\fC_{0}$), where $\fC_{0}$ is the compact Chevalley basis of $\fG_{0}$ given by $\fC_0 =$ \{$\tilde{\tau}_{\alpha}$: $\alpha \in \Delta$ ; A$_{\beta}$, B$_{\beta}$ : $\beta \in \Phi^{+}$\}.  It suffices to prove 1) for all $\beta \in \Phi^{+}$ since the metric transpose operation * is an involution. Now A$_{\beta} = X_{\beta} - X_{-\beta}$ and B$_{\beta} = i X_{\beta} + iX_{-\beta}$ for all $\beta \in \Phi^{+}$.  Hence X$_{\beta} = \frac{1}{2}(A_{\beta} - i B_{\beta})$ and X$_{-\beta} = - \frac{1}{2}(A_{\beta} + i B_{\beta})$.  Assertion 2) now follows immediately since A$_{\beta}$* $= -$ A$_{\beta}$ and B$_{\beta}$* $= -$ B$_{\beta}$.

    3)  It suffices to prove that $\tau_{\alpha}$* $ = \tau_{\alpha}$ for all $\alpha \in \Phi$ since these vectors span $\fH_{\br}$. As an element of End(V) we note that $\tau_{\alpha} =
$ [X$_{\alpha}$,X$_{-\alpha}$] $=$ X$_{\alpha} \circ$ X$_{-\alpha} - $ X$_{-\alpha} \circ$ X$_{\alpha}$.  Hence $\tau_{\alpha}$* $ = \tau_{\alpha}$ for all $\alpha \in \Phi$ by 2).

    4)  Let $\mu,\mu^{*}$ be distinct weights in $\Lambda$, and let v$_{\mu}$, v$_{\mu^{*}}$ be arbitrary vectors in V$_{\mu}$, V$_{\mu^{*}}$ respectively.  Since \{$\tau_{\alpha} : \alpha \in \Delta $\} spans $\fH$ there exists $\alpha \in \Phi$ such that $\mu(\tau_{\alpha}) \neq \mu^{*}(\tau_{\alpha})$.  Using 3) we compute $\mu(\tau_{\alpha}) \langle v_{\mu} , v_{\mu^{*}} \rangle = \langle \tau_{\alpha}(v_{\mu}) , v_{\mu^{*}} \rangle = \langle v_{\mu} , \tau_{\alpha}(v_{\mu^{*}}) \rangle = \overline{\mu^{*}(\tau_{\alpha})} \langle v_{\mu} , v_{\mu^{*}} \rangle = \mu^{*}(\tau_{\alpha}) \langle v_{\mu} , v_{\mu^{*}} \rangle$ since $\mu^{*}(\tau_{\alpha})$ is an integer.  We conclude that  $\langle v_{\mu} , v_{\mu^{*}} \rangle = 0$ since  $\mu(\tau_{\alpha}) \neq \mu^{*}(\tau_{\alpha})$.

    We show that V$_{0}$ is orthogonal to V$_{\mu}$ for all $\mu \in \Lambda$.  Given $\mu \in \Lambda$ choose $\alpha \in \Phi$ such that $\mu(\tau_{\alpha}) \neq 0$.  This can be done since $\mu \neq 0$ and \{$\tau_{\alpha} : \alpha \in \Delta$ \} spans $\fH$.  Given arbitrary elements v$_{0} \in$ V$_{0}$ and v$_{\mu} \in$ V$_{\mu}$ we compute $0 = \langle \tau_{\alpha}(v_{0}) , v_{\mu} \rangle = \langle v_{0} , \tau_{\alpha}(v_{\mu}) \rangle = \mu(\tau_{\alpha}) \langle v_{0} , v_{\mu} \rangle$ by 3), which completes the proof of 4).    \end{proof}

$\mathit{\bz-bases~ adapted~ to~ weight~ spaces}$

    Every complex $\fG$-module V admits a basis $\fB$* such that the elements of $\fU(\fG)_{\bz}$ (defined in section 3) leave invariant $\bz$-span($\fB$*).  In particular this applies to the elements of a Chevalley basis $\fC$ of $\fG$ since $\fC \subset \fU(\fG)_{\bz}$.  Moreover, the basis $\fB$* can be chosen to be a union of bases for the weight spaces V$_{0}$ and V$_{\mu}, \mu \in \Lambda$.  More precisely we have

\begin{proposition}  Let $\fG$ be a complex, semisimple Lie algebra, and let V be an irreducible $\fG$-module.  Let $\fH$ be a Cartan subalgebra of $\fG$, and let $\Delta$ be a base for the roots $\Phi$.  Let $\lambda \in \Lambda$ be the unique highest dominant weight of V, and let v be any nonzero vector in the weight space V$_{\lambda}$.  Then

\indent 1)  $\fU(\fG)_{\bz}$(v) is a finitely generated $\bz$-module.

\indent 2)  If $\fB$* is any $\bz$-basis for $\fU(\fG)_{\bz}$(v), then  $\fB$* is a $\bc$-basis for V.  Any element that lies in  $\fU(\fG)_{\bz} \subset$ End(V) has a $\bz$-matrix relative to $\fB$*.

\indent 3)  There exist $\bc$-bases $\fB_{0}$* for V$_{0}$ and $\fB_{\mu}$* for V$_{\mu}, \mu \in \Lambda$ such that the union $\fB$* of these bases is a $\bz$-basis for $\fU(\fG)_{\bz}$(v).
\end{proposition}
\begin{proof}  1) This assertion is proved, for example, in [Hu, p.156].

    2)  Let $\fB$* be any $\bz$-basis for $\fU(\fG)_{\bz}$(v).  The subring $\fU(\fG)_{\bz}$ of End(V) clearly leaves invariant $\bz$-span($\fB$*) = $\fU(\fG)_{\bz}$(v).  In particular the elements of the Chevalley basis $\fC$ of $\fG$ leave invariant $\bz$-span($\fB$*) since $\fC \subset \fU(\fG)_{\bz}$.  It follows that $\fG$ leaves invariant V$'  = \bc$-span($\fB$*), and hence V$'$ = V since V is an irreducible $\fG$-module.  The discussion in [Hu, p.156] shows that the $\bz$-rank of $\fU(\fG)_{\bz}$(v) equals dim V, and hence $\fB$* is a $\bc$-basis for V.

    3)  If M = $\fU(\fG)_{\bz}$(v), where v is a nonzero element of V$_{\lambda}$, then by [Hu, p.156] M = M$_{0} \oplus$ \(\sum_{\mu \in \Lambda} M_{\mu} \) (direct sum of $\bz$-modules), where M$_{0} =$ M $\cap$ V$_{0}$ and M$_{\mu} =$ M $\cap$ V$_{\mu}$ for all $\mu \in \Lambda$.  If $\fB_{0}$* and $\fB_{\mu}$* are $\bz$-bases for M$_{0}$ and M$_{\mu}$ respectively, then the union $\fB$* of the bases \{$\fB_{0}$*, $\fB_{\mu}$*, $\mu \in \Lambda$ \} is a $\bz$-basis for M.  By 2) $\fB$* is a $\bc$-basis for V, and hence the sets $\fB_{0}$* and  $\fB_{\mu}$* are $\bc$-bases for V$_{0}$ and V$_{\mu}$ since they are linearly independent spanning sets for  V$_{0}$ and V$_{\mu}$. \end{proof}

$\mathit{Orthogonality~ of~ weight~ spaces}$

\begin{proposition} Let V be a complex $\fG$-module such that V admits a
nondegenerate, $\fG$-invariant, symmetric,  complex, bilinear form B*. Then

\indent 1)  If $\lambda,\mu \in \Lambda$ , then B*(V$_\lambda$ , V$_\mu$) = $\{0\}$ unless $\lambda$ + $\mu$ = 0.  For every $\lambda \in \Lambda$ and every nonzero element v$_{\lambda} \in$ V$_{\lambda}$ there exists an element v$'_{- \lambda} \in$ V$_{- \lambda}$  such that B*(v$_\lambda$ , v$'_{-\lambda}$) $\neq \{0\}$.

\indent 2)  B*(V$_0$, V$_{\mu}$)  = $\{0\}$ for all $\mu \in \Lambda$.

\indent 3)  The bilinear form B* is nondegenerate on V$_0$, the zero weight space of V.

\indent 4)  If $\lambda \in \Lambda$, then $- \lambda \in \Lambda$.
\end{proposition}

\begin{proof}  1)  Let H $\in \fH$ , v$_{\lambda} \in$ V$_{\lambda}$ and v$_{\mu} \in$ V$_{\mu}$ be given.  Then $\lambda$(H) B*(v$_{\lambda}$ , v$_{\mu}$) = B*(H(v$_{\lambda}$) , v$_{\mu}$) = $-$B*(v$_{\lambda}$, H(v$_{\mu}$)) = $- \mu$(H) B*(v$_{\lambda}$ , v$_{\mu}$).  This proves the first assertion in 1).  This argument also proves 2) if we set $\lambda = 0$.

    If v$_{\lambda}$ is an element of V$_{\lambda}$ such that B*( v$_{\lambda}$, v$_{-\lambda}'$) = 0 for all v$_{-\lambda}' \in$ V$_{-\lambda}$, then B*(v$_{\lambda}$, v) = 0 for all v $\in$ V.  This shows that
v$_{\lambda}$= 0 by the nondegeneracy of B*, which completes the proof of 1).

    3)  By 2) the restriction of B* to V$_0$ must be nondegenerate since B* is nondegenerate on \ws.

    4)  Let $\fG_0$ be any compact real form of $\fG$.  If $\fH_0$ is a maximal abelian subalgebra of $\fG_0$ ,  then
 $\fH_0^{\bc}$ is a Cartan subalgebra of $\fG_0^{\bc} = \fG$.  The automorphism group of $\fG$ acts transitively on the
 set of Cartan subalgebras of $\fG$.  Hence if we replace $\fG_0$ by its image under some automorphism of $\fG$, then
 we may assume that $\fH_0^{\bc} = \fH$.  By (2.1) there exists a real $\fG_0$-module U
such that V = U$^{\bc}$.  Let J denote the conjugation of V induced by U.

    Now let $\lambda \in \Lambda$ and let v $\in$ V$_{\lambda}$.  We show that J(V$_{\lambda}$) = V$_{-\lambda}$.
By symmetry it suffices to show that J(V$_{\lambda}$) $\subseteq$ V$_{-\lambda}$.  Since $\fH_0^{\bc} = \fH$ it
suffices to prove that H$_0$(Jv) = $-\lambda$(H$_0$)(Jv ) for all H$_0 \in$ $\fH_0$.  Write v = u + i u$'$,where u,u$'$ are
elements of U.  Then $\lambda$(H$_0$)(u + i u$'$) = $\lambda$(H$_0$)v = H$_0$(v) = H$_0$(u) + i H$_0$(u$'$).  Since
$\lambda$(H$_{0}$) $\in$ i$\br$ by the second remark in (2.3) we conclude

\indent \indent (*) $\lambda$(H$_0$)u = i H$_0$(u$'$) and $\lambda$(H$_0$)u$' = -$ i H$_0$(u)

From (*) it follows immediately that H$_0$(Jv) = $-\lambda$(H$_0$)(Jv) for all H$_0 \in \fH_0$.
 This proves that J(V$_{\lambda}$) = V$_{-\lambda}$, and a similar argument shows that J(V$_0$) = V$_0$.
 \end{proof}

 $\mathit{Rootspace ~decomposition ~of ~\fG}$

    If V = $\fG$ and the representation is the adjoint representation, then the weights of $\fG$ are the
roots of $\fG$ and we obtain the root space decomposition determined by $\fH$ :

\indent \indent \indent \rs

The spaces $\{\fG_{\alpha}\}$ are known to have complex dimension 1.  If $\mu \in$ Hom($\fH,\bc$) $- \Phi$, then we define $\fG_{\mu}$
to be \{0\}.

    We collect some useful known facts.

\begin{proposition} Let \rs  be the root space decomposition of  a complex, semisimple Lie algebra $\fG$, and
let B : $\fG$ x $\fG \rightarrow \bc$ denote the Killing form on $\fG$.  Then

\indent 1)   If $\alpha,\beta \in \Phi$ , then B($\fG_{\alpha} , \fG_{\beta}) = \{0\}$ unless $\alpha + \beta$ = 0.
For every $\alpha \in \Phi$ , B($\fG_{\alpha}, \fG_{-\alpha}$) $\neq \{0\}$.

\indent 2)  The restriction of B to $\fH$ is nondegenerate.  In particular, for each $\alpha \in \Phi$ there exists a (root) vector
H$_{\alpha} \in \fH$ such that $\alpha$(H) = B(H, H$_{\alpha}$) for all H $\in \fH$.

\indent 3)  The restriction of B to $\fH_{\br} =$ \{H $\in \fH : \alpha(H) \in \br$ for all $\alpha \in \br$\} is
positive definite.  Moreover, H$_{\alpha}$ and $\tau_{\alpha} = 2H_{\alpha} / B(H_{\alpha},H_{\alpha})$ belong to $\fH_{\br}$.

\indent 4)  If X$_{\alpha}$, X$_{-\alpha}$ are nonzero elements of $\fG_{\alpha}, \fG_{-\alpha}$ respectively, then

\indent \indent \indent  [X$_{\alpha}$, X$_{-\alpha}$] = B(X$_{\alpha}$, X$_{-\alpha}$) H$_{\alpha} \neq$ 0
for all $\alpha \in \Phi$.
\end{proposition}
\begin{proof}
    Assertions 1) and 2) follow from 1) and 2) of (4.5) applied to the adjoint
representation.  Assertion 3) follows from [He, p.170].

    We prove 4).  The vectors $\xi_{\alpha}$ = [X$_{\alpha}$, X$_{-\alpha}$] and $\eta_{\alpha}$ =
B(X$_{\alpha}$, X$_{-\alpha}$) H$_{\alpha}$ lie in $\fH$, and the ad $\fG$-invariance of the Killing form B
implies that B(H, $\xi_{\alpha}$) = B(H, $\eta_{\alpha}$) for all H $\in \fH$.  Note that $\eta_{\alpha} \neq$ 0 by
1) and the fact that each $\fG_{\alpha}$is 1-dimensional. Now apply 2) to obtain 4).
\end{proof}

$\mathit{Bracket~ relations~ in~ \fN = V \oplus \fG}$

\begin{proposition} Let $\fG$ be a complex semisimple Lie algebra, and let $\fH$ be a Cartan subalgebra of $\fG$.  Let \{\nvg , B*\} be the complex 2-step nilpotent Lie algebra defined at the beginning of section 2.  Let $\Lambda \subset$ Hom($\fH$,$\bc$) be the set of weights determined by $\fH$, and let \ws ~be the corresponding weight space decomposition.  Then

\indent 1)  Let v$_{\lambda}$ be any nonzero element of V$_{\lambda}$.  Then
ad v$_{\lambda}(V_{-\lambda}) = \bc-span(t_{\lambda})$ for every
$\lambda \in \Lambda$, where  t$_{\lambda}$ is the element in $\fH$ such
 that $\lambda$(H) = B*(t$_{\lambda}$, H) for all H $\in \fH$.

\indent 2)  $[V_{\lambda},V_{\mu}] = \fG_{\lambda+\mu}$ \hspace{.42in} for all $\lambda,\mu \in \Lambda$ with $\lambda + \mu \neq 0$.

\indent 3)  $[V_{0},V_{\mu}] = \fG_{\mu}$ \hspace{.59in} for all $\mu \in \Lambda$

\indent 4)  $[V_0, V_0] =$ \{0\}.
\end{proposition}

\begin{proof}  We recall from the discussion preceding (4.6) that $\fG_{\mu}$ is defined to be \{0\} if $\mu \in$ Hom($\fH,\bc$) $- \Phi$.

    1)  Let v$_{\lambda}$ ~be any nonzero element of V$_{\lambda}$.  We show first that $[v_{\lambda},V_{\lambda}] \subset \fH$ for every $\lambda \in \Lambda$.  By 1) of (4.5) it suffices to show that for $\lambda \in \Lambda$ and $\beta \in \Phi, B^{*}([v_{\lambda},v_{-\lambda}], X_{\beta}) = 0$ for every $v_{\lambda} \in V_{\lambda}, v_{-\lambda} \in V_{-\lambda}$ and $X_{\beta} \in
\fG_{\beta}$.  Now observe that $B^{*}([v_{\lambda},v_{-\lambda}], X_{\beta}) = B^{*}(X_{\beta}(v_{\lambda}), v_{-\lambda}) \in B^{*}(V_{\lambda+\beta}, V_{-\lambda}) =$ \{0\} by 1) of
(4.5).

    Since $v_{\lambda} \neq 0$ we know by 1) of (4.5) that there exists $v_{-\lambda} \in V_{-\lambda}$ such that $B^{*}(v_{\lambda}, v_{-\lambda}) \neq 0$. To complete the proof it suffices to show

\indent \indent (*) $[v_{\lambda},v_{-\lambda}] = B^{*}(v_{\lambda} , v_{-\lambda}) t_{\lambda}$ for all  $v_{\lambda} \in V_{\lambda}$
and all $v_{-\lambda} \in V_{-\lambda}$.

If H $\in \fH$, then $B^{*}([v_{\lambda},v_{-\lambda}], H) = B^{*}(H(v_{\lambda}), v_{-\lambda}) = \lambda(H) B^{*}(v_{\lambda} , v_{-\lambda})
= B^{*}(B^{*}(v_{\lambda} , v_{-\lambda}) t_{\lambda}, H)$.  The proof of (*) is complete since H $\in \fH$ is arbitrary and B* is nondegenerate on $\fH$ by 3) of (4.5).

\indent 2)  Let $\beta \in \Phi$ be any element with $-~\beta \neq \lambda + \mu$. Consider arbitrary
elements $v_{\lambda} \in V_{\lambda}, v_{\mu} \in V_{\mu}$ and $X_{\beta} \in \fG_{\beta}$.  We compute $B^{*}([v_{\lambda},v_{\mu}], X_{\beta}) = B^{*}(X_{\beta}(v_{\lambda}), v_{\mu}) \in B^{*}(V_{\lambda+\beta}, V_{\mu}) =$ \{0\} by 1) of (4.5) and the hypothesis that $-~\beta \neq \lambda + \mu$.  If H $\in \fH$, then $B^{*}([v_{\lambda},v_{\mu}], H) = B^{*}(Hv_{\lambda},~v_{\mu}) = \lambda(H) B^{*}(v_{\lambda}, v_{\mu}) = 0$ by 1) of (4.5) since $\lambda~ + ~\mu \neq$ 0 by hypothesis.

    It follows from the root space decomposition of $\fG$ that $[V_{\lambda},V_{\mu}] \subseteq \fG_{\lambda+\mu}$.  To prove equality in 2) it suffices to consider the case that $\lambda + \mu = \beta \in \Phi$.  ~Let X$_{-\beta}$ be a nonzero element of $\fG_{-\beta}$.  By (4.5) we know that $-~\lambda$
 and $-~\mu$ belong to $\Lambda$, and $\lambda - \beta$ = $-~\mu \in \Lambda$.  By (4.1) we know that X$_{-\beta}$(V$_{\lambda}$) is a nonzero subspace of V$_{-\mu}$.  Choose v$_{\lambda} \in V_{\lambda}$ such that X$_{-\beta}$(v$_{\lambda}$) is nonzero.  It follows from 1) of (4.5)  that there exists an element v$_{\mu} \in$ V$_{\mu}$ such that $B^{*}(X_{-\beta}(v_{\lambda}), v_{\mu}) = B^{*}([v_{\lambda}, v_{\mu}], X_{-\beta})$ is nonzero.  Hence $[v_{\lambda}, v_{\mu}]$ is nonzero, and  $[v_{\lambda}, v_{\mu}]$ spans the 1-dimensional vector space $\fG_{\lambda+\mu}$.

\indent 3)  The proof of this assertion is a small modification of the proof of 2).  We omit the details.

\indent 4)  Let v$_0$ and v$_0'$ be arbitrary elements of V$_0$.  If $X_{\beta} \in \fG_{\beta}$ for $\beta \in \Phi$, then $B^{*}([v_{0}, v_{0}'], X_{\beta}) = B^{*}(X_{\beta}(v_{0}), v_{0}') \in B^{*}(V_{\beta}, V_{0}) =$ \{0\} by 2) of (4.5).  If H $\in \fH$, then $B^{*}([v_{0}, v_{0}'], H) = B^{*}(H(v_{0}), v_{0}') = 0$ since H(v$_0$) = 0.  Hence $B^{*}([v_{0}, v_{0}'] ,\fG) =$ \{0\}, and it follows that $[v_{0}, v_{0}'] = 0$ since B* is nondegenerate on $\fG$.
 \newline
\end{proof}

    Next we obtain some more precise information about the individual linear maps ad v : $\fN = \fV \oplus \fG \rightarrow \fG$.  We use the notation of the preceding result.
\newline

$\mathit{The~ range~ of~ ad v : v \in V}$

\begin{proposition}  Let $\{\fN = V \oplus \fG, B^{*}\}$ be as in (4.7).  If v $\in$ V, then im (ad v) = $\fG_{v}^{\perp}$, where $\fG_{v} = \{X \in \fG : X(v) = 0\}$ and $\fG_{v}^{\perp} = \{X \in \fG : B(X, \fG_{v}) = 0\}$.
\end{proposition}

\begin{proof}  Let v $\in$ V, and set A = im (ad v) $\subset \fG$.  We show that A$^{\perp} = \fG_{v}$.  Since B* is nondegenerate on $\fG$ it will then follow that $A = (A^{\perp})^{\perp}$, as desired.

    Let X $\in \fG$.  Then X $\in A^{\perp} \iff B^{*}([v,v'],X) = 0$ for all v$' \in$ V $\iff B^{*}(X(v),v') = 0$ for all v$' \in$ V $\iff$ X(v) $= 0$ by the nondegeneracy of B* on V.
\end{proof}


$\mathit{Surjectivity~ of~ ad~ v}$  For an element v $\in$ V we note that ad v($\fN$) = ad v(V) since $\fG$ lies in the center of $\fN = V \oplus \fG$.  If ad v : V $\rightarrow \fG$ is surjective for some v $\in$ V, then clearly we must have dim V $\geq$ dim $\fG$.  Conversely, if V is an irreducible $\fG$- module with dim V $\geq$ dim $\fG$, then with only a few exceptions there exists a Zariski open subset O $\subset$ V such that ad v : V $\rightarrow \fG$ is surjective for all v $\in$ V.  By the result above it suffices to show that except in a few cases there exists a Zariski open subset O $\subset$ V such that $\fG_{v} = \{0\}$ for all v $\in$ V.  A list of these exceptional cases may be found in [KL]  and [El].
\newline

$\mathit{The~ range~ of~ ad v_{\lambda} : v_{\lambda} \in V_{\lambda}}$

    For a nonzero element v$_{\lambda}$ in a weight space V$_{\lambda}$ we can say more about the image of ad v$_{\lambda}$.  We first introduce some terminology, and for each $\lambda \in \Lambda$ we define a certain proper Zariski closed subset $\Sigma_{\lambda} \subset V_{\lambda}$.

    Given $\lambda \in \Lambda$ we define $\Phi_{\lambda} = \{\alpha \in \Phi : V_{\lambda - \alpha} \neq \{0\}\}$. Recall that t$_{\lambda}$ is the element in the Cartan subalgebra $\fH$ such that B*(H,t$_{\lambda}$) $= \lambda$(H) for all H $\in \fH$.

    If $\alpha \in \Phi_{\lambda}$, then by definition V$_{\lambda - \alpha} \neq \{0\}$, and by (4.1) it follows that X$_{- \alpha} : V_{\lambda} \rightarrow V_{\lambda - \alpha}$ is a nonzero linear map.  If $\lambda \in \Lambda \cap \Phi$, then by (4.1) and (4.2) we know that V$_{0} \neq \{0\}$ and X$_{- \lambda} : V_{\lambda} \rightarrow V_{0}$ is a nonzero linear map.  We let Ker(X$_{- \lambda}$) denote the kernel of this map.  If $\lambda \in \Lambda \cap \Phi$, then we define a proper Zariski closed subset $\Sigma_{\lambda}$ of V$_{\lambda}$ by $\Sigma_{\lambda} = Ker(X_{- \lambda})~ \cup ~\bigcup_{\alpha \in \Phi_{\lambda}} Ker(X_{- \alpha})$.  If $\lambda \notin \Lambda \cap \Phi$, then we define a proper Zariski closed subset $\Sigma_{\lambda}$ of V$_{\lambda}$ by $\Sigma_{\lambda} = \bigcup_{\alpha \in \Phi_{\lambda}} Ker(X_{- \alpha})$.

\begin{proposition}  Let $\lambda \in \Lambda$ and let v$_{\lambda}$ be any element of V$_{\lambda}$.  Then ad v$_{\lambda}$(V) $\subset \bc t_{\lambda} \oplus \fG_{\lambda} \oplus$ \(\sum_{\alpha \in \Phi_{\lambda}} \fG_{\alpha} \).  Equality holds if $v_{\lambda} \in V_{\lambda} - \Sigma_{\lambda}$.
\end{proposition}

\begin{proof}  We write $V = V_{0}~ \oplus$ \(\sum_{\mu \in \Lambda} V_{\mu}\).  If v$_{\lambda} \in$ V$_{\lambda}$, then $ad(v_{\lambda})(V) \subset ad(v_{\lambda})(V_{0}) +$ \(\sum_{\mu \in \Lambda} ad(v_{\lambda})(V_{\mu})\) $\subset \fG_{\lambda} +$ \(\sum_{\mu \in \Lambda} \fG_{\lambda + \mu}\).  Now $\fG_{\lambda + \mu} = \{0\}$ unless $\lambda + \mu = \alpha \in \Phi$ but this occurs $\iff  \alpha \in \Phi_{\lambda}$.  From 1) of (4.7) we see that ad(v$_{\lambda}$)(V$_{- \lambda}$) $ = \bc t_{\lambda}$.  This proves that  ad v$_{\lambda}$(V) $\subset \bc t_{\lambda} \oplus \fG_{\lambda} \oplus$ \(\sum_{\alpha \in \Phi_{\lambda}} \fG_{\alpha} \).

    Now suppose that v$_{\lambda} \in$ V$_{\lambda} - \Sigma_{\lambda}$.  By the definition of $\Sigma_{\lambda}$ this means that $X_{- \lambda}(v_{\lambda}) \neq 0$ and
$X_{-\alpha}(v_{\lambda}) \neq 0$ for all $\alpha \in \Phi_{\lambda}$.  Since the subspaces $\fG_{\alpha}$ are 1-dimensional for all $\alpha \in \Phi$ the equality assertion will follow once we show that $ad(v_{\lambda})(V_{0}) \neq \{0\}$ for all $\lambda \in \Lambda~ \cap~ \Phi$ and $ad(v_{\lambda})(V_{\mu}) \neq \{0\}$ for all $\mu \in \Lambda$ such that $\lambda + \mu = \alpha \in \Phi_{\lambda}$.

    Let $\lambda \in \Lambda~ \cap~ \Phi$.  Then $ad(v_{\lambda})(V_{0}) = \{0\} \iff B^{*}([v_{\lambda},V_{0}], X_{- \lambda}) = 0 \iff B^{*}(X_{- \lambda}(v_{\lambda}), V_{0}) = 0
\iff X_{- \lambda}(v_{\lambda}) = 0$ by the nondegeneracy of B* on V.  Since $X_{- \lambda}(v_{\lambda}) \neq 0$ it follows that $ad(v_{\lambda})(V_{0}) \neq \{0\}$.

    Similarly, $ad(v_{\lambda})(V_{\mu}) = 0$ for $\lambda + \mu = \alpha \in \Phi_{\lambda} \iff
B^{*}([v_{\lambda}, V_{\mu}], X_{- \alpha}) = 0 \iff B^{*}(X_{- \alpha}(v_{\lambda}), V_{\mu}) = 0 \iff
X_{- \alpha}(v_{\lambda}) = 0$.  Since $X_{- \alpha}(v_{\lambda}) \neq 0$ it follows that $ad(v_{\lambda})(V_{\mu}) \neq \{0\}$.
\end{proof}

$\mathit{Automorphisms~ and~ derivations~ of~ \fN = V\oplus \fG}$

    We do not have a complete description of the group Aut($\fN$), where $\fN = V \oplus \fG$, but we show next that Aut($\fN$) contains
a subgroup whose Lie algebra is isomorphic to $\fG$.  See [L,Theorem 3.12] for a metric version of this result in the case
of a real 2-step nilpotent Lie algebra \nug.  In the remainder of this section we regard $\fG$ as a subalgebra
of End(V), where V is a finite dimensional, complex $\fG$-module.

\begin{proposition}  Let $\fN = V \oplus \fG$, and let G be the connected Lie subgroup  of SL(V) whose Lie algebra is $\fG$.  For g $\in$ G let T(g) denote the element of GL($\fN$) such that T(g) = g on V and T(g) = Ad(g) on $\fG$.  Then T(g) $\in$ Aut($\fN$) for all g $\in$ G, and T : G $\rightarrow$ Aut($\fN$) is an injective homomorphism.
\end{proposition}

$\mathit{Remark}$   We recall that Ad(g)(X) = gXg$^{-1}$ for g $\in$ G and X $\in \fG$ since $\fG$ is a
subalgebra of End(V).  Since G is a semisimple Lie group it is a standard fact that G $\subseteq$ SL(V).  Moreover,
Ad(G) = Aut($\fG$)$_0$ since these are connected subgroups of Aut($\fG$) with the same Lie algebra
([He, Corollary 6.5, p. 132]).  The kernel of Ad, which is the center of G, is discrete since G is semisimple.  If V is an irreducible
$\fG$-module of dimension n, then Ker Ad $\subseteq$ \{e$^{i2 \pi k/n}$ Id : 1 $\leq$ k $\leq$ n\} since
G $\subseteq$ SL(V) and $\fG$ leaves invariant the eigenspaces of elements in Ker Ad.

\begin{proof}  Let B* : $\fN$ x $\fN \rightarrow \bc$ be the symmetric nondegenerate, bilinear form that
defines the bracket structure on \nvg.  We shall need a preliminary result.
\newline

\begin{lem} (G-invariance of B*)

    1)  B*(v,w) = B*(gv, gw) for all v,w $\in$ V and g $\in$ G.

    2)  B*(Ad(g)X, Ad(g)Y) = B*(X,Y) for all X,Y $\in \fG$ and g $\in$ G.
\newline
\end{lem}

$\mathit{Proof~ of~ the~ Lemma}$

    We prove only 1) since the proof of 2) is virtually identical.  To prove 1) it suffices to consider elements of the form g = exp(X), X $\in \fG$, since such elements generate the connected group G.  Given X $\in \fG$ and v,w $\in$ V we define f(t) = B*(exp(tX)v, exp(tX)w) for t $\in \br$.  If we set v(t) = exp(tX)v and w(t) = exp(tX)w, then f$'$(t) = B*(Xv(t), w(t)) + B*(v(t), Xw(t)) = 0 since B* is $\fG$-invariant.  Hence B*(gv, gw) = f(1) = f(0) = B*(v,w).  This completes the proof of the Lemma.

    We now complete the proof of the Proposition.  Let g $\in$ G and v,w $\in$ V be given. To show that
T(g) $\in$ Aut($\fN$) it suffices to prove that [gv, gw] = Ad(g)[v,w] since $\fN$ is 2-step nilpotent and
[$\fN$, $\fN$] = [V,V] = $\fG$.  It is routine to show that T : G $\rightarrow$ Aut($\fN$) is a homomorphism, and obvious that T is injective.

    Let X $\in \fG$ be given.  By the observations at the beginning of the proof and the lemma above we obtain B*([gv, gw] , X) = B*(X(gv), gw) = B*(g \{Ad(g$^{-1}$)(X)\}(v), gw ) =
B*(Ad(g$^{-1}$)(X)(v), w) = B*([v, w], Ad(g$^{-1}$)(X)) = B*(Ad(g)([v, w]), X).  It follows that [gv, gw] = Ad(g)[v,w] since X $\in \fG$ was arbitrary and B* is nondegenerate on $\fG$.
\end{proof}

    We conclude this section with an immediate consequence of the preceding proposition.

\begin{corollary}  Let $\fG$ be a complex semisimple Lie algebra, and let V be a finite dimensional complex $\fG$-module.  Let $\fN = V~ \oplus~ \fG$.  For X $\in \fG$ let t(X) denote the element of End($\fN$) such that t(X) = X on V and t(X) = ad X on $\fG$.  Then t(X) $\in$ Der($\fN$) for all X, and t : $\fG  \rightarrow$ Der($\fN$) is an injective Lie algebra homomorphism.
 \end{corollary}

\section{Action~ of~ the~ Weyl~ group~ by~ unitary~ automorphisms}

    Let V be a complex $\fG$-module, and let G be the connected Lie subgroup of SL(V) whose Lie algebra is $\fG$, regarded as a subalgebra of End(V).  Every Cartan subalgebra $\fH$ and Chevalley basis $\fC$ of $\fG$ determines a discrete subgroup G$_{\bz}$ of G that is self adjoint with respect to any $\fG_{0}$-invariant Hermitian inner product on V.  Moreover, the group G$_{\bz}$ contains elements related to the Weyl group determined by $\fH$.  If V = U$^{\bc}$, where U is a real $\fG_{0}$-module, then we show that these elements leave U invariant.
\newline

$\mathit{The~ group~ G_{\bz}}$

     Let $\fC$ = \{X$_{\beta}$ : $\beta \in \Phi$ ; $\tau_{\alpha}$ : $\alpha \in \Delta$ \} be a Chevalley basis of $\fG$, and let $\fU(\fG)_{\bz}$ be the corresponding subring of the universal enveloping algebra of $\fG$ (cf.section 3).  Here we regard $\fU(\fG)_{\bz}$ as a subring of $\fU(\fG) \subset$ End(V).  We recall from (4.4) that if V is any irreducible complex $\fG$-module, then for any nonzero highest weight vector v there exists a basis $\fB$ of V such that $\bz$-span ($\fB$) $= \fU(\fG)_{\bz}$(v).

     We define G$_{\bz}$ to be the subgroup of G generated by the elements \{g$_{\alpha} =$ exp(X$_{\alpha}$) : $\alpha \in \Phi$ \}.  By 1) of (4.3) it follows that g$_{\alpha}$* $=$ g$_{-\alpha}$ for all $\alpha \in \Phi$.  Hence G$_{\bz}$ is a self adjoint subgroup of G since the generating set of G$_{\bz}$ is invariant under *.  By a remark preceding (3.3) the generating set for G$_{\bz}$ is contained in $\fU(\fG)_{\bz}$, and hence G$_{\bz}$ leaves invariant $\bz$-span ($\fB$) $= \fU(\fG)_{\bz}$(v).  In particular G$_{\bz}$ is a discrete subgroup of G.  If $\langle , \rangle$ is any $\fG_{0}$ - invariant Hermitian inner product on V, then the subgroup G$'_{\bz}$ of G$_{\bz}$ consisting of unitary elements is therefore a finite subgroup of G.  We show next that  G$'_{\bz}$ contains elements that permute the weight spaces of V like elements of the Weyl group.
\newline

$\mathit{The~ Weyl~ group~ and~ complex~ \fG-modules}$

    We recall that for $\alpha \in \Phi$ the Weyl reflection $\sigma_{\alpha}$ acts on $\fH$
by $\sigma_{\alpha}$(H) = H $- ~ \alpha(H) \tau_{\alpha}$ and on Hom ($\fH$,$\bc$) by $\sigma_{\alpha}$($\lambda$) = $\lambda~ - ~\lambda(\tau_{\alpha})\alpha$.  More generally, the elements of W act on Hom ($\fH$,$\bc$) by $\sigma$($\lambda$)(H) = $\lambda$($\sigma^{-1}$H) for all $\sigma \in$ W, $\lambda \in$ Hom($\fH$, $\bc$) and H $\in \fH$.

\begin{proposition}  Let V be a complex $\fG$ - module.  Let $\fH$ be a Cartan subalgebra of $\fG$ with root system $\Phi \subset$ Hom($\fH$,$\bc$).  Let W $\subset$ GL($\fH$) be the Weyl group.  For each $\sigma \in$ W there exists an element T$_{\sigma}$ in G$_{\bz}$ such that

\indent 1)  T$_{\sigma} \circ$ H $\circ$ T$_{\sigma}^{-1}  = \sigma$(H) on V for all H $\in \fH$.

\indent 2)  T$_{\sigma}$(V\submu) = V$_{\sigma(\mu)}$   for all $\mu \in \Lambda$.

$\hspace{.15in}$T$_{\sigma}$(V$_0$) = V$_0$.

\indent 3)   T$_{\sigma}$ is unitary with respect to any $\fG_{0}$-invariant Hermitian inner product
$\< , \>$ on V.
\end{proposition}

$\mathit{Remark}$ T$_{\sigma}$ is not uniquely determined by properties 1), 2) and 3).  For example, let H$'$ be any element of  $\bz$-span \{$\tau_{\alpha} : \alpha \in \Delta$\}. Let H $= \pi i$ H$'$, and let a = exp(H).)  If $\mu \in \Lambda$, then a = exp($\mu$(H) Id $= \pm$ Id on the weight space V$_{\mu}$ since $\mu$(H) $\in \pi i \bz$.  Clearly a commutes with the elements of $\fH$.  Moreover, a is unitary with respect to any $\fG_{0}$-invariant Hermitian inner product $\< , \>$ on V since the weight spaces of V are orthogonal relative to $\< , \>$ by 4) of (4.3).  If T$_{\sigma}' =$ T$_{\sigma} \circ$ a, then it is easy to check that T$_{\sigma}'$  satisfies the three conditions of the proposition.   However, this is essentially the only way that T$_{\sigma}$ fails to be uniquely determined by 1), 2) and 3).  It follows from 3) of Lemma B below that T$_{\sigma}$ is unique up to composition with a transformation a that is $\pm$ Id on each weight space of each irreducible $\fG$-submodule of V.

\begin{proof}
    It suffices to prove this result in the case that V is an irreducible $\fG$-module.  Moreover, it suffices to prove the proposition for the generating elements of W, the reflections S =
\{$\sigma_{\alpha}~:~\alpha \in \Phi \}$.  If $\sigma \in$ W is arbitrary we then write
$\sigma$ = $\sigma_{1} \circ ...\circ \sigma_{k}$ , where $\sigma_{i} \in$ S, and define T$_{\sigma}$ =
T$_{\sigma_{1}} \circ ... \circ$ T$_{\sigma_{k}}$.

    Given $\alpha \in \Phi$ we let X$_{\alpha} \in \fG_{\alpha}$ be given arbitrarily and choose an element X$_{-\alpha} \in \fG_{-\alpha}$ such that  [X$_{\alpha}$, X$_{-\alpha}$] = $\tau_{\alpha}$.  Define T$_{\alpha}$ = exp(X$_{\alpha}$) exp($-$X$_{-\alpha}$) exp(X$_{\alpha}$) $\in$ G.  By the remarks in section 3 the transformation T$_{\alpha}$ belongs to G $\cap~ \fU(\fG)_{\bz}$, and T$_{\alpha}$ belongs to G$_{\bz}$ by the definitions of T$_{\alpha}$ and G$_{\bz}$.

    We recall that for X, Y $\in \fG$ one has

\indent ({\#}) \hspace{.5in}        exp(X) $\circ$ Y $\circ$ exp($-$X) = exp(ad X) (Y)

\noindent since Ad(exp(X)) = exp(ad X).

    1)   For each $\beta \in \Phi$ it follows from (\#) above that exp(X$_{\beta}$) $\circ$ Y $\circ$ exp($-$X$_{\beta}$)
 = exp(ad X$_{\beta}$)(Y) for all Y $\in \fG$.  Let $\varphi_{\alpha}$= exp(ad X$_{\alpha}$) $\circ$
 exp(ad $-$X$_{-\alpha}$) $\circ$ exp(ad X$_{\alpha}$).  Note that $\varphi_{\alpha}$ lies in Aut($\fG$) since ad X$_{\alpha}$ and ad $-$X$_{-\alpha}$ are derivations of $\fG$.  By the definition of T$_{\alpha}$ we obtain T$_{\alpha} \circ$ Y $\circ$ T$_{\alpha}^{-1}$ = $\varphi_{\alpha}$(Y) for
 all Y $\in \fG$.

    To complete the proof of 1) it remains to show that $\varphi_{\alpha}$(H) = $\sigma_{\alpha}$(H) for all H $\in \fH$.
For completeness we include a slightly expanded version of the proof
of [FH, p. 497].  It is clear that  $\varphi_{\alpha}$(H)
 = $\sigma_{\alpha}$(H) = H if $\alpha$(H) = 0 by the definition of $\sigma_{\alpha}$ and the fact that  ad X$_{\alpha}$(H)
 = ad$-X_{-\alpha}$(H) = 0.
 It suffices to show that $\varphi_{\alpha}$($\tau_{\alpha}$) = $\sigma_{\alpha}$($\tau_{\alpha}$) = $-\tau_{\alpha}$
since $\fH$ is spanned by $\tau_{\alpha}$ and those vectors H satisfying $\alpha$(H) = 0.

    Since ad X$_{\alpha}$($\tau_{\alpha}$) = $-$2X$_{\alpha}$ it follows that (ad X$_{\alpha}$)$^m$($\tau_{\alpha}$) = 0
for all m $\geq$ 2.  Hence

\indent (i)  exp(ad X$_{\alpha}$)($\tau_{\alpha}$) = $\tau_{\alpha} - $ 2X$_{\alpha}$

Since $\alpha$($\tau_{\alpha}$) = 2 and [X$_{\alpha}$, X$_{-\alpha}$] = $\tau_{\alpha}$ we obtain
(ad $-$X$_{-\alpha}$)($\tau_{\alpha} -2$ X$_{\alpha}$) = $-$2$\tau_{\alpha}$ $-$2X$_{-\alpha}$ and \linebreak
(ad $-$X$_{-\alpha}$)$^{2}$($\tau_{\alpha} -2$ X$_{\alpha}$) = 4X$_{-\alpha}$.  It follows that
(ad $-$X$_{-\alpha}$)$^{m}$($\tau_{\alpha} -2$ X$_{\alpha}$) = 0 for all m $\geq$ 3. From (i) and the discussion above we obtain

\indent (ii) exp(ad $-$X$_{-\alpha}$) $\circ$ exp(ad X$_{\alpha}$) ($\tau_{\alpha}$) = $-\tau_{\alpha}$ $-$2X$_{\alpha}$

Finally ad X$_{\alpha}$($-\tau_{\alpha}-$2X$_{\alpha}$) =
2X$_{\alpha}$, which shows that (ad
X$_{\alpha}$)$^{m}$($-\tau_{\alpha}-$2X$_{\alpha}$) = 0 for all m
$\geq$ 2.  From (ii) we obtain $\varphi_{\alpha}(\tau_{\alpha}) =$
exp(ad X$_{\alpha})(-\tau_{\alpha}-2X_{\alpha}) =
(-\tau_{\alpha}-2X_{\alpha}) + 2X_{\alpha} = -\tau_{\alpha} =
\sigma_{\alpha}(\tau_{\alpha})$ The proof of 1) is complete.

\indent 2)  From 1) we obtain

\indent \indent     a) T$_{\alpha} \circ$ H $\circ$ T$_{\alpha}$$^{-1}$ = $\sigma_{\alpha}$(H) =
H $- \alpha$(H) $\tau_{\alpha}$  in V for all H $\in \fH$.

From a) and the fact that $\alpha$($\tau_{\alpha}$) = 2 we obtain

\indent \indent     b)  T$_{\alpha}~ \circ~ \tau_{\alpha} = -~\tau_{\alpha}~ \circ$ T$_{\alpha}$ in V.

From a) and b) it follows that H $\circ$ T$_{\alpha}$ = T$_{\alpha}~ \circ H +~ \alpha(H) (\tau_{\alpha}~ \circ$ T$_{\alpha}$)
= T$_{\alpha}~ \circ H-~ \alpha(H)$ (T$_{\alpha} \circ \tau_{\alpha}$) for all H $\in \fH$.  If v $\in V_{0}$ is arbitrary,
then H(T$_{\alpha}$v) = 0 by the equation above, which proves that T$_{\alpha}(V_{0}) \subseteq V_{0}$.  Equality holds since T$_{\alpha}$ is invertible.

    Given $\mu \in \Lambda$ and H $\in \fH$ we know that H = $\mu$(H) Id on V$_{\mu}$. If we set H* = $\sigma_{\alpha}$(H), then from 1) and the discussion preceding the statement of (5.1) we see that on the subspace T$_{\alpha}$(V$_{\mu}$) one has H* = T$_{\alpha}~ \circ$ H $\circ$ T$_{\alpha}$$^{-1}$ = $\mu$(H) Id = ($\mu~ \circ~ \sigma_{\alpha} ^{-1}$)(H*) Id = $\sigma_{\alpha}$($\mu$)(H*) Id.  Hence T$_{\alpha}$(V$_{\mu}$) $\subseteq$ V$_{\sigma_{\alpha}(\mu)}$ and dim (V$_{\mu}$) $\leq$ dim (V$_{\sigma_{\alpha}(\mu)}$). Equality holds by symmetry, which completes the proof of 2).

\indent 3)  We shall need two preliminary results.  We recall from (4.4) that if V is any irreducible complex $\fG$-module, then for any nonzero highest weight vector v there exists a basis $\fB$ of V such that $\bz$-span ($\fB$) $= \fU(\fG)_{\bz}$(v).

    Assertion 3) of the proposition will follow immediately from 2) of Lemma A and 2) of Lemma B.

$\mathbf{Lemma~ A}$  Let V be a complex $\fG$-module, and let $\< , \>$ be a $\fG_{0}$-invariant Hermitian inner product on V.  Let $\fB$ be a basis of V such that $\fU(\fG)_{\bz}$ leaves invariant $\bz$-span ($\fB$). Let $\alpha \in \Phi$.  Then

\indent 1) T$_{\alpha}^{*} = $ T$_{-\alpha}$.

\indent 2)  If $\psi_{\alpha} =$ T$_{\alpha} \circ$ T$_{\alpha}^{*}$, then $\psi_{\alpha}$ is self adjoint and positive definite, commutes with $\fH$ and leaves invariant $\bz$-span ($\fB$).

$\mathbf{Lemma~B}$  Let V be an irreducible complex $\fG$ - module, and let $\fB$ be as in Lemma A.  Let g $\in$ G be an element that commutes with $\fH$ and leaves invariant $\bz$-span ($\fB$).  Then

\indent 1)  For $\alpha \in \Phi$ and X$_{\alpha} \in \fG_{\alpha}$ we have g $\circ$ X$_{\alpha} \circ$ g$^{-1} = \pm$ X$_{\alpha}$ on V.

\indent 2)  If g is self adjoint and positive definite, then g = Id.

\indent 3)  If g is unitary, then g $= \pm$ Id on each weight space of V.

$\mathit{Proof~ of~ Lemma~ A}$  1)  For $\alpha \in \Phi$ we know that X$_{\alpha}^{*} =$ X$_{-\alpha}$ by 2) of (4.3).  Hence T$_{\alpha}^{*} =$  \{exp(X$_{\alpha}$)~ exp($-$X$_{-\alpha}$)~ exp(X$_{\alpha}$)\} $^{*} =$ exp(X$_{\alpha}$*) exp($-$X$_{-\alpha}$*) exp(X$_{\alpha}$*) $=$ exp(X$_{-\alpha}$) exp($-$X$_{\alpha}$) exp(X$_{-\alpha}$) $=$ T$_{-\alpha}$.

    2)  Clearly  $\psi_{\alpha} =$ T$_{\alpha} \circ$ T$_{\alpha}^{*}$ is self adjoint and positive definite.  By 1) it follows that $\psi_{\alpha} =$ T$_{\alpha} \circ$ T$_{-\alpha}$.  Hence  $\psi_{\alpha}$ commutes with $\fH$ by condition 1) of the proposition.  Moreover, by the remarks in section 3 both T$_{\alpha}$ and T$_{-\alpha}$ belong to G $\cap~ \fU(\fG)_{\bz}$, and hence $\psi_{\alpha} \in$ G $\cap~ \fU(\fG)_{\bz}$.  It follows that $\psi_{\alpha}$ leaves invariant $\bz$-span ($\fB$) $= \fU(\fG)$(v).

$\mathit{Proof~ of~ Lemma~ B}$  1)  Since $\fG_{\alpha}$ is 1-dimensional it suffices to consider the case that X$_{\alpha} \in \fC~ \cap~ \fG_{\alpha}$, where $\fC$ is a Chevalley basis of $\fG$.  Clearly we may assume that X$_{\alpha}$ is nonzero on V, for otherwise 1) is immediate.

    We show first that g $\circ$ X$_{\alpha} \circ$ g$^{-1}$ leaves invariant $\bz$-span ($\fB$).  Since X$_{\alpha} \in \fC \subset \fU(\fG)_{\bz}$ it follows that X$_{\alpha}$ leaves invariant $\bz$-span ($\fB$) $= \fU(\fG)$(v).  Note that 1 = det g since g $\in$ G $\subset$ SL(V), and g leaves invariant  $\bz$-span ($\fB$) by hypothesis.  Hence g$^{-1}$ leaves invariant $\bz$-span ($\fB$), and we conclude that g $\circ$ X$_{\alpha} \circ$ g$^{-1}$ leaves invariant $\bz$-span ($\fB$).

    We observe next that  g $\circ$ X$_{\alpha} \circ$ g$^{-1} =$ c$_{\alpha}$ X$_{\alpha}$ for some complex number c$_{\alpha}$.  If c$_{g}$ denotes conjugation by g, then c$_{g} =$ Ad(g) defines an automorphism of $\fG$, the Lie algebra of $\fG$.  By hypothesis c$_{g}$ is the identity on $\fH$.  Hence c$_{g}$ leaves invariant each root space $\fG_{\alpha}$, and since each $\fG_{\alpha}$ is 1-dimensional it follows that c$_{g} =$ c$_{\alpha}$ Id on $\fG_{\alpha}$ for some complex number c$_{\alpha}$.

    To complete the proof of 1) we need to show that c$_{\alpha} = \pm 1$.  Recall that $\bz$-span ($\fB$) is a free $\bz$-module with basis $\fB$, and $\fB$ is also a $\bc$ - basis of V.  Since X$_{\alpha}$ is nonzero on V there exists v $\in \fB$ such that $0 \neq X_{\alpha}(v) = v'$.  By the discussion above both X$_{\alpha}$ and g $\circ$ X$_{\alpha} \circ$ g$^{-1}$ leave invariant  $\bz$-span ($\fB$).  In particular, $v'$ and c$_{\alpha} v' =$ (g $\circ$ X$_{\alpha} \circ$ g$^{-1}$)(v) lie in $\bz$-span ($\fB$), which proves that $c_{\alpha} \in \bq$.  Note that g$^{n} \circ X_{\alpha} \circ$ g$^{-n} = (c_{\alpha})^{n} X_{\alpha}$ for all positive integers n.  Hence $c_{\alpha} = \pm 1$ since $X_{\alpha}$ and  g$^{n} \circ X_{\alpha} \circ$ g$^{-n}$ have the same eigenvalues for all n.

    We prove 2).  Let g be a positive definite, self adjoint element of G that commutes with $\fH$ and leaves invariant $\bz$-span ($\fB$).  We show first that g commutes with the elements of $\fG$.  By the root space decomposition of $\fG$ and the fact that g commutes with $\fH$ it suffices to show that g commutes with the elements of $\fG_{\alpha}$ for every $\alpha \in \Phi$.

    Let $\alpha \in \Phi$, and let X$_{\alpha} \in \fC~ \cap~ \fG_{\alpha}$ be given.  We may assume that X$_{\alpha}$ is nonzero on V since otherwise g clearly commutes with X$_{\alpha}$.  Note that V is a direct sum of eigenspaces of g since g is self adjoint.  Hence there exists a subspace V$'$ of V such that X$_{\alpha}$ is nonzero on V$'$ and g = c Id on V$'$ for some positive real number c.  Choose v$' \in$ V$'$ so that X$_{\alpha}$(v$'$) $\neq 0$.  From 1) we know that g $\circ$ X$_{\alpha} =$ c$_{\alpha}$ X$_{\alpha} \circ$ g, where c$_{\alpha} = \pm 1$.  Hence g(X$_{\alpha}$(v$'$)) $=$ c$_{\alpha}$ X$_{\alpha} \circ$ (gv$'$) $=$
(c c$_{\alpha}$)(X$_{\alpha}$(v$'$)).  The eigenvalues c and cc$_{\alpha}$ for g are both positive since g is self adjoint and positive definite.  Hence c$_{\alpha}$ is positive, and we conclude that c$_{\alpha} = 1$ for all $\alpha \in \Phi$.  Hence g commutes with $\fG$.

    Since g commutes with $\fG$ it follows that $\fG$ leaves invariant every eigenspace of g in V.  Since V is an irreducible $\fG$-module and g is self adjoint and  positive definite we conclude that g = c Id for some positive real number c.  Finally c =1 since 1 = det g = c$^{dim~ V}$.  This completes the proof of 2).

    3)  Let g $\in$ G be a unitary element that commutes with $\fH$ and leaves invariant $\bz$-span ($\fB$).    By 1) we know that g $\circ$ X$_{\alpha} \circ$ g$^{-1} =$ c$_{\alpha}$ X$_{\alpha}$, where c$_{\alpha} = \pm 1$.  Let $\Delta =$ \{$\alpha_{1}, \alpha_{2},~...~, \alpha_{n}$\} be the given base for $\Phi$.  For $1 \leq k \leq n$ let N$_{k} = 1$ if c$_{\alpha_{k}} = -1$ and let  N$_{k} = 2$ if c$_{\alpha_{k}} = 1$.  Let H $\in \fH$ be that element such that  $\alpha_{k}$(H) $=$ N$_{k} \pi i$ for $1 \leq k \leq n$.  This can be done since H $\rightarrow (\alpha_{1}(H),~ ...,~ \alpha_{n}(H))$ is a linear isomorphism of $\fH$ onto $\bc^{n}$.  Let a = exp(H).  We assert

    (i)  Ad(a) = c$_{\alpha_{k}}$ Id on $\fH$ and b = ga$^{-1}$  commutes with the elements of $\fG$.

    (ii)  There exists a constant c$_{1}$ such that a $= \pm $  c$_{1}$ Id on each nonzero weight space of V.

    (iii)  Let b = ga$^{-1}$ as in (i).  There exists a constant c$_{2}$ such that b = c$_{2}$ Id on V.

    Assuming (i), (ii) and (iii) for the moment, we complete the proof of 3).  It follows from (i), (ii) and (iii) that g $= \pm$ c$_{3}$ Id on each nonzero weight space of V, where  c$_{3} =$ c$_{1}$ c$_{2}$.  By the hypotheses of Lemma B the element g leaves invariant $\bz$-span ($\fB$) $= \fU(\fG)_{\bz}(v)$ for some nonzero highest weight vector v.  By 3) of (4.4) we may assume that $\fB$ is a union of bases $\fB_{0}$ for V$_{0}$ and $\fB_{\mu}$ for V$_{\mu}, \mu \in \Lambda$.  Since g commutes with $\fH$ it leaves invariant V$_{\mu}$ and $\bz$-span ( $\fB_{\mu}$) $=$ V$_{\mu}~ \cap~ \bz$-span ( $\fB$).  Hence c$_{3} \in \bz$.  Applying the same argument to g$^{-1}$, which also satisfies the hypotheses of Lemma B since det(g) $= 1$, we find that g$^{-1} = \pm$ c$_{-3}$ Id on each nonzero weight space of V, where c$_{-3} \in \bz$. Since c$_{3}$ c$_{-3} = \pm 1$ it follows that c$_{3} = \pm 1$.  Since a = exp(H) is the identity on V$_{0}$ we conclude from (iii) that g $=$ c$_{2}$ Id on V$_{0}$.  Similarly, it follows that c$_{2} = \pm 1$ since g and g$^{-1}$ leave invariant $\bz$-span ($\fB_{0}$).  This completes the proof of 3).

    We now prove (i), (ii) and (iii).

$\mathbf{Proof~ of~ (i)}$  By the definition of N$_{k}$ and H it follows that exp($\alpha_{k}$(H)) $=$ c$_{\alpha_{k}}$ for $1 \leq k \leq n$.  On the root space $\fG_{\alpha_{k}} \subset$ End(V) we have Ad(a) = exp(ad(H)) =  exp($\alpha_{k}$(H)) Id = c$_{\alpha_{k}}$ Id = Ad(g).  Hence b commutes with $\fG_{\alpha_{k}}$ for $1 \leq k \leq n$ since Ad(b) = Ad(a) Ad(g)$^{-1}$ = Id on $\fG_{\alpha_{k}}$.  Note that b commutes with $\fH$ since g and a have this property.  In particular, if c$_{b}$ denotes conjugation by b, then c$_{b}$ is the identity on $\fH$ and must therefore leave invariant each root space of $\fG$.  For $\alpha \in \Phi$ let  d$_{\alpha}$ be that complex number such that c$_{b} =$ d$_{\alpha}$ Id on $\fG_{\alpha}$.

    The Lie algebra $\fG$ is generated by the root spaces \{$\fG_{\alpha} : \alpha \in \Delta \cup -\Delta$ \} (cf. section 14.2 of  [Hu]).  To show that b commutes with $\fG$ it therefore suffices to show that b commutes with $\fG_{-\alpha}$ for all $\alpha \in \Delta$. Let $\alpha \in \Delta$ be given, and let X$_{\alpha}$ and X$_{-\alpha}$ be the elements of $\fG_{\alpha}$ and $\fG_{-\alpha}$ that belong to the fixed Chevalley basis $\fC$.  Since c$_{b}$ fixes X$_{\alpha}$ and $\tau_{\alpha} = [X_{\alpha}, X_{-\alpha}]$ we obtain $\tau_{\alpha} = c_{b}(\tau_{\alpha}) = [c_{b}(X_{\alpha}), c_{b}(X_{-\alpha})] =
[X_{\alpha},  d_{- \alpha} X_{-\alpha}] = d_{- \alpha} \tau_{\alpha}$.  Hence $ d_{- \alpha} = 1$, which proves that c$_{b}$ is the identity on $\fG_{- \alpha}$.  This proves (i).

$\mathbf{Proof~ of~ (ii)}$  Since V is an irreducible $\fG$-module there is a unique highest dominant weight $\lambda \in \Lambda$.  We show that c$_{1} =$ e$^{\lambda(H)}$ satisfies the assertion of (ii), where a $=$ exp(H).

    If $\mu \in \Lambda$, then recall that $\mu = \lambda -$ \(\sum_{j=1}^{n} m_{j} \alpha_{j} \), where m$_{j} \in \bz^{+}$ for all j (cf. (20.2) of [Hu]).  By the definition of H we know that $\alpha_{k}$(H) $= N_{k} \pi i$ for suitable integers $N{_k} , 1 \leq k \leq n$.  Hence $\mu$(H) $= \lambda$(H) $ -$ \(\sum_{j=1}^{n} m_{j} \alpha_{j}(H) \) $= \lambda$(H) $- N \pi i$, where N $=$ \(\sum_{j=1}^{n} m_{j} N_{j} \) $\in \bz$.  Hence e$^{\mu(H)} = \pm$ e$^{\lambda(H)}$ for all $\mu \in \Lambda$.

    On V$_{\mu}$ we have a = e$^{\mu(H)}$ Id $= \pm$  e$^{\lambda(H)}$ Id since a = exp(H).  This proves (ii).

$\mathbf{Proof~ of~ (iii)}$  Every eigenspace of b in V is invariant under $\fG$ since b commutes with $\fG$ by (i).  Hence b $=$  c$_{2}$ Id on V for some complex number c$_{2}$ since V is an irreducible $\fG$-module.
\end{proof}

$\mathit{The~ Weyl~ group~ and~ real~ \fG_{0}-modules}$

    We now extend the previous proposition to the setting of real $\fG_{0}$-modules.  Let U be a real $\fG_{0}$-module, and let V $=$ U$^{\bc}$ be the corresponding $\fG = \fG_{0}^{\bc}$-module.  Let J : V $\rightarrow$ V denote the conjugation defined by U.

    For $\mu \in \Lambda$ we define U$_{\mu} =$ (V$_{\mu}~ \oplus~$V$_{-\mu}$) $\cap$ U, and we define U$_{0} =$ V$_{0}~ \cap~$U.  These " weight "  spaces for U are discussed in more detail in the next section.

\begin{proposition}  Let \{T$_{\sigma} : \sigma \in$ W\} be the elements of G$_{\bz}$ constructed in (5.1).  Let $\sigma \in$ W.  Then

\indent 1)  T$_{\sigma}$(U) $=$ U.

\indent 2)  T$_{\sigma}$(U$_{\mu}$) $=$ U$_{\sigma(\mu)}$ for all $\mu \in \Lambda$.

\indent 3)  T$_{\sigma}$ is an isometry of U with respect to any $\fG_{0}$-invariant inner product ( , ) on U.
\end{proposition}

\begin{proof}  1)  We shall need two preliminary results.

$\mathbf{Lemma~ A}$  Let $\fC$ = \{X$_{\beta}$ : $\beta \in \Phi$ ; $\tau_{\alpha}$ : $\alpha \in \Delta$ \} be a Chevalley basis for $\fG$ corresponding to a Cartan subalgebra $\fH$ and a choice of base $\Delta$ for the roots $\Phi$ of $\fH$.  Then the conjugation operator J : V $\rightarrow$ V has the following properties :

\indent 1)  J $\circ$ X$_{\alpha} \circ$ J$^{-1} = -$X$_{-\alpha}$ for all $\alpha \in \Phi$.

\indent 2)  J $\circ$ H $ \circ$ J$^{-1} = -$ H for all H $\in \fH_{\br}$.

\indent 3)  J normalizes $\fU(\fG)_{\bz}$ in End(V).

$\mathbf{Lemma~ B}$  Let $\langle , \rangle$ be a $\fG_{0}$-invariant Hermitian inner product on V.  Then J $\circ$ T$_{\sigma} \circ$ J$^{-1} =$ (T$_{\sigma}$*)$^{-1}$ for all $\sigma \in$ W.

    We assume Lemmas A and B for the moment and complete the proof of 1).  It suffices to prove 1) for the transformations T$_{\alpha} =$ T$_{\sigma_{\alpha}}$ since every T$_{\sigma}$ is a composition of transformations from \{T$_{\alpha} : \alpha \in \Phi$ \}.  Since U is the +1 eigenspace of the $\br$-linear operator J it suffices to prove that J commutes with T$_{\alpha}$ for all $\alpha \in \Phi$.  From Lemma B we obtain (J $\circ$ T$_{\alpha}$) (T$_{\alpha} \circ$ J)$^{-1} =$ J $\circ$ T$_{\alpha} \circ$ J$^{-1} \circ$ T$_{\alpha}^{-1} =$ (T$_{\alpha}$*)$^{-1}$(T$_{\alpha}$)$^{-1} =$ (T$_{\alpha} \circ$ T$_{\alpha}$*)$^{-1}$.  We know that T$_{\alpha} \circ$ T$_{\alpha}$* $=$ Id by the assertions 2) of Lemmas A and B in the proof of (5.1).  Hence J $\circ$ T$_{\alpha} =$ T$_{\alpha} \circ$ J for all $\alpha \in \Phi$.
\newline

$\mathit{Proof~ of~ Lemma~ A}$  1)  We recall from section 3 that the elements \{A$_{\alpha}$, B$_{\alpha} : \alpha \in \Phi^{+}$ \} of the compact Chevalley basis $\fC_{0}$ of $\fG_{0}$ are given by A$_{\alpha} =$ X$_{\alpha} -$ X$_{-\alpha}$ and B$_{\alpha} =$ i X$_{\alpha} +$ iX$_{-\alpha}$.  Hence
X$_{\alpha} = \frac{1}{2}$(A$_{\alpha} -$ i B$_{\alpha}$) and X$_{-\alpha} = - \frac{1}{2}$(A$_{\alpha} +$ i B$_{\alpha}$).  Note that J commutes with the elements of $\fG_{0}$ since $\fG_{0}$ leaves U invariant, and in particular J commutes with A$_{\alpha}$ and B$_{\alpha}$ for all $\alpha \in \Phi^{+}$.  The assertion of 1) now follows since J is conjugate linear on V.

    2)  For any $\alpha \in \Delta$ J commutes with $\tilde{\tau_{\alpha}} = i \tau_{\alpha}$ since $\tilde{\tau_{\alpha}} \in \fC_{0} \subset \fG_{0}$.  This, together with the conjugate linearity of J, shows that J $\circ~ \tau_{\alpha}~ \circ$ J$^{-1} = - \tau_{\alpha}$ for all $\alpha \in \Delta$.  The assertion 2) now follows since $\fH_{\br} = \br$-span \{$\tau_{\alpha} : \alpha \in \Delta$ \}.

    3)  This assertion follows from 1) since $\fU(\fG)_{\bz}$ is generated in End(V) by Id and \{(X$_{\alpha}$)$^{n} / n !$ \}.
\newline

$\mathit{Proof~ of~ Lemma~ B}$  Let $\sigma \in$ W be given.  The map C : g $\rightarrow$ (g*)$^{-1}$ is an isomorphism of GL(V).  By its definition in (5.1) T$_{\sigma}$ is a composition of transformations of the form T$_{\alpha} =$ exp(X$_{\alpha}$) $\circ$ exp($-$X$_{-\alpha}$) $\circ$ exp(X$_{\alpha}$), where $\alpha \in \Phi$.  Hence it suffices to consider the case that T$_{\sigma} =$ T$_{\alpha}$ for $\alpha \in \Phi$.

    Using properties of the matrix exponential map and Lemma A we compute J $\circ$ T$_{\alpha}~ \circ$  J$^{-1} =$ (J $\circ$ exp(X$_{\alpha}$) $\circ$ J$^{-1}$) $\circ$ (J $\circ$ exp($-$X$_{-\alpha}$) $\circ$ J$^{-1}$) $\circ$ (J $\circ$ exp(X$_{\alpha}$) $\circ$ J$^{-1}$) $=$ (exp(J $\circ$ X$_{\alpha} \circ$ J$^{-1}$)) $\circ$ (exp($-$J $\circ$ X$_{-\alpha} \circ$ J$^{-1}$)) $\circ$ (exp(J $\circ$ X$_{\alpha} \circ$ J$^{-1}$)) $=$ exp($-$X$_{-\alpha}$) $\circ$ exp(X$_{\alpha}$) $\circ$ exp($-$X$_{-\alpha}$).  A similar argument together with 2) of (4.3) shows that T$_{\alpha}$* $=$ exp(X$_{-\alpha}$) $\circ$ exp($-$X$_{\alpha}$) $\circ$ exp(X$_{-\alpha}$) $=$ (J $\circ$ T$_{\alpha} \circ$ J$^{-1}$)$^{-1}$, which completes the proof of Lemma B.

    We prove assertion 2) of the proposition.  By definition U$_{\mu} =$ (V$_{\mu}~ \oplus~$V$_{-\mu}$) $\cap$ U and U$_{0} =$ V$_{0}~ \cap~$U.  Assertion 2) is now an immediate consequence of 1) in (5.2) and 2) in (5.1).

    3)  Let $( , )$ be a $\fG_{0}$ - invariant inner product on U, and let $\langle , \rangle$ denote the extension of $( , )$ to a $\fG_{0}$ - invariant Hermitian inner product on V $=$ U$^{\bc}$.  By 3) of (5.1) each transformation T$_{\sigma}$ preserves  $\langle , \rangle$ and hence T$_{\sigma}$ preserves its restriction $( , )$ on U $=$ T$_{\sigma}$(U).
\end{proof}

$\mathit{The~ Weyl~ group~ and~ Aut(\fN)}$

    For each element $\sigma$ of W we now extend the transformation T$_{\sigma} \in$ GL(V) that was constructed in (5.1) to an automorphism of $\fN = V \oplus \fG$ with extra symmetry properties.  For the case that V $= \fG$, a G-module with respect to the homomorphism Ad : G $\rightarrow$ GL($\fG$), we define $\varphi_{\sigma} = T_{\sigma} :  \fG \rightarrow \fG$.  Let $\fG_{0}$ be a compact real form for $\fG$.

\begin{proposition}  For each element $\sigma$ of W let $\zeta_{\sigma} : \fN = V \oplus \fG \rightarrow \fN$ be the linear map given by $\zeta_{\sigma} = T_{\sigma}$ on V and $\zeta_{\sigma} = \varphi_{\sigma}$ on $\fG$.  Then $\zeta_{\sigma}$ is an automorphism of $\fN$ that is unitary with respect to any $\fG_{0}$ invariant Hermitian inner product on $\fN$ for which V and $\fG$ are orthogonal.
\end{proposition}

\begin{proof}  It suffices to prove the result in the case that $\sigma$ is one of the Weyl reflections $\sigma = \sigma_{\alpha}, \alpha \in \Phi$.  In this case we use the simpler notation $T_{\alpha}, \varphi_{\alpha}$ and $\zeta_{\alpha}$.

    Let $\ip$ be a $\fG_{0}$ invariant Hermitian inner product on $\fN = V \oplus \fG$ for which V and $\fG$ are orthogonal.  Then $T_{\alpha}$ is unitary on V by (5.1), and $\varphi_{\alpha}$ is unitary on $\fG$ since $\{\varphi_{\alpha}, \fG \}$ is a special case of $\{T_{\alpha}, V \}$.  Hence $\zeta_{\alpha}$ is unitary on $\fN = V \oplus \fG$.

    To prove that $\zeta_{\alpha}$ is an automorphism of $\fN$ we need a preliminary result similar to the one used in the proof of (4.10).

$\mathbf{Lemma}$

\indent 1)  $B^{*}(T_{\alpha}(v), T_{\alpha}(w)) = B^{*}(v,w)$ for all v,w $\in$ V.

\indent 2)   $B^{*}(\varphi_{\alpha}(X), \varphi_{\alpha}(Y)) = B^{*}(X,Y)$ for all X,Y $\in \fG$.
\newline

$\mathit{Proof~ of~ the~ Lemma}$  We note that 2) is a special case of 1) so it suffices to prove 1).  We define $\langle v,w \rangle = B^{*}(v,Jw)$ for v,w $\in$ V.  By a) of (2.2) we conclude that $\ip$ is a $\fG_{0}$ invariant Hermitian inner product on V.  The proof of (5.2) shows that T$_{\alpha}$ commutes with J : V $\rightarrow$ V.  Since $T_{\alpha} : V \rightarrow V$ is unitary with respect to $\ip$ by (5.1) we conclude that $B^{*}(v,w) = \langle v, Jw \rangle = \langle T_{\alpha}(v), T_{\alpha}(Jw) \rangle = \langle T_{\alpha}(v), JT_{\alpha}(w) \rangle = B^{*}(T_{\alpha}(v), T_{\alpha}(w))$.  This completes the proof of the Lemma.

    We now complete the proof of the Proposition, essentially in the same way as the completion of the proof of (4.10).   By the nondegeneracy of B* it suffices to show that $B^{*}(\varphi_{\alpha}([v,w]), X) = B^{*}([T_{\alpha}v, T_{\alpha}w], X)$ for all v,w $\in$ V and all X $\in \fG$.  The proof of 1) of (5.1) shows that $T_{\alpha}~ \circ~ Y~ \circ T_{\alpha}^{-1} = \varphi_{\alpha}(Y)$ for all Y $\in \fG$, or equivalently that $\varphi_{\alpha}^{-1}(X) = T_{\alpha}^{-1}~ \circ~ X~ \circ~ T_{\alpha}$ for all X $\in \fG$.  Using both parts of the lemma above we compute $B^{*}([T_{\alpha}v, T_{\alpha}w], X) = B^{*}(X(T_{\alpha}v), T_{\alpha}w) = B^{*}((T_{\alpha}^{-1}~ \circ~ X~ \circ~ T_{\alpha})(v), w) = B^{*}((\varphi_{\alpha}^{-1}(X))(v),w) = B^{*}([v,w], \varphi_{\alpha}^{-1}(X)) = B^{*}(\varphi_{\alpha}([v,w]), X)$.

\end{proof}

\section{Weight~ space~ decomposition~ of~ a~ real~ $\fG_{0}$~ module}

    Let $\fG_0$ be a real, compact, semisimple Lie algebra.   In this section we let U be a finite dimensional real $\fG_{0}$-module, and we set V = U$^\bc$, the corresponding finite dimensional complex $\fG$ module, where $\fG = \fG_{0}^\bc$.  Let $\fH_0$ be a maximal abelian subalgebra of $\fG_{0}$, and let $\fH = \fH_{0}^\bc$ be the corresponding Cartan subalgebra of $\fG$.  We show how the weight space decomposition \ws determined by $\fH$ induces a (weight space) decomposition U = U$_0$ + \(\sum_{\lambda \epsilon \Lambda^{+}}\) U$_\lambda$, where $\Lambda^{+}$ is any subset of
$\Lambda$ such that $\Lambda$ is the disjoint union of $\Lambda^{+}$ and $-\Lambda^{+}$.  We also show that if $\ip$ is any inner product on U for which the elements of $\fG_0 \subset$ End(U) are skew symmetric, then the decomposition of U is an orthogonal direct sum decomposition.

\begin{proposition}  Let $\fH_0$ be a maximal abelian subalgebra of $\fG_{0}$, and let $\fH = \fH_{0}^{\bc}$ be the corresponding Cartan subalgebra of $\fG = \fG_{0}^{\bc}$.  Let $\fH_{\br} = $\{H $\in \fH : \alpha(H) \in \br$ for all $\alpha \in \Phi \}$.  Then i $\fH_0  = \fH_{\br} =  \br$-span $\{\tau_{\alpha}: \alpha \in \Delta\}$
\end{proposition}

\begin{proof}  Recall that $\{\tau_{\alpha}: \alpha \in \Delta \}$ is a $\bc$ - basis of $\fH$.  Since $\alpha(\tau_{\beta}) \in \bz$ for all $\alpha, \beta \in \Phi$ it follows that $\br$-span $\{\tau_{\alpha}: \alpha \in \Delta \} \subseteq \fH_{\br}$.  To prove the reverse inclusion let H = \(\sum_{i=1}^{n} b_{i} \tau_{\alpha_{i}} \in \fH_{\br}$, where $b_{i} \in \bc$ and $\Delta = \{\alpha_{1},~ ...~,\alpha_{n} \}$.  If $A =   (\alpha_{1}(H),~ ...~, \alpha_{n}(H)) , B = (b_{1},~ ...~, b_{n})$ and C is the Cartan matrix given by $C_{ij} = \alpha_{i}(\tau_{\alpha_{j}})$, then $A = BC^{t}$.  Since C is invertible with integer entries and A by hypothesis has real entries it follows that $B = A(C^{t})^{-1}$ has real entries.  Hence $\fH_{\br} \subseteq \br$-span $\{\tau_{\alpha}: \alpha \in \Delta \}$.  Finally, $\alpha(\fH_{0}) \subset i\br$ for all $\alpha \in \Phi$ by the second remark preceding (2.4).  It follows that $i\fH_{0} \subseteq \fH_{\br} = \br$-span $\{\tau_{\alpha}: \alpha \in \Delta \}$ and equality holds since $(i\fH_{0})^{\bc} = \fH_{\br}^{\bc} = \fH$.
\end{proof}

$\mathit{Remark}$  Note that $\mu(H) \in \br$ for all H $\in \fH_{\br}$ and all $\mu \in \Lambda$.  This follows from (6.1) and the fact that $\mu(\tau_{\alpha})$ is an integer for all $\mu \in \Lambda$ and $\alpha \in \Phi$.
\newline

$\mathit{Conjugations}$

\begin{proposition}  Let $\fH_{0}, \fG_{0},$U and V be as above.  Let J : V $\rightarrow$ V and J$_{0}$ : $\fG_{0} \rightarrow \fG_{0}$ denote the conjugations determined by U and $\fG_{0}$ respectively.  Then

    1)     If $\lambda \in \Lambda$, then $-\lambda \in \Lambda$ and J(V$_{\lambda}$) =
V$_{-\lambda}$.  In particular, dim V$_{\lambda}$ = dim V$_{-\lambda}$.  Moreover, J(V$_0$) = V$_0$.

    2)  $J_{0}(\fG_{\beta}) = \fG_{- \beta}$

    3)  $J \circ X = J_{0}(X) \circ J$ for all $X \in \fG$.
\end{proposition}

\begin{proof}  The proof of 4) of (4.5) contains a proof of 1).  Assertion 2) is assertion 1) applied to the adjoint representation of $\fG$.  Assertion 3) is a routine computation.
\end{proof}

$\mathit{Weight~ space~ decomposition~ of~ U}$

    For a real finite dimensional $\fG_0$-module U we construct the analogue of the weight space decomposition of the complex finite dimensional $\fG$-module V = U$^{\bc}$.

\begin{proposition} For each $\lambda \in \Lambda$ let U$_{\lambda}$ = (V$_{\lambda} \oplus$ V$_{-\lambda}$) $\cap$ U.  Let U$_0$ =  V$_0~ \cap$ U.  Then

\indent 1)  If H$_0 \in \fH_0$, then H$_0$(U$_{\lambda}$) $\subseteq$U$_{\lambda}$.

\indent 2) U$_{\lambda}$ = Re(V$_{\lambda}$) = Im(V$_{\lambda}$) = Re(V$_{-\lambda}$) = Im(V$_{-\lambda}$).
U$_0$ = Re(V$_0$) = Im(V$_0$).

\indent 3)  U$_{\lambda}^{\bc}$ = V$_{\lambda} \oplus$ V$_{-\lambda}$.  U$_0^{\bc}$= V$_0$.

\indent 4)  dim$_{\br}$U$_{\lambda}$ = 2 dim$_{\bc}$V$_{\lambda}$

\indent 5)  Define $\varphi : V_{\lambda} \rightarrow U_{\lambda}$ by $\varphi(v_{\lambda}) = v_{\lambda} + J(v_{\lambda})$ for all v$_{\lambda} \in$ V$_{\lambda}$.  Then $\varphi$ is an $\br$ - linear isomorphism, where we regard $V_{\lambda}$ as a real vector space.

\indent 6)  Let $\Lambda_0 =$ i $\Lambda$.  Then

\indent \indent     a)  $\Lambda_0 \subset$ Hom ($\fH_0, \br$)

\indent \indent     b)  Let $\lambda \in \Lambda$ and let $\lambda_0 =$ i $\lambda \in \Lambda_0$. Then

\hspace{.26in} $U_{\lambda} = \{u \in U : H_0^{2}(u) = -~\lambda_{0}(H_{0})^{2}u$ \hspace{.05in}  for all $H_{0} \in \fH_{0}\}$.

\end{proposition}

\begin{proof}  1) Clearly H$_0$(U) $\subseteq$ U and since $\fH_0 \subset \fH$ it follows that H$_0$(V$_{\lambda}$) $\subseteq$ V$_{\lambda}$ and H$_0$(V$_{-\lambda}$) $\subseteq$ V$_{-\lambda}$.  The assertion now follows from the definition of U$_{\lambda}$.

\indent 2)  If v $\in$ V$_{\lambda}$, then Jv and i Jv are elements of V$_{-\lambda} =$ J(V$_{\lambda}$) by (6.2).  Since Re(v) = Re(Jv) = Im(i Jv) it follows that Re(V$_{\lambda}$) = Re(V$_{-\lambda}$) = Im(V$_{\lambda}$) = Im(V$_{-\lambda}$).

    If v $\in$ V$_{\lambda}$, then Re(v) = (1/2)(v + Jv) $\in$ (V$_{\lambda} \oplus$ V$_{-\lambda}$) $\cap$ U = U$_{\lambda}$
by (6.2).  If u = v$_1$ + v$_2$ $\in$ U$_{\lambda}$ , where v$_1 \in$ V$_{\lambda}$ and v$_2 \in$ V$_{-\lambda}$
= J(V$_{\lambda}$), then  u = Re(u) = Re(v$_1$ + Jv$_2) \in$ Re(V$_{\lambda}$).  This shows that U$_{\lambda}$ = Re(V$_{\lambda}$).
 Since J(V$_0$) = V$_0$ a similar argument shows that U$_0$ = Re(V$_0$) = Im(V$_0$).

\indent 3)  Observe that  U$_{\lambda}^{\bc} \subseteq$ V$_{\lambda} \oplus$ V$_{-\lambda}$ since V$_{\lambda} \oplus$ V$_{-\lambda}$ is a complex subspace of V that contains U$_{\lambda}$.  Note that V$_{\lambda} \subset$ U$_{\lambda}^{\bc}$ since U$_{\lambda}$ = Re(V$_{\lambda}$) = Im(V$_{\lambda}$).  Similarly V$_{-\lambda} \subset$ U$_{\lambda}^{\bc}$ since U$_{\lambda}$ =
Re(V$_{-\lambda}$) = Im(V$_{-\lambda}$) , and it follows that V$_{\lambda} \oplus$ V$_{-\lambda} \subseteq$ U$_{\lambda}^{\bc}$ . A similar argument shows that U$_0^{\bc}$ = V$_0$.

\indent 4)  This follows immediately from 3) and (6.2) since  dim$_{\br}$U$_{\lambda}$ = dim$_{\bc}$U$_{\lambda}^{\bc}$.

\indent 5)  Clearly $\varphi(V_{\lambda}) \subseteq U_{\lambda}$ since J fixes each element $\varphi(v_{\lambda}), v_{\lambda} \in V_{\lambda}$.  If $\varphi(v_{\lambda}) = 0$ for some $v_{\lambda} \in V_{\lambda}$, then $J(v_{\lambda}) = - v_{\lambda}$, which implies that $v_{\lambda} \in V_{\lambda} \cap J(V_{\lambda}) = V_{\lambda}~ \cap~ V_{- \lambda} = \{0\}$.  Hence $\varphi$ is injective.

    It remains only to show that dim$_{\br}U_{\lambda} =$ dim $_{\br}V_{\lambda} = 2$dim $_{\bc}V_{\lambda}$, but this is 4).

\indent 6)  a)  This is an immediate consequence of the fact that $\lambda$($\fH_0$) $\subset$ i $\br$ for all $\lambda \in \Lambda$  by (6.1) and the remark that follows it.

\indent b)  If $\lambda_0 =$ i $\lambda$ for some $\lambda \in \Lambda$, then from the definition of U$_{\lambda}$ it is clear that U$_{\lambda} \subseteq$ U$'_{\lambda}$= $\{\mathrm{u} \in \mathrm{U : H}_0^{2}\mathrm{(u)} = -~ \lambda_{0}$ ({H}$_0)^{2}~ \mathrm{u~ for~ all~ H}_0 \in \fH_0\}$ since $\fH_0 \subset \fH = \fH_{0}^{\bc}$.  To prove the reverse inclusion we use the following

$\mathbf{Lemma}$  There exists an element H$_{0}$ of $\fH_{0}$ such that a) $\lambda_{0}(H_{0}) \neq 0$ for all $\lambda_{0} \in \Lambda_{0}$ \hspace{.2in} b)  if $\sigma_{0}(H_{0})^{2} = \mu_{0}(H_{0})^{2}$ for elements $\sigma_{0}, \mu_{0}$ of $\fH_{0}$, then $\sigma_{0} = \mu_{0}$ or $\sigma_{0} = - \mu_{0}$
\begin{proof}  Let S$_{1} = \{\sigma_{0} - \mu_{0} : \sigma_{0} \neq\mu_{0} \in \Lambda_{0} \}$ and let S$_{2} = \{\sigma_{0} + \mu_{0} : \sigma_{0} \neq - \mu_{0} \in \Lambda_{0} \}$.  If W$_{0}$ is the union of the kernels of elements in $S_{1}~ \cup~ S_{2}~ \cup~ \Lambda_{0}$, then W$_{0}$ is a proper subset of $\fH_{0}$ since the kernels of nonzero elements of Hom($\fH_{0},\bc$) are proper subspaces of $\fH_{0}$.  If H$_{0} \in \fH_{0} - W_{0}$, then H$_{0}$ satisfies the conditions of the lemma.
\end{proof}
    We now complete the proof of b).  Let u be any element of U$'_{\lambda}$and write u = v$_0$ + \(\sum_{\mu \epsilon \Lambda}\) v$_\mu$, where v$_0 \in$ V$_0$ and v$_\mu \in$ V$_\mu$ for all $\mu \in \Lambda$.  Then $-~ \lambda_{0}(H_{0})^{2}v_{0} +$ \(\sum_{\mu \epsilon \Lambda}\) $-~ \lambda_{0}(H_{0})^{2} v_{\mu} =
(H_{0})^{2}(u) =$ \(\sum_{\mu \epsilon \Lambda}\) $\mu(H_{0})^{2} v{_\mu} =$
\(\sum_{\mu \epsilon \Lambda}\)$ -~\mu_{0}(H_{0})^{2} v{_\mu}$, where $\mu_{0} =$ i $\mu$.  Hence v$_0$ = 0 and $-~ \lambda_{0}(H_{0})^{2} = -~ \mu_{0}(H_{0})^{2}$ whenever v$_{\mu}$ is nonzero.  By the choice of H$_0 \in \fH_0$ it follows that v$_{\mu}$ = 0 if $\mu \notin \{\lambda,-~ \lambda\}$.  Hence u = v$_{\lambda}$ + v$_{-\lambda} \in$ (V$_{\lambda} \oplus$ V$_{-\lambda}$) $\cap$ U = U$_{\lambda}$ , which proves that U$'_{\lambda} \subseteq$ U$_{\lambda}$.
\end{proof}

\begin{proposition}  (Weight space decomposition of U)  Let $\Lambda^{+}$ be any subset of $\Lambda$ such
that $\Lambda$ is the disjoint union of $\Lambda^{+}$ and $-\Lambda^{+}$.   Then

\indent \indent     (*) U = U$_0$ + \(\sum_{\lambda \epsilon \Lambda^{+}}\) U$_\lambda$  (direct sum)

If $\ip$ is any inner product on U such that the elements of $\fG_0 \subset$ End(U) are skew symmetric, then the
decomposition of U above is an orthogonal direct sum decomposition.
\end{proposition}

\begin{proof}  Let U$'$ = U$_0$ + \(\sum_{\lambda \epsilon \Lambda^{+}}\) U$_\lambda$.  This sum is direct by the definitions of U$_0$ and U$_\lambda$ and the fact that the weight space decomposition \ws~is a direct sum.  Since U$_{\lambda}$ = U$_{-\lambda}$ for all $\lambda \in \Lambda$ it follows from 3) of (6.3) that  V = (U$'$)$^{\bc} \subseteq$ U$^{\bc}$ = V, which proves (*).

    Now let $\ip$ be any inner product on U such that the elements of $\fG_0 \subset$ End(U) are skew symmetric.  Choose a nonzero element H$_0$ of $\fH_0$ such that if $\sigma_{0}$(H$_0$)$^2$ = $\mu_{0}$(H$_0$)$^2$ for elements $\sigma_{0}$, $\mu_{0}$ of $\Lambda_{0}$, then $\sigma_{0}$ = $\mu_{0}$ or $\sigma = -\mu_{0}$.  The element H$_0^{2}$ is symmetric and negative semidefinite relative to $\ip$ since H$_0$ is skew symmetric.  By the choice of H$_0$ and 6 b) of (6.3) the  vector spaces U$_0$ and U$_{\lambda}~, ~\lambda \in \Lambda^{+}$, are distinct eigenspaces of H$_0^{2}$ and must therefore be orthogonal.
\end{proof}

$\mathit{Root~ space~ decomposition~ of~ \fG_0}$

    If U = $\fG_0$ , then U is a real $\fG_0$-module with respect to the adjoint representation : X(Y) =
\linebreak ad X(Y) = [X, Y] for all X,Y $\in \fG_0$.  In this case V = U$^{\bc}$ = $\fG_0^{\bc}$ = $\fG$  describes the adjoint representation for $\fG$.

    Let $\fH_0$ be a maximal abelian subalgebra of $\fG_0$, and let $\fH$ = $\fH_0^{\bc}$ be the corresponding Cartan subalgebra of $\fG$ = $\fG_0^{\bc}$.  In this case the weights $\Lambda$ are identical to the roots $\Phi$, and the weight = root spaces $\fG_{\beta}$ are 1-dimensional complex vector spaces.  By (6.2)  ~J$_0$($\fG_{\beta}$) = $\fG_{-\beta}$ for each $\beta \in \Phi$, where J$_0$ denotes the conjugation in $\fG$ determined by $\fG_0$.  If $\fG_{0,\beta}$ = ($\fG_{\beta} \oplus \fG_{-\beta}$) $\cap~ \fG_0$, then by 3) and 4) of (6.3)  each $\fG_{0,\beta}~,~\beta \in \Phi$, is a 2-dimensional subspace of $\fG_0$, and $\fG_{0,\beta}^{\bc} = \fG_{\beta} \oplus \fG_{-\beta}$.  Note that $\fG_{0,\beta} = \fG_{0,-\beta}$ for all $\beta \in \Phi$ so we only need to consider $\beta \in \Phi^{+}$.

    The zero weight space for the adjoint representation of $\fG$ is the Cartan subalgebra
$\fH$, and $\fH_0 = \fH \cap \fG_0$.  Hence from the root space decomposition \rs we obtain

\begin{proposition}  (Rootspace decomposition of $\fG_{0}$)  Let $\fG_{0}$ be a compact, semisimple real Lie algebra, and let $\fH_{0}$ be a maximal abelian subspace of $\fG_{0}$. Then

\indent \indent $\fG_{0}$ = $\fH_{0}$ + \(\sum_{\alpha \epsilon \Phi^{+}}$ $\fG_{0,\alpha}$ \hspace{.5in} (direct sum)

where $\fG_{0,\alpha}$ = ($\fG_{\alpha} \oplus \fG_{-\alpha}$) $\cap~ \fG_0$, $\fG = \fG_0^{\bc}$  and
$\fH = \fH_0^{\bc}$.  Each subspace $\fG_{0,\alpha}$ is 2-dimensional.
\end{proposition}

    If $\ip = -B_{0}$ , then $\ip$ is an inner product on $\fG_{0}$ since the Killing form $B_0$ is negative definiteon $\fG_0$.  Moreover, the elements ad X, X $\in \fG_0$, are skew symmetric elements of End($\fG_0$) by properties ofthe Killing form.  Hence by (6.4) the root space decomposition of $\fG_0$ is orthogonal relative to $\ip = -B_{0}$.

    Our next observation is the real analogue of the fact that $\fG_{\beta}$(V$_{\lambda}$) $\subseteq$
V$_{\lambda+\beta}$ for all complex $\fG$-modules V, where $\lambda \in \Lambda$ and $\beta \in \Phi$.

\begin{proposition}  Let  (*)   U = U$_0$ + \(\sum_{\lambda \epsilon \Lambda^{+}}\) U$_\lambda$
(direct sum)

be the weight space decomposition of U.  Let $\fG_{0,\beta}$ = ($\fG_{\beta} \oplus \fG_{-\beta}$) $\cap~ \fG_0$
for all $\beta \in \Phi^{+}$.  Then $\fG_{0,\beta}(U_{\lambda}) \subseteq U_{\lambda+\beta} \oplus
U_{\lambda-\beta}$ for all $\lambda \in \Lambda~,~ \beta \in \Phi^{+}$.
\end{proposition}

$\mathit{Remark}$  Our convention is that U$_{\mu}$= \{0\}  if $\mu \in$ Hom($\fH$, $\bc$) $-~ \Lambda$.

\begin{proof}   Let X $\in \fG_{0,\beta}$, and $u_{\lambda} \in U_{\lambda}$ be given for $\beta \in \Phi^{+}$
and $\lambda \in \Lambda$.  Write X = X$_{\beta}$ + X$_{-\beta}$ and u$_{\lambda}$ = v$_{\lambda}$ + v$_{-\lambda}$,
where X$_{\beta} \in$ $\fG_{\beta}$, X$_{-\beta} \in \fG_{-\beta}$, v$_{\lambda} \in$ V$_{\lambda}$ and
v$_{-\lambda} \in$ V$_{-\lambda}$.  Then X(u$_{\lambda}$) = $\xi_{1}$ + $\xi_{2}$, where $\xi_{1}$ =
X$_{\beta}$(v$_{\lambda}$) + X$_{-\beta}$(v$_{-\lambda}$) $\in$ V$_1$ = V$_{\lambda+\beta}~ \oplus$ V$_{-\lambda-\beta}$
and $\xi_{2}$ = X$_{\beta}$(v$_{-\lambda}$) + X$_{-\beta}$(v$_{\lambda}$) $\in$ V$_2$ = V$_{\lambda-\beta}~ \oplus$~
V$_{\beta-\lambda}$.   By (6.2) the spaces V$_1$ and V$_2$ are invariant under the conjugation J
of V induced by U.  Since X(u$_{\lambda}$) $\in$ U we have $\xi_{1} + \xi_{2}$ = X(u$_{\lambda}$) = J(X(u$_{\lambda}$))
 = J($\xi_{1}$) + J($\xi_{2}$), where J($\xi_{1}$) $\in$ V$_1$ and J($\xi_{2}$) $\in$ V$_2$.
 Since V$_1~ \cap$ V$_2 =$ \{0\} it follows that $\xi_{1}$ = J($\xi_{1}$) and $\xi_{2} =$ J($\xi_{2}$), which implies
that $\xi_{1} \in$ V$_1~ \cap$ U = U$_{\lambda+\beta}$ and $\xi_{2} \in$ V$_2~ \cap$ U = U$_{\lambda-\beta}$.
Hence X(u$_{\lambda}$) $\in$ U$_{\lambda+\beta}~ \oplus$ U$_{\lambda-\beta}$ for all X $\in \fG_{0,\beta}$,
and all u$_{\lambda} \in$ U$_{\lambda}$.
\end{proof}

$\mathit{Bracket~ relations~ in~ \fN_{0} = U \oplus \fG_{0}}$

\begin{proposition}   Let $\fG_0$ be a compact, semisimple real Lie algebra, and let $\fH_0$ be a maximal abelian subalgebra of $\fG_0$.
Let U be a real $\fG_0$-module, and let $\{\fN_{0} = U \oplus \fG_{0}, \langle~ ,~ \rangle \}$ be the real 2-step nilpotent Lie algebra defined in (2.4).  Let V = U$^{\bc}$ , $\fH = \fH_0^{\bc}$ and $\fG = \fG_0^{\bc}$ and let $\Lambda \subset$ Hom($\fH,\bc$) be the set of weights for V determined by $\fH$.  Let

\indent (*) U = U$_0$ + \(\sum_{\lambda \epsilon \Lambda^{+}}\) U$_\lambda$  (direct sum)

be the corresponding weight space decomposition of U (cf. 6.4).  Then

\indent 1)  [U$_0$ , U$_{\lambda}$] = $\fG_{0,\lambda}$ \hspace{1.3in} for all $\lambda \in \Lambda$

\indent 2)  [U$_{\mu}$ , U$_{\lambda}$] = $\fG_{0,\mu+\lambda} \oplus \fG_{0,\mu-\lambda}$ \hspace{.52in} for all $\mu~, \lambda$ with
$\mu~+~\lambda \neq$ 0

\indent 3)   [U$_0$ , U$_0$] = \{0\}
\end{proposition}

$\mathit{Remark}$  Our convention is that $\fG_{0,\lambda}$ = \{0\}  if $\lambda \notin \Phi$.

\begin{proof}

    We prove 2) but not 1) since the proof of 1) is almost identical.  Let $\mu~,~\lambda \in \Lambda$  be given.  From the definitions and (4.7) we have [U$_{\mu}$, U$_{\lambda}$] $=$
[(V$_{\mu}~ \oplus$ V$_{-\mu}$) $\cap$ U , (V$_{\lambda}~ \oplus$ V$_{-\lambda}$) $\cap$ U] $\subseteq \fP \cap \fG_0$, where $\fP$ =
($\fG_{\mu+\lambda} \oplus \fG_{-\mu-\lambda} \oplus \fG_{\mu-\lambda} \oplus \fG_{\lambda-\mu}$).  By (6.5) $\fP_0 := \fG_{0,\mu+\lambda} \oplus \fG_{0,\mu-\lambda} \subseteq \fP \cap \fG_0$ , and equality holds since $(\fP_{0})^{\bc} \subseteq (\fP \cap \fG_{0})^{\bc}  \subseteq \fP = (\fP_{0})^{\bc}$ by the discussion preceding (6.5).  This proves that [U$_{\lambda}$ , U$_{\mu}$] $\subseteq\fP \cap \fG_{0} = \fP_{0} = \fG_{0,\mu+\lambda} \oplus \fG_{0,\mu-\lambda}$

    To prove equality in the inclusion above it suffices to show that $[U_{\lambda} , U_{\mu}]^{\bc} = \fP = (\fP_{0})^{\bc}$.  Using
(4.7) and (6.3) we compute $[U_{\lambda} , U_{\mu}]^{\bc} = [U_{\mu}^{\bc} , U_{\lambda}^{\bc}] =
[V_{\mu}~ \oplus V_{-\mu} , V_{\lambda}~ \oplus V_{-\lambda}] =
[V_{\mu} , V_{\lambda}]~ \oplus [V_{-\mu} , V_{\lambda}] \oplus [V_{\mu} , V_{-\lambda}] \oplus~ [V_{-\mu} , V_{-\lambda}] ~=
\fG_{\mu+\lambda}~ \oplus~ \fG_{-\mu-\lambda}~ \oplus ~\fG_{\mu-\lambda}~ \oplus~ \fG_{\lambda-\mu} = \fP$.

\indent 3)  Let $\fG = \fG_0^{\bc}$ and V = U$^{\bc}$ and consider $\fN^{\bc}$ = V$\oplus \fG$.  Since U$_0$ = U $\cap$ V$_0$ it follows
 from (6.3) and 4) of (4.7) that [U$_0$ , U$_0$] $\subset$ [V$_0$ , V$_0$] = \{0\}.
 \newline
 \end{proof}

$\mathit{Abstract~ weights~ and~ real~ weight~ vectors}$

        We let $\tilde{\Lambda}$ denote the set of $\mathit{abstract~ weights}$ determined by $\fH$; that is, $\tilde{\Lambda} =$ \{$\lambda \in$ Hom$(\fH,\bc) : \lambda(\tau_{\alpha}) \in \bz$ for all $\alpha \in \Phi$\}.  Let $\fN$ = $\{$U$ \oplus \fG_{0} , \ip \}$ be as in (6.7).   For $\lambda \in \tilde{\Lambda}$ we deduce from (6.1) and the remark that follows it that i $\lambda \in$ Hom($\fH_0,\br$).   Let $\tilde{H}_{\lambda}$ be the unique vector in $\fH_0$ such that $\langle \tilde{H}, \tilde{H}_{\lambda} \rangle$ = $-$ i $\lambda$($\tilde{H}$) for all $\tilde{H} \in \fH_0$.  We call $\tilde{H}_{\lambda}$ the $\mathit{real~ weight~ vector}$ determined by $\lambda$ and $\ip$.
\newline

$\mathit{Relationship~ to~ complex~ weight~ vectors}$

    Let $\fG = \fG_0^{\bc}$ and $\fH~ = ~\fH_0^{\bc}$ as in (6.7), and let B denote the Killing form of $\fG$.  For $\lambda \in$ Hom($\fH,\bc$) let H$_{\lambda} \in \fH$ be the unique vector such that B(H, H$_{\lambda}$) = $\lambda$(H) for all H $\in \fH$.  We call H$_{\lambda}$ the $\mathit{complex~ weight~ vector}$ determined by $\lambda$.  The existence of H$_{\lambda}$ follows from the nondegeneracy of B on H.

\begin{proposition}  Let $\ip = - B_{0}$ on $\fG_{0}$, where $B_{0}$ denotes the Killing form of $\fG_0$.  For $\lambda \in \tilde{\Lambda}$ let $\tilde{H}_{\lambda} \in \fH_0$ and H$_{\lambda} \in \fH$ be the weight vectors defined above.  Then $\tilde{H}_{\lambda}$ = i H$_{\lambda}$.
\end{proposition}

\begin{proof}   Note that  $B(\tau_{\alpha}, H_{\lambda}) = \lambda(\tau_{\alpha})$ is an integer for all $\alpha \in \Phi$, which implies that $\alpha$(H$_{\lambda}$) = B(H$_{\alpha}$, H$_{\lambda}$) $\in \br$ since H$_{\alpha}$ is a real multiple of $\tau_{\alpha}$.  Hence  $H_{\lambda} \in \fH_{\br} = \{H \in \fH : \alpha(H) \in \br$ for all $\alpha \in \Phi \}$. It follows that $i H_{\lambda} \in i \fH_{\br} = \fH_0$ by (6.1).  Since $B = B_{0}$ on $\fG_0$ we obtain, for every $\tilde{H} \in \fH_0,
\langle \tilde{H} , i H_{\lambda} \rangle = - B_{0}(\tilde{H} , i H_{\lambda}) = - B(\tilde{H} , i H_{\lambda}) =~ -~
i~\lambda(\tilde{H})~ =~ \langle \tilde{H}~ , ~\tilde{H}_{\lambda} \rangle$.
\end{proof}

$\mathit{Rationality~ of~ weight~ vectors}$

    If $\tilde{\tau}_{\alpha}~ = ~ $i $\tau_{\alpha}$  for $\alpha \in \Phi$, then $\tilde{\tau}_{\alpha} \in$ i $\fH_{\br}~ =~ \fH_0$ by (6.1), and $\fH_0~ =~ \br$-span \{$\tilde{\tau}_{\alpha}$ : $\alpha \in \Delta \}$ since $\fH_{\br} = \br$-span $\{\tau_{\alpha} : \alpha \in \Delta \}$.  The next result is useful.

\begin{proposition}  Let $\ip = -$ B$_0$ on $\fG_{0}$ , where B$_0$ denotes the Killing form of $\fG_0$.  For any $\lambda \in \Lambda$ the real weight vector $\tilde{H}_{\lambda}$ lies
in $\bq$-span \{$\tilde{\tau}_{\alpha}$ : $\alpha \in \Delta \}$.
\end{proposition}
\begin{proof}

    We shall need the following observation whose proof we leave as an exercise

\begin{lem} Let $\{$U,~ $\ip$ \} be a finite dimensional real inner product space, and let $\fB$ be a basis of U
such that $\langle \xi , \xi' \rangle \in \bq$ for all $\xi~,~\xi' \in \fB$.  Let u be a vector of U such that
$\langle$ u, $\xi \rangle \in \bq$ for all $\xi \in \fB$.  Then u $\in \bq$-span($\fB$).
\end{lem}

    We complete the proof of the Proposition.  Note that for $\alpha,\beta \in \Phi$ we have ad $\tau_{\alpha}$ =$\beta(\tau_{\alpha})$ Id on $\fG_{\beta}$, and $\beta(\tau_{\alpha}) \in \bz$.  Hence
$\langle \tilde{\tau}_{\alpha}~ ,~ \tilde{\tau}_{\beta} \rangle~ =~ -$ B$_0$(i $\tau_{\alpha}$ , i $\tau_{\beta}$) = $-$ B(i $\tau_{\alpha}$ , i $\tau_{\beta}$) = B($\tau_{\alpha}$ , $\tau_{\beta}$) = trace (ad $\tau_{\alpha}~ \circ$ ad $\tau_{\beta}$) = \(\sum_{\lambda \epsilon \Phi} \lambda(\tau_{\alpha}) \lambda(\tau_{\beta}) \) $\in \bz$.   If $\fB$ = \{$\tilde{\tau}_{\alpha}$ : $\alpha \in \Delta \}$, then $\fB$ is a basis for $\fH_{\br} = \br$-span \{$\tilde{\tau}_{\alpha}$ : $\alpha \in \Delta \}$. For $\alpha \in \Delta~,~\lambda \in \Lambda$we obtain $\langle \tilde{\tau}_{\alpha}~ ,~ \tilde{H}_{\lambda} \rangle$  = ($-$ i $\lambda$)(i $\tau_{\alpha}$) = $\lambda$($\tau_{\alpha})$ $\in \bz$.  Hence $\tilde{H}_{\lambda} \in \bq$-span \{$\tilde{\tau}_{\alpha}$ : $\alpha \in \Delta \}$ by the lemma above.
 \end{proof}

\begin{proposition}  Let \{$\fN_{0} = U \oplus \fG_{0}, \langle~ ,~ \rangle$ \} be as in (6.7).  Then

\indent 1) $[U_{\lambda},U_{\lambda}] = \br$-span $\{\tilde{H}_{\lambda}\}$ \hspace{.5in}   if $\lambda \in \Lambda~,~ 2\lambda \notin \Phi$.

\indent 2) $[U_{\beta},U_{\beta}] = \br$-span $\{\tilde{H}_{\beta}\}$ \hspace{.5in} if $\beta \in \Lambda \cap \Phi$
\end{proposition}

\begin{proof}  1)  If $\lambda \in \Lambda$, and 2$\lambda \notin \Phi$ , then [U$_{\lambda}$,U$_{\lambda}$] $\subseteq \fG_{0,0}~
 =~ \fH \cap \fG_0~ =~ \fH_0$ by (6.5) and 2) of (6.7).  Next, we show that $[U_{\lambda},U_{\lambda}]~\neq \{0\}$.
 By 1) and 2) of (4.7), 3) of (6.3) and the hypothesis that 2$\lambda \notin \Phi$ ( hence $-~ 2\lambda \notin \Phi$)
 we obtain $[U_{\lambda},U_{\lambda}]^{\bc} = [U_{\lambda}^{\bc}, U_{\lambda}^{\bc}] =
 [V_{\lambda} \oplus V_{-\lambda}~,~V_{\lambda} \oplus V_{-\lambda}] = [V_{\lambda}, V_{-\lambda}] \neq$ \{0\}.

    Let u$_{\lambda}$ and u$_{\lambda}'$ be arbitrary elements of U$_{\lambda}$, and let $\tilde{H}$ be an element of $\fH_0$
such that $\langle \tilde{H}~ ,~ \tilde{H}_{\lambda} \rangle~ =~ 0$.  To finish the proof of 1) it suffices to show that
$\langle [u_{\lambda}~, u_{\lambda}'] ,  \tilde{H} \rangle~ =~ 0$.

    We show first that $\tilde{H}$(u$_{\lambda}$)  = 0.  Note that $\lambda$($\tilde{H}$) = 0 since
$- i \lambda(\tilde{H}) = \langle \tilde{H}~ ,~ \tilde{H_{\lambda}}\rangle = 0$.  By the definition of U$_{\lambda}$
in (6.3) there exist elements v$_{\lambda} \in$ V$_{\lambda}$ and v$_{-\lambda} \in$ V$_{-\lambda}$ such that
u$_{\lambda}$ = v$_{\lambda}$ + v$_{-\lambda}$.  Hence $\tilde{H}$(u$_{\lambda}$) $=
\lambda$($\tilde{H}$)v$_{\lambda}~-~ \lambda$($\tilde{H}$)v$_{-\lambda}~=~ 0$.

    Finally, $\langle [u_{\lambda}~, u_{\lambda}'] ,  \tilde{H} \rangle~ =
~ \langle \tilde{H}(u_{\lambda})~,u_{\lambda}' \rangle~ =~ 0$.

\indent 2)  If $\beta \in \Lambda \cap \Phi$ , then 2$\beta \notin \Phi$  by a basic property of root systems.
The assertion now follows immediately from 1)
\end{proof}

    Next we obtain the analogues for $\fN_{0} = U \oplus \fG_{0}$ of (4.8) and (4.9).
\newline

$\mathit{The~ range~ of~ad~u~ :~ u~ \in U}$

\begin{proposition}  Let $\{\fN_{0} = U \oplus \fG_{0}, \ip \}$ be as in (6.7).  If u $\in$ U, then $im (ad~ u) = \fG_{0,u}^{\perp}$, where $\fG_{0,u} = \{X \in \fG_{0} : X(u) = 0 \}$ and $\fG_{0,u}^{\perp} = \{X \in \fG_{0} : \langle X, \fG_{o,u} \rangle = 0 \}$.
\end{proposition}

\begin{proof}  This is a trivial modification of the proof of (4.8).
\end{proof}

$\mathit{Surjectivity~ of~ ad~ u}$  The remarks following (4.8) also apply here, and the statements here can be deduced routinely from those following (4.8).  For completeness we assert that when dim U $\geq$ dim $\fG_{0}$, then ad u : U $\rightarrow \fG_{0}$ is surjective for all elements u in some Zariski open subset of U.
\newline

$\mathit{The~ range~ of~ ad~ u_{\lambda}~ :~  u_{\lambda} \in U_{\lambda}}$

    Fix a maximal abelian subalgebra $\fH_{0}$ of $\fG_{0}$ and let $\fH = \fH_{0}^{\bc}$ be the corresponding Cartan subalgebra of $\fG = \fG_{0}^{\bc}$.  Let $\lambda \in \Lambda$ and let $\varphi : V_{\lambda} \rightarrow U_{\lambda}$ be the $\br$ - linear isomorphism of (6.3) given by $\varphi(v_{\lambda}) = v_{\lambda} + J(v_{\lambda})$.  Let $\Sigma_{\lambda}$ be the proper Zariski closed subset of V$_{\lambda}$ defined following (4.8), and let $\Sigma'_{\lambda} = \varphi(\Sigma_{\lambda})$, a proper Zariski closed subset of U$_{\lambda}$.  Let $\tilde{H}_{\lambda}$ be the unique element of $\fH_{0}$ such that $\langle \tilde{H}, \tilde{H}_{\lambda} \rangle = - i \lambda(\tilde{H})$ for all $\tilde{H} \in \fH_{0}$.  As in the discussion preceding (4.9) we define $\Phi_{\lambda} = \{\alpha \in \Phi : V_{\lambda - \alpha} \neq \{0\} \}$.

\begin{proposition}  Let $\lambda \in \Lambda^{+}$ be an element such that there are no solutions to the equations $2\lambda = \alpha \in \Phi$ or $2\lambda = \alpha + \beta$, where $\alpha, \beta \in \Phi$.  Let u$_{\lambda}$ be any element of U$_{\lambda}$.  Then ad u$_{\lambda}(U) \subseteq \br\tilde{H}_{\lambda} \oplus \fG_{0,\lambda} \oplus$ \(\sum_{\alpha \in \Phi_{\lambda}} \fG_{0,\alpha}\).  Equality holds if u$_{\lambda} \in$ U$_{\lambda} - \Sigma'_{\lambda}$.
\end{proposition}

$\mathit{Remark}$  Very few weights $\lambda$ satisfy either of the equations $2\lambda = \alpha \in \Phi$ or $2\lambda = \alpha + \beta$, where $\alpha, \beta \in \Phi$.  The weights that do satisfy one of these equations are contained in the lists of (11.4) and (11.6).

\begin{proof}  We begin with the following observation.

\indent 1)  Let $\lambda \in \Lambda^{+}$ satisfy the hypotheses of the proposition.  If $\mu \in \Lambda$ is arbitrary, then the elements $\lambda + \mu$ and $\lambda - \mu$ cannot both lie in $\Phi$.

     If there exists $\mu \in \Lambda$ such that $\lambda + \mu = \alpha \in \Phi$ and $\lambda - \mu = \beta \in \Phi$, then $2\lambda = \alpha + \beta$, which was ruled out by hypothesis.  This proves 1).

    Next, let u $\in$ U be given.  By (6.4) there exist elements $u_{0}' \in U_{0}$ and $u_{\mu}' \in U_{\mu}, \mu \in \Lambda^{+}$, such that $u = u_{0}' +$ \(\sum_{\mu \in \Lambda^{+}} u_{\mu}' \).  Then ad $u_{\lambda}(u) = [u_{\lambda}, u_{0}'] +$ \(\sum_{\mu \in \Lambda^{+}} [u_{\lambda}, u_{\mu}'] \).  By  (6.7) and (6.10) we know that

\indent 2)$[u_{\lambda}, u_{0}'] \in \fG_{0,\lambda}$ and $[u_{\lambda},u_{\lambda}'] \in \br \tilde{H}_{\lambda}$.

    Now let $\mu \in \Lambda^{+}$ be an element distinct from $\lambda$.  By (6.7) we have $[u_{\lambda}, u_{\mu}'] \in \fG_{0,\lambda + \mu} \oplus \fG_{0,\lambda - \mu}$.  Now $\fG_{0,\lambda + \mu} = \{0\}$ unless $\lambda + \mu = \alpha \in \Phi$.  In this case $\alpha \in \Phi_{\lambda}$ and $\fG_{0,\lambda + \mu} = \fG_{0,\alpha}$.  Similarly, if $\fG_{0,\lambda - \mu} \neq \{0\}$, then $\lambda - \mu = \alpha \in \Phi_{\lambda}$.  By 1) $\lambda + \mu$ and $\lambda - \mu$ cannot both lie in $\Phi$.  We have proved

\indent 3)a) If $\lambda + \mu = \alpha \in \Phi_{\lambda}$, then $ [u_{\lambda}, u_{\mu}'] \in \fG_{0,\alpha}$

\hspace{.12in}b) If $\lambda - \mu = \alpha \in \Phi_{\lambda}$, then $[u_{\lambda}, u_{\mu}'] \in \fG_{0,\alpha}$

\hspace{.12in}c)   $[u_{\lambda}, u_{\mu}'] = 0$ in all other cases.

    The containment assertion of the proposition now follows from 1), 2) and 3).

    Next, let $u_{\lambda} \in U_{\lambda} - \Sigma_{\lambda}' $ be given.  Let J : $V \rightarrow V = U^{\bc}$ denote the conjugation map determined by U.
\newline

$\mathbf{Lemma~ 1}$  Let $u_{\lambda} \in U_{\lambda} - \Sigma'_{\lambda}$.  Then there exists an element $v_{\lambda} \in V_{\lambda} - \Sigma_{\lambda}$ such that

\indent     1) $u_{\lambda} = v_{\lambda} + J(v_{\lambda})$

\indent     2)  $X_{- \alpha}(v_{\lambda}) \neq 0$ and $X_{\alpha}(J(v_{\lambda})) \neq 0$ for all $\alpha \in \Phi_{\lambda}$.

\indent 3   $X_{- \lambda}(v_{\lambda}) \neq 0$ and $X_{\lambda}(J(v_{\lambda})) \neq 0$ if $\lambda \in \Lambda \cap \Phi$.

$\mathit{Proof~ of~ Lemma~ 1}$  Let $\varphi : V_{\lambda} \rightarrow U_{\lambda}$ be the $\br$ - linear isomorphism of 5) of (6.3) defined by $\varphi(v_{\lambda}) = v_{\lambda} + J(v_{\lambda})$.  If we choose $v_{\lambda} \in V_{\lambda}$ such that $u_{\lambda} = \varphi(v_{\lambda})$, then $v_{\lambda} \in V_{\lambda} - \Sigma_{\lambda}$ since $\varphi(\Sigma_{\lambda}) = \Sigma'_{\lambda}$ by the definition of $\Sigma'_{\lambda}$.  By the definition of $\Sigma_{\lambda}$, which precedes the statement of (4.9), we have  $X_{- \alpha}(v_{\lambda}) \neq 0$ for all $\alpha \in \Phi_{\lambda}$ and  $X_{- \lambda}(v_{\lambda}) \neq 0$ if $\lambda \in \Lambda \cap \Phi$.

    If $J_{0} : \fG \rightarrow \fG = \fG_{0}^{\bc}$ denotes the conjugation determined by $\fG_{0}$, then $J \circ X = J_{0}(X) \circ J$ for all $X \in \fG$.  For $\alpha \in \Phi_{\lambda}$ we obtain $0 \neq J(X_{- \alpha}(v_{\lambda})) = \{J_{0}(X_{- \alpha})\}(J(v_{\lambda}))$. which implies that $X_{\alpha}(J(v_{\lambda}) \neq 0$ since $J_{0} : \fG_{- \alpha} \rightarrow \fG_{\alpha}$ is an isomorphism of 1-dimensional complex vector spaces.  Similarly, if $\lambda \in \Lambda \cap \Phi$, then $0 \neq J(X_{- \lambda}(v_{\lambda})) = \{J_{0}(X_{- \lambda})\}(J(v_{\lambda}))$, which implies that $X_{\lambda}(J(v_{\lambda})) \neq 0$.  This proves Lemma 1.
\newline

$\mathbf{Lemma~ 2}$  Let $u_{\lambda} \in U_{\lambda} - \Sigma'_{\lambda}$.  If $\alpha \in \Phi_{\lambda}$, then ad $u_{\lambda}(U_{- \alpha + \lambda}) = \fG_{0,\alpha}$.

$\mathit{Proof~ of~ Lemma~ 2}$  We know from 2) of (6.7) that ad $u_{\lambda}(U_{- \alpha + \lambda}) \subseteq \fG_{0,\alpha} \oplus \fG_{0,2\lambda - \alpha} =  \fG_{0.\alpha}$ since $2\lambda - \alpha \notin \Phi$ by the hypotheses on $\lambda$.  To prove equality we show that \{ad~ u$_{\lambda}(U_{-\alpha + \lambda})\}^{\bc} = \fG_{0,\alpha}^{\bc}$.  By 3) of (6.3) and the remarks preceding (6.5) we have $\fG_{0,\alpha}^{\bc} = \fG_{\alpha} \oplus \fG_{- \alpha}$ and  \{ad u$_{\lambda}(U_{-\alpha + \lambda})\}^{\bc} =$ ad u$_{\lambda} \{U_{- \alpha + \lambda}\}^{\bc} =$ ad u$_{\lambda}(V_{- \alpha + \lambda} \oplus V_{\alpha - \lambda})$.  Hence it suffices to show that ad $u_{\lambda}(V_{- \alpha + \lambda}) = \fG_{- \alpha}$ and ad $u_{\lambda}(V_{\alpha - \lambda}) = \fG_{\alpha}$.

    As in Lemma 1 we write $u_{\lambda} = v_{\lambda} + J(v_{\lambda})$, where $v_{\lambda} \in V_{\lambda} - \Sigma_{\lambda}$.  Then by 2) of (4.7) ad v$_{\lambda}(V_{- \alpha + \lambda}) \subseteq \fG_{2\lambda - \alpha} = \{0\}$ since $2\lambda - \alpha \notin \Phi$.  Hence ad u$_{\lambda}(V_{- \alpha + \lambda}) =$ ad (J$(v_{\lambda})(V_{- \alpha + \lambda}) \subseteq \fG_{- \alpha}$ by 2) of (4.7) since J$(v_{\lambda}) \in V_{- \lambda}$ by (6.2).  To obtain equality we need only show that ad (J$(v_{\lambda})(V_{- \alpha + \lambda}) \neq \{0\}$ since $\fG_{- \alpha}$ is 1 - dimensional.  For $X_{\alpha} \in \fG_{\alpha}$ we have $B^{*}([J(v_{\lambda}), V_{- \alpha + \lambda}], X_{\alpha}) = B^{*}(X_{\alpha}(J(v_{\lambda})), V_{- \alpha + \lambda}) \neq 0$ by 1) of (4.5) since $X_{\alpha}(J(v_{\lambda}))$ is a nonzero element of $V_{\alpha - \lambda}$ by Lemma 1.  We conclude that $\{0\} \neq [J(v_{\lambda}), V_{- \alpha + \lambda}] = \fG_{- \alpha}$.

    Similarly, $[J(v_{\lambda}), V_{\alpha - \lambda}] \subseteq \fG_{0,\alpha - 2\lambda} = \{0\}$ since $2\lambda - \alpha \notin \Phi$, and hence ad$(u_{\lambda})(V_{\alpha - \lambda}) =$ ad$(v_{\lambda})(V_{\alpha - \lambda}) \subseteq \fG_{\alpha}$.  Since $X_{- \alpha}(v_{\lambda})$ is a nonzero element of $V_{\lambda - \alpha}$ by Lemma 1 it follows that $B^{*}([v_{\lambda}, V_{\alpha - \lambda}], X_{- \alpha}) = B^{*}(X_{- \alpha}(v_{\lambda}), V_{\alpha - \lambda}) \neq 0$, and we conclude that $\{0\} \neq$ ad $v_{\lambda}(V_{\alpha - \lambda}) = \fG_{\alpha}$.  This completes the proof of Lemma 2.

    We now complete the proof of (6.12).  We note that if $\lambda \in \Lambda \cap \Phi$, then ad $(u_{\lambda})(U_{0}) = \fG_{0,\lambda}$ by Lemma 2 and 3) of Lemma 1.

    In 2) above we showed that ad u$_{\lambda}(U_{\lambda}) \subseteq \br \tilde{H_{\lambda}}$.  It remains only to show equality, or equivalently that ad $u_{\lambda}(U_{\lambda}) \neq \{0\}$.  If $u_{\lambda} = v_{\lambda} + J(v_{\lambda})$ as above, then by (6.2) and the proof of 1) of (4.5) there exists $v'_{\lambda} \in V_{\lambda}$ such that $B^{*}(v_{\lambda}, J(v'_{\lambda})) = 1$.  Let $u'_{\lambda} = v'_{\lambda} + J(v'_{\lambda})$.  Since $2\lambda \notin \Phi$ by the hypotheses on $\lambda$ it follows from 2) of (4.7) that $[v_{\lambda}, v'_{\lambda}] = 0$ and  $[J(v_{\lambda}), J(v'_{\lambda})] = 0$.  Hence from 2) of (2.1) and the proof of 1) of (4.7) we obtain $[u_{\lambda},u'_{\lambda}] = [v_{\lambda}, J(v'_{\lambda})] + [J(v_{\lambda}), v'_{\lambda}] = B^{*}(v_{\lambda}, J(v'_{\lambda})) t_{\lambda} + B^{*}(J(v_{\lambda}), v'_{\lambda}) t_{\lambda} = \{B^{*}(v_{\lambda}, J(v'_{\lambda})) + \overline{B^{*}(v_{\lambda}, J(v'_{\lambda}))} \} t_{\lambda} = 2t_{\lambda}$.  Hence $[u_{\lambda}, U_{\lambda}] \neq 0$.
\end{proof}

$\mathit{Automorphisms~ and~ derivations~ of~ \fN_{0} = U \oplus \fG_{0}}$

    By the proofs of (4.3) and (4.10) and the definition of $\fN_{0}$ in section 2 we obtain

\begin{proposition}  Let $\fG_0$ be a compact, semisimple Lie algebra, and let U be a finite dimensional real $\fG_0$-module.  Let $\fN_{0} = U~ \oplus~ \fG_{0}$, and let G$_0$ be the connected subgroup of GL(U) whose Lie algebra is $\fG_0$.  For g $\in$ G$_0$ let T(g) denote the element of GL($\fN_{0}$) such that T(g) = g on U and T(g) = Ad(g) on $\fG_0$.  Then T(g) $\in$ Aut($\fN_{0}$) $\cap~$ I($\fN_{0}$) for all g $\in$ G$_0$, where  Aut($\fN_{0}$) and  I($\fN_{0}$) denote the automorphism and linear isometry groups of $\fN_{0}$.   The map T : G$_0 \rightarrow$ Aut($\fN_{0}$) $ \cap$ I($\fN_{0}$) is an injective homomorphism.
\end{proposition}

    For a complete description of Aut($\fN_{0}$) $\cap~$ I($\fN_{0}$) and its Lie algebra see Theorem 3.12 of [L].  As an immediate consequence of the result above we obtain

\begin{corollary}  Let $\fG_0$ be a compact, semisimple Lie algebra, and let U be a finite dimensional real $\fG_0$-module.  Let $\fN_{0} = U~ \oplus~ \fG_{0}$.  For X $\in \fG_{0}$ let t(X) denote the element of End($\fN_{0}$) such that t(X) = X on U and t(X) = ad X on $\fG_{0}$.  Then t(X) is a skew symmetric derivation of $\fN_{0}$ for all X $\in \fG_{0}$.  The map t : $\fG_{0} \rightarrow$ Der($\fN_{0}$) is an injective homomorphism.
 \end{corollary}

 $\mathit{The~ Weyl~ group~ and~ Aut(\fN_{0})}$

    We conclude this section with the real analogue of (5.3).

\begin{proposition}  Let $V = U^{\bc}$.  Let $\{T_{\sigma} : \sigma \in W\}$ be the elements of GL(V) constructed in (5.1), and let  $\{\varphi_{\sigma} : \sigma \in W\}$ be the corresponding elements of Aut($\fG$) $\subset$ GL($\fG$).  Let $\fN = V \oplus \fG$ and let $\zeta_{\sigma} \in$ Aut($\fN$) be the element of (5.3) such that $\zeta_{\sigma} = T_{\sigma}$ on V and $\zeta_{\sigma} = \varphi_{\sigma}$ on $\fG$.  Then for $\fN_{0} = U \oplus \fG_{0}$ we have $\zeta_{\sigma}(\fN_{0}) = \fN_{0}$ and $\zeta_{\sigma} \in$ Aut($\fN_{0}$).  Moreover,  $\zeta_{\sigma}$ is an isometry of $\fN_{0}$ with respect to any $\fG_{0}$ - invariant inner product on $\fN_{0}$.
\end{proposition}

\begin{proof}  By (5.2) T$_{\sigma}$(U) = U, and hence in the special case that U $= \fG_{0}$ and V $= U^{\bc} = \fG$ it follows that $\varphi_{\sigma} \in$ Aut($\fG_{0})$.  Hence $\zeta_{\sigma}(\fN_{0}) = \fN_{0} = U \oplus \fG_{0}$.  By (5.3) $\zeta_{\sigma} \in$ Aut($\fN$) and it follows that the restriction of $\zeta_{\sigma}$ to the real subalgebra $\fN_{0}$ lies in Aut($\fN_{0}$).  Finally, let $\ip$ be a $\fG_{0}$ - invariant inner product on $\fN_{0}$.  Then $\ip$ extends to a $\fG_{0}$ - invariant Hermitian inner product on $\fN = V \oplus \fG = \fN_{0}^{\bc}$.  By (5.3) $\zeta_{\sigma}$ is unitary with respect to $\ip$ on $\fN$, and hence $\zeta_{\sigma}$ is an orthogonal linear transformation of $\fN_{0}$ relative to $\ip$.
\end{proof}

\section{Adapted~ bases~ of~ $\fG_0$}

$\mathit{The~ elements~ A_{\beta}, B_{\beta} : \beta \in \Phi}$

    Let  \rs ~be the root space decomposition of $\fG~ = ~\fG_0^{\bc}$ determined by $\fH~ = ~\fH_0^{\bc}$.  Let J$_0$ : $\fG \rightarrow \fG$ denote the conjugation in $\fG$ determined by $\fG_0$.  For $\beta \in \Phi^{+}$ let X$_{\beta}$ be a nonzero vector in $\fG_{\beta}$, and set A$_{\beta}$ = X$_{\beta}$ + J$_0$(X$_{\beta}$) and B$_{\beta}$ = i (X$_{\beta}~ -~ $J$_0$(X$_{\beta}$)).  The elements  \{A$_{\beta}$, B$_{\beta}$\} are fixed by J$_{0}$ and hence form an $\br$ - basis of $\fG_{0,\beta}$.  By (6.1) $\tilde{\tau}_{\alpha} =$ i $\tau_{\alpha} \in \fH_0 \subset \fG_0$.

    If $\fC_{0} = \{ \tau_{\alpha}~:~ \alpha \in \Delta~; A_{\beta},B_{\beta}~ :~ \beta \in \Phi^{+}\}$, then
$\fG_0'~ =~ \br-$span$ (\fC_0) \subseteq \fG_0$ and equality holds since $(\fG_0')^{\bc}~ =~ (\fG_0)^{\bc}~ =~ \fG$.
It is easy to see that the elements of $\fC_{0}$ are linearly independent over $\br$, and hence $\fC_{0}$  is a basis of $\fG_{0}$.
We call $\fC_{0}$ an $\mathit{adapted~ basis}$ of $\fG_{0}$.
\newline

$\mathit{Bracket~ relations~ between~ A_{\beta}}~ and~ B_{\beta}$

    In the next result we list the bracket relations between A$_{\beta}$, B$_{\beta}$ and iH for each $\beta \in \Phi^{+}$ and each H $\in \fH$.  We recall from (6.1) and (6.2) that i $\fH_0 =
\fH_{\br}~ =~ \{$H$~ \in \fH~ =~ \fH_0^{\bc}~ :~ \beta($H$) \in \br$ for all $\beta \in \Phi \}$ and J$_0(\fG_{\beta})~ =~ \fG_{-\beta}$.  The proof of the result below now follows routinely from the definitions, the facts above and 4) of (4.6).

\begin{proposition} Let A$_{\beta}$ = X$_{\beta}$ + J$_0$(X$_{\beta}$) and B$_{\beta}$= i (X$_{\beta}$ $-$ J$_0$(X$_{\beta}$))
for $\beta \in \Phi^{+}$,  and let H $\in \fH~ =~ \fH_0^{\bc}$.  Then

\indent 1)  [i H, A$_{\beta}$] = $\beta$(H) B$_{\beta}$.

\indent 2)  [i H, B$_{\beta}$] = $-~ \beta$(H) A$_{\beta}$.

\indent 3)  [A$_{\beta}$ , B$_{\beta}$] = $-$2 B(X$_{\beta}$, J$_0$(X$_{\beta}$)) (i H$_{\beta}$)
\end{proposition}

$\mathit{Remark}$  If X$_{\beta} \in \fG_{\beta}$ is chosen carefully for each $\beta \in \Phi^{+}$, then the structure constants of the $\br~ -~$basis $\fC_0$ defined above are integers and J$_0$(X$_{\beta}$) = $-$ X$_{- \beta}$ for all $\beta \in \Phi$.  See 3) of (3.2) for further details.  In this case A$_{\beta}$ = X$_{\beta}$ + J$_0$(X$_{\beta}$) =  X$_{\beta}~-~$X$_{-\beta}$ and B$_{\beta}$ = i (X$_{\beta}~ -~ $J$_0$(X$_{\beta}$)) =  i X$_{\beta}$ + i X$_{-\beta}$ for all $\beta \in \Phi$.
\newline

$\mathit{Kernel~ of~ A_{\beta}, B_{\beta}}$

\begin{proposition}  Let $\beta \in \Phi^{+}~ , ~\lambda \in \Lambda$.  Let A$_{\beta}$
= X$_{\beta}$ + J$_0$(X$_{\beta})$ and B$_{\beta}$ =  i (X$_{\beta}~ -~$J$_0$(X$_{\beta}$)).  Let Ker (A$_{\beta}$) and
Ker (B$_{\beta}$) denote the kernels of the linear maps A$_{\beta}$ : U$_{\lambda}\rightarrow$ U$_{\lambda+\beta}~ \oplus$
U$_{\lambda-\beta}$ and B$_{\beta}$ : U$_{\lambda}\rightarrow$ U$_{\lambda+\beta}~ \oplus$ U$_{\lambda-\beta}$ respectively.  Then

\indent 1)  Ker (A$_{\beta}$) = Ker (B$_{\beta}$).  The maps A$_{\beta}$ and B$_{\beta}$ are identically zero on U$_{\lambda}$
$\iff$ $V_{\beta + \lambda} = \{0\}$ and $V_{\beta - \lambda} = \{0\}$.

\indent 2)  If $\lambda$(H$_{\beta}$) $\neq$ 0 , then Ker (A$_{\beta}$) = Ker (B$_{\beta}$) = \{0\}.

\indent 3)  If u$_{\lambda} \notin$ Ker (A$_{\beta}$) = Ker (B$_{\beta}$) , then A$_{\beta}$(u$_{\lambda}$) and
B$_{\beta}$(u$_{\lambda}$) are linearly independent.
\end{proposition}

$\mathit{Remark}$  The analogue of this result for $\lambda~ =~ 0$ also holds, and the proofs in this case are trivial modifications of the proofs of the result above.

\begin{proof}  1)  Let u$_{\lambda}$ be any element of U$_{\lambda}$, and write u$_{\lambda}$ = v$_{\lambda}$ + v$_{-\lambda}$,
where v$_{\lambda} \in$ V$_{\lambda}$ and v$_{-\lambda} \in$ V$_{-\lambda}$.  We compute

\indent A$_{\beta}$(u$_{\lambda}$) =  X$_{\beta}$(v$_{\lambda}$) + X$_{\beta}$(v$_{-\lambda}$) + J$_0$(X$_{\beta})$(v$_{\lambda}$)
+ J$_0$(X$_{\beta})$(v$_{-\lambda}$).

\indent B$_{\beta}$(u$_{\lambda}$) =  i X$_{\beta}$(v$_{\lambda}$) + i X$_{\beta}$(v$_{-\lambda}$)
 $-$ i J$_0$(X$_{\beta})$(v$_{\lambda}$) $-$ i J$_0$(X$_{\beta}$)(v$_{-\lambda}$).

The four components of A$_{\beta}$(u$_{\lambda}$) lie in different weight spaces of V, and the same is true for
 B$_{\beta}$(u$_{\lambda}$).  Hence A$_{\beta}$(u$_{\lambda}$) = 0 $\iff$ the vectors
 \{X$_{\beta}$(v$_{\lambda}$), X$_{\beta}$(v$_{-\lambda}$), J$_0$(X$_{\beta})$(v$_{\lambda}$), J$_0$(X$_{\beta})$(v$_{-\lambda}$)\}
are all zero $\iff$ B$_{\beta}$(u$_{\lambda}$) = 0.  This proves that Ker (A$_{\beta}$) = Ker (B$_{\beta}$).

    If $V_{\beta + \lambda} = \{0\}$ and $V_{\beta - \lambda} = \{0\}$. then $V_{- \beta - \lambda} = \{0\}$ and $V_{- \beta + \lambda} = \{0\}$ by (6.2).  It follows that A$_{\beta}(U_{\lambda})$ and B$_{\beta}(U_{\lambda})$
are subspaces of $U_{\beta + \lambda}~ \oplus~ U_{\beta - \lambda} \subset V_{\beta + \lambda}~ \oplus~ V_{-\beta - \lambda} \oplus~
V_{\beta - \lambda} \oplus~ V_{-\beta + \lambda}  =$ \{0\}.  Hence
 A$_{\beta}$ and B$_{\beta}$ are identically zero on $U_{\lambda}$.  Conversely, suppose that A$_{\beta}$ is identically zero
 on $U_{\lambda}$ (and hence B$_{\beta}$ is identically zero on $U_{\lambda}$).  Let X$_{\beta}$ be a nonzero element of $\fG_{\beta}$. Let $v_{\lambda} \in V_{\lambda}$ be given, and let u$_{\lambda} = v_{\lambda} + J(v_{\lambda}) \in V_{\lambda}~ \oplus~V_{-\lambda}$, where J : V $\rightarrow$ V denotes the conjugation determined by U.  Then $u_{\lambda} \in U_{\lambda}$ since J$(u_{\lambda}) =
 u_{\lambda}$.  From the discussion above it follows that $0 = X_{\beta}(v_{\lambda}) = J_{0}(X_{\beta})(v_{\lambda}))$ since $0 = A_{\beta}(u_{\lambda})$.  By (6.2) $J_{0}(X_{\beta}) \in
 \fG_{- \beta}$ and since $\fG_{\beta}$ and $\fG_{- \beta}$ are 1-dimensional it follows from (4.1) that
 $V_{\lambda + \beta} = \{0 \}$ and $V_{\lambda - \beta} = \{0 \}$.  Hence $V_{\beta - \lambda} = J(V_{\lambda - \beta}) = \{0 \}$, which proves 1).

\indent 2)  Let u$_{\lambda}$ be an element of Ker (A$_{\beta}$) = Ker (B$_{\beta}$), and write u$_{\lambda}$ = v$_{\lambda}$ + v$_{-\lambda}$ as above.    By the proof of 1) X$_{\beta}$(v$_{\lambda}$) = J$_0$(X$_{\beta})$(v$_{\lambda}$) = 0.  By 1) and 4) of (4.6)
[X$_{\beta}$, J$_0$(X$_{\beta}$)] = B(X$_{\beta}$, J$_0$(X$_{\beta}$)) H$_{\beta} \neq~0$.  Hence
$\lambda$(H$_{\beta}$) v$_{\lambda}$ = H$_{\beta}$(v$_{\lambda}$) = (1 / B(X$_{\beta}$, J$_0$(X$_{\beta}$))) [X$_{\beta}$, J$_0$(X$_{\beta}$](v$_{\lambda}$) = (1 / B(X$_{\beta}$, J$_0$(X$_{\beta}$))) \{X$_{\beta}$ J$_0$(X$_{\beta}$)(v$_{\lambda}$) $-$ J$_0$(X$_{\beta}$)X$_{\beta}$(v$_{\lambda}$)\} = 0.  It follows that v$_{\lambda}$ = 0 since $\lambda$(H$_{\beta}$) $\neq$ 0.  A similar argument shows that v$_{-\lambda}$ = 0.  Hence u$_{\lambda}$ = v$_{\lambda}$ + v$_{-\lambda}$ = 0.

\indent 3)  The proof of 1) expresses A$_{\beta}(u_{\lambda})$ and B$_{\beta}(u_{\lambda})$ as a sum of components from four distinct weight spaces of V, at least two of which are nonzero since A$_{\beta}(u_{\lambda})$ and B$_{\beta}(u_{\lambda})$ are both nonzero.  Moreover, each component of B$_{\beta}(u_{\lambda})$ is $\pm$ i times the corresponding component of A$_{\beta}(u_{\lambda})$.  Hence it is impossible to satisfy the relation B$_{\beta}(u_{\lambda}) = c A_{\beta}(u_{\lambda})$ for some nonzero real number c.  This proves 3).
\end{proof}


$\mathit{Range~ of~ A_{\beta}, B_{\beta}}$

    We now establish some further properties of the maps A$_{\beta}$ : U$_{\lambda} \rightarrow~ $U$_{\lambda+\beta}~ \oplus$
U$_{-\lambda-\beta}$ and B$_{\beta}$ : U$_{\lambda} \rightarrow~ $U$_{\lambda+\beta}~ \oplus$
U$_{-\lambda-\beta}$.  The next result sharpens (6.6).

\begin{proposition} Let $\beta \in \Phi^{+}$ and $\lambda \in \Lambda$. Then $\fG_{0,\beta}$(U$_{\lambda}$) = Re \{$\fG_{\beta}$(V$_{\lambda}~ \oplus$ V$_{-\lambda}$)\} =
Re \{$\fG_{-\beta}$(V$_{\lambda}~ \oplus$ V$_{-\lambda}$)\} = Im \{$\fG_{\beta}$(V$_{\lambda}~ \oplus$ V$_{-\lambda}$)\} = Im \{$\fG_{-\beta}$(V$_{\lambda}~ \oplus$ V$_{-\lambda}$)\}
\end{proposition}

\begin{proof}  We recall from (6.2) that J $\circ$ X = J$_{0}$(X) $\circ$ J on V for all X $\in \fG$.  If X = X$_{\beta}$, a nonzero element of $\fG_{\beta}$, then it follows from this observation and (6.2) that
 J\{$\fG_{\beta}(V_{\lambda}~ \oplus V_{-\lambda}$)\}$= \fG_{-\beta}(V_{\lambda}~ \oplus V_{-\lambda})$.  Hence, the second equality in the statement follows from the first, and the last two equalities follow immediately from the first two.

    It suffices to prove the first equality in the statement.  Since A$_{\beta}$ = X$_{\beta}$ + J$_0$(X$_{\beta}$) and B$_{\beta}$ = i (X$_{\beta} -$ J$_0$(X$_{\beta}$)) it follows that A$_{\beta} -$ i  B$_{\beta}$ = 2X$_{\beta}$.  Let  v$_{\lambda} \in$ V$_{\lambda}$ and recall that J(v$_{\lambda}$) $\in$ J(V$_{\lambda}$) = V$_{-\lambda}$ by (6.2).  From the definitions and 2) of (6.3) we obtain

\indent 1)  2 Re (X$_{\beta}$(v$_{\lambda}$)) = A$_{\beta}$(Re v$_{\lambda}$) + B$_{\beta}$(Im v$_{\lambda}$) $\in$
$\fG_{0,\beta}$(U$_{\lambda}$)

\indent 2)  2 Re (X$_{\beta}$(J(v$_{\lambda}$))) = A$_{\beta}$(Re v$_{\lambda}$) $-$ B$_{\beta}$(Im v$_{\lambda}$)
$\in$ $\fG_{0,\beta}$(U$_{\lambda}$)

It follows that Re \{$\fG_{\beta}$(V$_{\lambda}~ \oplus$ V$_{-\lambda}$)\} $\subseteq \fG_{0,\beta}$(U$_{\lambda}$).  From 1)
and 2) we obtain

\indent 3)  A$_{\beta}$(Re v$_{\lambda}$) = Re X$_{\beta}$(v$_{\lambda}$ + J(v$_{\lambda}$)) $\in$
Re \{$\fG_{\beta}$(V$_{\lambda}~ \oplus$ V$_{-\lambda}$)\}

\indent 4)  B$_{\beta}$(Re v$_{\lambda}$) = Re X$_{\beta}$(i v$_{\lambda}-$ J(i v$_{\lambda}$))
 $\in$ Re \{$\fG_{\beta}$(V$_{\lambda}~ \oplus$ V$_{-\lambda}$)\}

Hence $\fG_{0,\beta}$(U$_{\lambda}$) $\subseteq$ Re \{$\fG_{\beta}$(V$_{\lambda}~ \oplus$ V$_{-\lambda}$)\ since A$_{\beta}$ and B$_{\beta}$ span $\fG_{0,\beta}$ and U$_{\lambda}$ = Re (V$_{\lambda}$) by 2) of (6.3).
\end{proof}

$\mathit{Nonsingular~ subspaces~ for~ A_{\beta}, B_{\beta}}$

    If $V_{\beta + \lambda} = \{0\}$ and $V_{\beta - \lambda} = \{0\}$,  then by 1) of (7.2) both A$_{\beta}$ and B$_{\beta}$ are identically zero on U$_{\lambda}$.  If either $V_{\beta + \lambda} \neq \{0\}$ or $V_{\beta - \lambda} \neq \{0\}$,  then we show next that A$_{\beta}$ and B$_{\beta}$ are both nonsingular on a naturally defined nonzero subspace U$_{\lambda,\beta}$ of U$_{\lambda}$.

\begin{proposition}  Let $\beta \in \Phi^{+}, \lambda \in \Lambda$ and suppose that $V_{\beta+\lambda} \neq \{0\}$ or  $V_{\beta - \lambda} \neq \{0\}$.  Let k $\geq$ 0 and j $\geq$ 0 be the largest integers such that $V_{\lambda + k\beta} \neq \{0\}$ and $V_{\lambda - j\beta} \neq \{0\}$.  Let X$_{\beta}$ and X$_{-\beta}$ be nonzero elements of $\fG_{\beta}$ and $\fG_{-\beta}$ respectively.  Then

\indent 1)  (X$_{-\beta}$)$^k$(V$_{\lambda + k\beta}$) = (X$_{\beta}$)$^j$(V$_{\lambda - j\beta}$) : = V$_{\lambda,\beta} \subset$ V$_{\lambda}$.
The subspace V$_{\lambda,\beta}$ is nonzero and independent of the choice of X$_{\beta}$ and X$_{-\beta}$.

\indent 2) If k $\geq$ 1, then X$_{\beta}$ : V$_{\lambda,\beta} \rightarrow$ V$_{\lambda+\beta}$ is injective.  If j $\geq 1$, then
X$_{-\beta}$ : V$_{\lambda,\beta} \rightarrow$ V$_{\lambda-\beta}$ is injective.

\indent 3)  J(V$_{\lambda,\beta}$) = V$_{-\lambda,\beta}$, where J : V $\rightarrow$ V is the conjugation defined by U.

\indent 4)  Let U$_{\lambda,\beta}$ = (V$_{\lambda,\beta}~ \oplus~$V$_{-\lambda,\beta}$) $\cap$ U.   Then

\indent \indent     a) $U_{\lambda,\beta} = Re (V_{\lambda,\beta}) = Im (V_{\lambda,\beta})  = Re (V_{- \lambda,\beta}) = Im (V_{- \lambda,\beta}) = U_{- \lambda, \beta} \subset U_{\lambda}$.

\indent \indent     b)  $(U_{\lambda,\beta})^{\bc} = V_{\lambda,\beta} \oplus  V_{- \lambda,\beta}$.

\indent \indent     c)  dim$_{\br}$U$_{\lambda,\beta}$ = dim$_{\bc}$V$_{\lambda+k\beta}$ = dim$_{\bc}$V$_{\lambda-j\beta}$   $>$ 0.

 \indent    5)  A$_{\beta}$ : U$_{\lambda,\beta}~\rightarrow~$U$_{\lambda+\beta}~ \oplus~$U$_{\lambda-\beta}$
 and B$_{\beta}$ : U$_{\lambda,\beta}~\rightarrow~$U$_{\lambda+\beta}~ \oplus~$U$_{\lambda-\beta}$ are injective.
 \end{proposition}

\begin{proof}  By the lemma in the proof of (4.1) the elements \{$\lambda + m \beta : m \in \bz$\} that lie in $\Lambda$ are precisely the elements $\{\lambda + m \beta : - j \leq m \leq k\}$.  Since the Weyl reflection $\sigma_{\beta}$ leaves invariant this weight string it follows that  $\sigma_{\beta}(\lambda+k\beta) = \lambda-j\beta$ and $\sigma_{\beta}(\lambda-j\beta) = \lambda+k\beta$.  It is known that dim V$_{\mu} =$ dim V$_{\sigma(\mu)}$ for all $\mu \in \Lambda$ and all $\sigma \in$ W.  (For example, this follows from 2) of (5.1)).   In particular V$_{\lambda+k\beta}$ and V$_{\lambda-j\beta}$ have the same dimension.

    Recall that $\lambda(\tau_{\beta}) + k = j$ by the lemma in the proof of (3.3).  Let N $= \lambda(\tau_{\beta}) + 2k = k+j \geq k$. In the proof of (4.1) we showed that
    (X$_{-\beta}$)$^N$(V$_{\lambda+k\beta}$) $\neq \{0\}$ but (X$_{-\beta}$)$^{N+1}$(V$_{\lambda+k\beta}$) = $\{0\}$.   That proof also shows that
(X$_{-\beta}$)$^N$ : V$_{\lambda+k\beta} \rightarrow$ V$_{\lambda - j\beta}$ is injective and hence an isomorphism since V$_{\lambda+k\beta}$ and V$_{\lambda - j\beta}$ have the same dimension.

    In the proof of (4.1) the relations in a) show that
(X$_{\beta}$)$^j$(X$_{-\beta}$)$^N$(V$_{\lambda+k\beta}$) = (X$_{-\beta}$)$^{N-j}$(V$_{\lambda+k\beta}$) = (X$_{-\beta}$)$^k$(V$_{\lambda+k\beta}$). We conclude that
(X$_{-\beta}$)$^k$(V$_{\lambda+k\beta}$) = (X$_{\beta}$)$^j$(V$_{\lambda-j\beta}$) : = V$_{\lambda,\beta}$ since (X$_{-\beta}$)$^N$ : V$_{\lambda+k\beta} \rightarrow$ V$_{\lambda-j\beta}$ is an isomorphism.  The subspace V$_{\lambda,\beta}$ is independent of the choice of X$_{\beta}$ and X$_{-\beta}$ since $\fG_{\beta}$ and $\fG_{-\beta}$ are 1-dimensional.  The subspace V$_{\lambda,\beta} =$ (X$_{-\beta}$)$^k$(V$_{\lambda+k\beta}$) is nonzero since k $\leq$ N and (X$_{-\beta}$)$^{N}$(V$_{\lambda+k\beta}$) $=$ V$_{\lambda-j\beta}$ is nonzero. This proves 1).

\indent 2)  Note that k $\geq 1 \iff V_{\beta + \lambda} \neq \{0 \}$and j $\geq 1 \iff V_{\beta - \lambda} \neq \{0 \}$.  If k $\geq$ 1, then the proof of (4.1) shows that X$_{\beta} : V_{\lambda,\beta} \rightarrow V_{\lambda+\beta}$ is injective, and we observed above that
$(X_{-\beta})^{N} : V_{\lambda+k\beta} \rightarrow V_{\lambda-j\beta}$ is an isomorphism.  In particular, if j $\geq$ 1, then
X$_{-\beta}$ is injective on V$_{\lambda,\beta} = (X_{-\beta})^{k}(V_{\lambda+k\beta})$ since N = k+j $\geq$ k+1.

\indent 3)  Recall that J$_0$ : $\fG \rightarrow \fG$ is  the conjugation defined by $\fG_0$.  By (6.2) and induction on k we obtain

\indent \indent     i)  J $\circ$ X$_{\alpha}^k$ = J$_0$(X$_{\alpha}$)$^k \circ$ J on V for all integers k $\geq$ 1 and all $\alpha \in \Phi$.

    By the lemma in (4.1) $V_{\lambda + m\beta} \neq \{0 \}  \iff - j \leq m \leq k$.  It follows from (6.2) that

\indent \indent     ii)  $V_{-\lambda+m\beta} \neq \{0 \} \iff - k \leq m \leq j$

From 1), 3i) and (6.2) we obtain J(V$_{\lambda,\beta}$) = (J $\circ$ (X$_{-\beta}$)$^k)$(V$_{\lambda+k\beta}$) = (J$_0$(X$_{-\beta}$)$^k~ \circ$ J)
(V$_{\lambda+k\beta}$) = (J$_0$(X$_{-\beta}$)$^k$)(V$_{-\lambda-k\beta}$).  Note that J$_0$(X$_{-\beta}$) $\in \fG_{\beta}$ by (6.2).  Similarly, J$_0$(X$_{\beta}$) $\in \fG_{-\beta}$ and J(V$_{\lambda,\beta}$) = (J $\circ$ (X$_{\beta}$)$^j$)(V$_{\lambda-j\beta}$) = (J$_0$(X$_{\beta}$)$^j \circ$ J)(V$_{\lambda-j\beta}$) = J$_0$(X$_{\beta}$)$^j$(V$_{-\lambda+j\beta}$).  Hence
J$_0$(X$_{-\beta}$)$^k$(V$_{-\lambda-k\beta}$) = J$_0$(X$_{\beta}$)$^j$(V$_{-\lambda+j\beta}$) = J(V$_{\lambda,\beta}$).  From ii) and the definition of V$_{-\lambda,\beta}$ in 1) we see that
the proof of 3) is complete.

    4)  Clearly U$_{\lambda,\beta} \subseteq$ (V$_{\lambda}~ \oplus$ V$_{-\lambda}$) $\cap$ U = U$_{\lambda}$ since V$_{\lambda,\beta}
\subset$ V$_{\lambda}$ and V$_{-\lambda,\beta} \subset$ V$_{-\lambda}$.  The proof of 4a) and 4b) is very similar to the proof of 2) of (6.3), using the fact from 3) that J(V$_{\lambda,\beta}$) = V$_{-\lambda,\beta}$.  We omit the details.

    We prove 4c).  It follows from 4b) that dim$_{\br}$U$_{\lambda,\beta}$ = dim$_{\bc}$V$_{\lambda,\beta}$ = dim$_{\bc}$V$_{\lambda+k\beta}$ since
(X$_{-\beta}$)$^k$ : V$_{\lambda+k\beta} \rightarrow$ V$_{\lambda,\beta}$ is injective, hence an isomorphism, by the proof of (4.1).  Similarly the proof of (4.1) shows that (X$_{\beta}$)$^j$ : V$_{\lambda-j\beta} \rightarrow$ V$_{\lambda,\beta}$ is also an isomorphism.  This proves 4c).

\indent 5)  By 1) of (7.2) Ker (A$_{\beta}$) = Ker (B$_{\beta}$) in U$_{\lambda}$.  Let u$_{\lambda} \in$ U$_{\lambda,\beta}$, and choose v$_{\lambda}\in$ V$_{\lambda,\beta}$ such that u$_{\lambda}$ = Re(v$_{\lambda}$).  If 0 = A$_{\beta}$(u$_{\lambda}$) = B$_{\beta}$(u$_{\lambda}$), then by the proof of 1) of (7.2) we see that $X_{\beta}(v_{\lambda}) = 0$ and $(J_{0}(X_{\beta}))(v_{\lambda}) = 0$.  Hence v$_{\lambda}$ = 0 by 2) since $X_{- \beta} = J_{0}(X_{\beta}) \in \fG_{- \beta}$.
\end{proof}

    We obtain an analogue of (4.1).

\begin{proposition}  Let $\beta \in \Phi^{+}$ and $\lambda \in \Lambda$.  Let X$_{\beta}$ be a nonzero element of $\fG_{\beta}$. Let A$_{\beta}$ = X$_{\beta} +$ J$_0$(X$_{\beta}$) $\in \fG_0$ and B$_{\beta}$ = i (X$_{\beta} -$ J$_0$(X$_{\beta}$)) $\in \fG_0$.  Then the following are equivalent :

\indent 1)  $V_{\beta + \lambda} \neq \{0\}$ or $V_{\beta - \lambda} \neq \{0\}$.

\indent 2)  $\fG_{0,\beta}$(U$_{\lambda}$) is a nonzero subspace of U$_{\beta + \lambda}~ \oplus~ $U$_{\beta - \lambda}$.

\indent 3)  $Ker(A_{\beta})~ \cap~ U_{\lambda} =  Ker(B_{\beta})~ \cap~ U_{\lambda} \neq U_{\lambda}$.  If u$_{\lambda} \notin Ker(A_{\beta})~ \cap~ U_{\lambda} =  Ker(B_{\beta})~ \cap~ U_{\lambda}$ , then A$_{\beta}$(u$_{\lambda}$) and B$_{\beta}$(u$_{\lambda}$) are linearly independent in U$_{\lambda + \beta} \oplus$ U$_{\lambda - \beta}$.
\end{proposition}

\begin{proof} The assertion 1) \imp 3) follows from 1) and 3) of (7.2).  The assertion 3) \imp 2)  is obvious since A$_{\beta}$ and B$_{\beta}$ span $\fG_{0,\beta}$.  If 2) holds, then either V$_{\beta+\lambda} \neq \{0\}$ or V$_{\beta - \lambda} \neq \{0\}$ by (7.3). It follows that 2) \imp 1).
 \end{proof}

$\mathit{\fG_0-modules~ with~ nontrivial~ zero~ weight~ space}$

    The zero weight space U$_0$ of the real $\fG_0$-module U is nonzero $\iff$ the zero weight space V$_0$ of the complex $\fG$-module V = U$^{\bc}$ is nonzero by 2) of (6.3).  By (4.2) V$_{0}$ is nonzero if $\Lambda \cap \Phi$ is nonempty.  As we remarked following (4.2), if $\Lambda \cap \Phi$ is nonempty, then $\Phi \subset \Lambda$, with only a few exceptions.  See [D] for details.

    We now obtain some sharper versions of (4.10) and (7.5) that are useful for applications.  See [D] for an application to the density of closed geodesics on a compact nilmanifold $\Gamma_{0} \backslash N_{0}$, where N$_{0}$ is a simply connected 2-step nilpotent Lie group with a
 left invariant metric and Lie algebra $\fN_{0} = U \oplus \fG_{0}$ and $\Gamma_{0}$ is a discrete, cocompact subgroup of N$_{0}$.

    If V$_0$ and U$_0$ are nonzero , then we obtain the following analogue for $\lambda$ = 0 of (7.5).  We omit the proof, which is virtually identical if one uses the analogue for $\lambda = 0$ of (7.2).

\begin{proposition}  Suppose that the zero weight space U$_0$ is nonzero.  Let $\beta \in \Phi^{+}$ and let X$_{\beta}$ be a nonzero element of
$\fG_{\beta}$.  Let A$_{\beta}$ = X$_{\beta} +$ J$_0$(X$_{\beta}$) $\in \fG_{0,\beta}$ and B$_{\beta}$ = i (X$_{\beta} -$ J$_0$(X$_{0,\beta}$))
 $\in \fG_{0,\beta}$.
Then the following are equivalent :

\indent 1)  $V_{\beta} \neq \{0\}$

\indent 2)  $\fG_{0,\beta}$(U$_0$) is a nonzero subspace of U$_{\beta}$.

\indent 3)  $Ker(A_{\beta})~ \cap~ U_{0} =  Ker(B_{\beta})~ \cap~ U_{0}  \neq U_0$.  If u$_{0} \notin Ker(A_{\beta})~ \cap~ U_{0} =  Ker(B_{\beta}) ~ \cap~ U_{0}$ , then A$_{\beta}$(u$_0$) and
B$_{\beta}$(u$_0$) are linearly independent in U$_{\beta}$.
\end{proposition}

    In the next three results we consider roots that are also weights, and we obtain more detailed information about the bracket operation in
$\fN = U \oplus \fG_{0}$.

\begin{proposition}  Let $\beta \in \Phi^{+} \cap \Lambda$.  Let u$_0 \in$ U$_0$ with  u$_0 \notin$
Ker(A$_{\beta}$) $\cap~ U_{0} =$  Ker(B$_{\beta}$) $\cap~ U_{0}$ ,  and  consider  ad u$_0$ : U$_{\beta} \rightarrow \fG_{0,\beta}$.  Then

\indent 1)  Ker (ad u$_0$) = \{$\br$-span (A$_{\beta}$(u$_0$), B$_{\beta}$(u$_0$))\}$^{\perp}$, where $^{\perp}$ denotes the orthogonal complement.

\indent 2)  The elements [u$_0$, A$_{\beta}$(u$_0$)] and [u$_0$, B$_{\beta}$(u$_0$)] span $\fG_{0,\beta}$.  In particular,  ad u$_0$ : U$_{\beta} \rightarrow \fG_{0,\beta}$ is surjective.
\end{proposition}

\begin{proof}   1)  Recall from the discussion at the beginning of section 7 that A$_{\beta}$ and B$_{\beta}$ form a basis of $\fG_{0,\beta}$.  If u$_{\beta} \in$ U$_{\beta}$, then u$_{\beta} \in$ Ker ad u$_0$ $\iff$ 0 = $\langle$ [u$_0$, u$_{\beta}$], A$_{\beta} \rangle =
\langle$ A$_{\beta}$(u$_0$), u$_{\beta} \rangle$ and 0 = $\langle$ [u$_0$, u$_{\beta}$], B$_{\beta} \rangle = \langle$ B$_{\beta}$(u$_0$), u$_{\beta} \rangle$.  The proof is complete.  Note that $\br$-span \{A$_{\beta}$(u$_0$), B$_{\beta}$(u$_0$)\} is 2-dimensional by 3) of (7.6).

\indent 2)  By (6.6) and 1) of (6.7) the elements [u$_0$, A$_{\beta}$(u$_0$)] and
[u$_0$, B$_{\beta}$(u$_0$)] lie in $\fG_{0,\beta}$.  It suffices to show that they are linearly independent since $\fG_{0,\beta}$
is 2-dimensional.  Suppose that 0 = a [u$_0$, A$_{\beta}$(u$_0$)] + b [u$_0$, B$_{\beta}$(u$_0$)] for some real numbers a,b.  Taking inner products with both A$_{\beta}$ and B$_{\beta}$ yields  C$\begin{pmatrix} a \\ b \end{pmatrix}$ = $\begin{pmatrix} 0 \\ 0 \end{pmatrix}$ , where C = $\begin{pmatrix} |$A$_{\beta}$(u$_0$)$|^2 & \langle$ A$_{\beta}$(u$_0$)~,~
B$_{\beta}$(u$_0$) $\rangle \\ \langle$ A$_{\beta}$(u$_0$)~,~B$_{\beta}$(u$_0$) $\rangle & |$B$_{\beta}$(u$_0$)$|^2 \end{pmatrix}$.  By 3) of (7.6) det C $>$ 0, and hence a = b = 0.
\end{proof}

    Recall that an element $\beta \in \Phi \cap \Lambda$ determines a real root vector $\tilde{H}_{\beta} \in \fH_{0}$, where $\fH = \fH_{0}^{\bc}$ is the Cartan subalgebra of $\fG$ that defines $\Phi$ and $\Lambda$.  See the discussion preceding (6.8).

\begin{proposition}  Let $\beta \in \Phi^{+} \cap \Lambda$.  Let u$_{\beta}$ be a nonzero element of
U$_{\beta}$, and consider ad u$_{\beta}$ : U$_{\beta} \rightarrow \br$-span\{$\tilde{H}_{\beta}$\}, where $\tilde{H}_{\beta} \in \fH_{0}$
is the real root vector determined by $\beta$.  Then

\indent 1)  Ker (ad u$_{\beta}$) = \{$\tilde{H}_{\beta}$(u$_{\beta}$)\}$^{\perp}$.

\indent 2)  ad u$_{\beta}$ : U$_{\beta} \rightarrow \br$-span\{$\tilde{H}_{\beta}$\} is surjective.
\end{proposition}

\begin{proof}  We recall from (6.10) that  ad u$_{\beta}$(U$_{\beta}$) $\subseteq \br$-span\{$\tilde{H}_{\beta}\}$.  Hence if u$_{\beta}' \in$ U$_{\beta}$, then ad u$_{\beta}$(u$_{\beta}'$) = 0 $\iff 0 = \langle$ [u$_{\beta}$ , u$_{\beta}'$], $\tilde{H}_{\beta}
\rangle = \langle \tilde{H}_{\beta}$(u$_{\beta}$), u$_{\beta}' \rangle$.  This proves 1), and 2) follows
immediately from 1) since dim U$_{\beta} \geq$ 2 by 4) of (6.3).
\end{proof}

\begin{proposition}  Let $\beta \in \Phi^{+} \cap \Lambda$.  Let u$_0 \in$ U$_0$ with  u$_0 \notin$
Ker (A$_{\beta}$) =  Ker (B$_{\beta}$) , and let u$_{\beta}$ be a nonzero element of U$_{\beta}$.  Then Ker (ad u$_0$) $\subset$
Ker (ad u$_{\beta}$) $\iff$ u$_{\beta} \in \br$-span\{(A$_{\beta}$(u$_0$), B$_{\beta}$(u$_0$)\}.
\end{proposition}

\begin{lem}  The map $\tilde{H}_{\beta} : U_{\beta} \rightarrow U_{\beta}$ is an isomorphism.
\end{lem}
\begin{proof}  Note that $\tilde{H}_{\beta}$ leaves U$_{\beta}$ invariant by 1) of (6.3) and the fact that $\tilde{H}_{\beta} \in \fH_{0}$.  It suffices to prove that $\tilde{H}_{\beta}$ is injective.  Let u$_{\beta}$ be a nonzero element of U$_{\beta}$ and write $u_{\beta} = v_{\beta} + v_{- \beta}$, where $v_{\beta} \in V_{\beta}$ and $v_{- \beta} \in V_{- \beta}$.  Then $\tilde{H}_{\beta}(u_{\beta}) = \beta(\tilde{H}_{\beta})(v_{\beta} - v_{- \beta}) = i \langle \tilde{H}_{\beta}, \tilde{H}_{\beta} \rangle (v_{\beta} - v_{- \beta})$ by the definition of $\tilde{H}_{\beta}$ in the discussion preceding (6.8).  It follows that $\tilde{H}_{\beta}(u_{\beta}) \neq 0$ since $\tilde{H}_{\beta} \neq 0$ and u$_{\beta} \neq 0$.
\end{proof}

\begin{proof}  We now complete the proof of the proposition.  From (7.7) and (7.8) it follows that  Ker (ad u$_0$) $\subseteq$ Ker (ad u$_{\beta}$) $\iff \br$-span \{A$_{\beta}$(u$_0$), B$_{\beta}$(u$_0$)\} = Ker (ad u$_0$)$^{\perp} \supseteq$ Ker (ad u$_{\beta}$)$^{\perp} =
\tilde{H}_{\beta}$(u$_{\beta}$). By the lemma above $\tilde{H}_{\beta} : U_{\beta} \rightarrow U_{\beta}$ is an isomorphism.  It therefore suffices to prove that $\tilde{H}_{\beta}$ leaves invariant $\br$-span \{A$_{\beta}$(u$_0$), B$_{\beta}$(u$_0$)\} $\subseteq$ U$_{\beta}$.  By (6.1)  we may write $\tilde{H}_{\beta}$ = i H for some element H $\in \fH_{\br}$.  Hence $\beta$(H) $= -$ i $\beta$($\tilde{H}_{\beta}$) $\in \br$  by (6.1). We compute $\tilde{H}_{\beta}$(A$_{\beta}$(u$_0$)) = A$_{\beta}$($\tilde{H}_{\beta}$(u$_0$)) + [$\tilde{H}_{\beta}$, A$_{\beta}$](u$_0$) =
[$\tilde{H}_{\beta}$, A$_{\beta}$](u$_0$) = $\beta$(H) B$_{\beta}$(u$_0$) by (7.1) and the fact that u$_0 \in$ U$_0$.  Similarly, we obtain $\tilde{H}_{\beta}$(B$_{\beta}$(u$_0$)) = $- \beta$(H) A$_{\beta}$(u$_0$).  Hence $\tilde{H}_{\beta}$ leaves invariant $\br$-span (A$_{\beta}$(u$_{0}$), B$_{\beta}$(u$_{0}$)).
\end{proof}

\section{Rational structures on $\fN_{0} = U \oplus ~\fG_{0}$}

    Let $\fH$ be a finite dimensional real Lie algebra, and let $\fB = \{\xi_{1},~ ... ~, \xi_{n}\}$ be a basis
of $\fH$ with rational structure constants; that is, [$\xi_{i} , \xi_{j}$] = \(\sum_{k=1}^{n} C^{k}_{ij} \xi_{k} \)
where the constants $\{C^{k}_{ij}\}$ lie in $\bq$.  The $\bq$-Lie algebra $\fH_{\bq}$ = $\bq$-span ($\fB$) is said to be a $\mathit{rational~ structure}$ for $\fH$.  A subspace $\fH'$ of $\fH$ is said to be a $\mathit{rational~ subspace}$ of $\fH$ (relative to $\fH_{\bq}$) if $\fH'$ admits a basis contained in $\fH_{\bq}$.

    Let $\fN_{0}$ be a real nilpotent Lie algebra (not necessarily 2-step nilpotent), and let N$_{0}$ denote the simply connected Lie group with Lie algebra $\fN_{0}$.  It is known that the exponential map exp : $\fN_{0} \rightarrow N_{0}$ is a diffeomorphism.  Let log : $N_{0} \rightarrow \fN_{0}$ denote the inverse of exp.  If $\fN_{0}'$ is a subalgebra of $\fN_{0}$, then $N_{0}' =$ exp($\fN_{0}'$) is a simply connected subgroup of N$_{0}$ with Lie algebra $\fN_{0}'$.

    We now consider rational structures on $\fN_{0}$, and we state some basic facts
without proof.  See [R2] or [CG] for further explanation.

    A subgroup $\Gamma_{0}$ of N$_{0}$ is called a $\mathit{lattice}$ if $\Gamma_{0}$ is discrete and the manifold of left cosets $\Gamma_{0} \backslash N_{0}$ is compact.  A criterion of Mal'cev [Ma] states that N$_{0}$ admits a lattice $\Gamma_{0} \iff \fN_{0}$ admits a rational structure $\fN_{0, \bq}=\bq$-span ($\fB$).  More precisely, if $\fN_{0}$ admits a rational structure $\fN_{0, \bq}$, and L is a vector lattice in $\fN_{0,\bq}$, then exp(L) generates a lattice $\Gamma_{0}$ in N$_{0}$, and $\bq$-span (log $\Gamma_{0}$) = $\fN_{0, \bq}$.  If L$_1$ and L$_2$ are vector lattices in $\fN_{0, \bq}$ with corresponding lattices $\Gamma_{1}$ and $\Gamma_{2}$ in N$_{0}$, then $\Gamma_{1}~ \cap ~\Gamma_{2}$ has finite index in both $\Gamma_{1}$ and $\Gamma_{2}$.  Conversely, if $\Gamma_{0}$ is a lattice in N$_{0}$, then $\fN_{0, \bq} = \bq$-span (log $\Gamma_{0}$) is a rational structure for $\fN_{0}$.  Most nilpotent Lie algebras $\fN_{0}$ , even 2-step nilpotent ones, do not admit a rational structure.  For specific examples and more precise statements see [R2] or [Eb 1].

    Let N$_{0}$ be a simply connected nilpotent Lie group with Lie algebra $\fN_{0}$.  Suppose that  N$_{0}$ admits a lattice $\Gamma_{0}$, and equip $\fN_{0}$ with the corresponding rational structure $\fN_{0, \bq}=\bq$-span (log $\Gamma_{0}$).  Let $\fN_{0}'$ be a rational subalgebra of $\fN_{0}$,
 and let $N_{0}' =$ exp($\fN_{0}'$) be the simply connected Lie subgroup of N$_{0}$ with Lie algebra $\fN_{0}'$.  If $\Gamma_{0}' = \Gamma_{0}~ \cap N_{0}'$ , then $\Gamma_{0}'$ is a lattice in $N_{0}'$.  See [CG,Theorem 5.1.11] for a proof.
 \newline

$\mathit{Existence~ of~ rational~ structures~ on~ \fN_{0} = U \oplus \fG_{0}}$

    Let $\fG_0$ be a compact, semisimple real Lie algebra, and let U be a finite dimensional real $\fG_{0}$-module.  Let $\fN_{0} = U \oplus \fG_{0}$ be equipped with the 2-step nilpotent Lie algebra structure defined in (2.3).  Let $\fH_{0}$ be a maximal abelian subalgebra of $\fG_{0}$ , and let $\fH = \fH_{0}^{\bc}$ be the corresponding Cartan subalgebra of $\fG = \fG_{0}^{\bc}$.  Let
U = U$_0$ + \(\sum_{\lambda \epsilon \Lambda^{+}}\) U$_\lambda$ (direct sum) be the weightspace decomposition of U determined by $\fH$ as in (6.4).

    To construct a rational structure on $\fN_{0}  = U \oplus \fG_{0}$  we use the following result of Raghunathan.  For a generalization of this result due to D. Morris (= D. Witte) see [Mo].

\begin{theorem} ([R1])  Let $\fG_{0}$ be a compact, semisimple real Lie algebra, and let $\fC_{0}$ be a compact Chevalley basis of $\fG_{0}$.  Let U be a finite dimensional, real $\fG_{0}$-module.  Then there exists a basis $\fB_{0}$ of U such that $\fG_{0}$ leaves invariant $\bq$-span ($\fB_{0}$).
\end{theorem}

$\mathit{Rationality~ of~ weight~ spaces}$

    As a consequence of the preceding result we obtain

\begin{proposition}  Let $\fC_{0}$ be the compact Chevalley basis of $\fG_{0}$ as defined in (3.3).  Let $\fB_{0}$ be any basis of U such that $\fG_{0}$ leaves invariant $\bq$-span($\fB_{0}$).  Let $\ip$ be an inner product on $\fN_{0} = U \oplus \fG_{0}$ satisfying the conditions of (2.4) such that $\langle \xi, \xi' \rangle \in \bq$ for all $\xi, \xi' \in \fB_{0} \cup \fC_{0}$.  Then

\indent 1)  $\fL_{0} = \fB_{0} \cup \fC_{0}$ is a basis of $\fN_{0} = U \oplus \fG_{0}$ that has rational structure constants.  Furthermore, $\fB_{0}$ may be chosen so that  $\fL_{0}$ is an orthogonal basis of $\{\fN_{0}, \langle , \rangle \}$.

\indent 2)  The weight spaces U$_0$ and U$_{\lambda}, \lambda \in \Lambda^{+}$, are rational with respect to the rational structure $\fN_{0, \bq}=\bq$-span($\fL_{0}$).
\end{proposition}

$\mathit{Remark}$  We call $\fN_{0, \bq}=\bq$-span($\fL_{0}$) a $\mathit{Chevalley~ rational~ structure}$ for $\fN_{0}$.  The most interesting case occurs when $\fB_{0}$ is the union of bases for the weight spaces U$_0$, U$_{\mu}$, which is possible by 2) of the proposition above.

\begin{proof}

    1) By (2.5) it suffices to consider the case that U is an irreducible $\fG_0$-module.  We first prove
the existence of an inner product $\ip$ on $\fN$ such that $\langle \xi~,~\xi' \rangle \in \bq$ for all $\xi~,~\xi' \in \fB_{0} \cup \fC_{0}$ . It suffices to define $\ip$ independently on U and $\fG_0$ so that $\langle \xi~,~\xi' \rangle \in \bq$ for all $\xi~,~\xi' \in \fB_{0}$ and $\langle \xi~,~\xi' \rangle \in \bq$ for all $\xi~,~\xi' \in \fC_{0}$.  We then let $\ip$ denote the inner product on $\fN_{0} = U \oplus \fG_0$ that makes U and $\fG_0$ orthogonal and restricts to the inner products on these two factors.

    On U we choose a $\fG_0$-invariant inner product $\ip$ such that $\langle \xi~,~\xi' \rangle \in \bq$ for all $\xi~,~\xi' \in \fB_{0}$.  This can be done by Lemma 3 in the proof of (2.7) of [Eb 2].  On $\fG_0$ we may define the inner product $\langle Z, Z' \rangle = -$ trace (ZZ$'$), where we regard Z, Z$'$ as elements of End(U).  It is easy to check that $\ip$ is $\fG_0$-invariant.  By hypothesis $\fG_0$ leaves invariant $\bq$-span($\fB_{0}$), or equivalently, every element of $\fG_0$ has a matrix with rational entries with respect to $\fB_{0}$.  It follows that $\langle Z, Z' \rangle \in \bq$ if Z, Z$'  \in \fG_0$ .  We could also define $\ip = -$ B$_0$ on $\fG_0$, where B$_0$ denotes the Killing form of $\fG_0$.  In this case $\langle Z, Z' \rangle \in \bz$ if Z, Z$'$ $\in \fC_0$ by (3.3).

    If $\fL_{0} = \fB_{0} \cup \fC_0$, then by the discussion above $\langle \xi~,~\xi' \rangle \in \bq$ for all $\xi~,~\xi' \in \fL_{0}$.  If X,Y $\in \fB_{0}$ then [X,Y] $\in \fG_0$ is orthogonal to $\fB_{0} \subset$ U and $\langle [X,Y] , Z \rangle = \langle Z(X), Y \rangle \in \bq$ for all Z $\in \fC_0$ since the elements of $\fG_0$ leave invariant $\bq$-span($\fB_{0}$) by hypothesis.  Hence [$\fL_{0}, \fL_{0}$] $\subset
\bq$-span ($\fC_0$) $\subset \bq$-span ($\fL_{0}$) by the lemma in (6.9).

    Now let $\fB_{0} = \{u_{1},~ ...~,u_{q} \}$ be any basis of U such that $\fC_{0}$ leaves invariant $\bq$ - span ($\fB_{0}$).  Use the Gram - Schmidt process  to define inductively an orthogonal set $\fB_{0}' = \{u_{1}',~ ...~, u_{q}'  \}$ given by $u_{1}' = u_{1}$ and $u_{r}' = u_{r} -$ \(\sum_{k=1}^{r-1} \{\langle u_{r}, u_{k}' \rangle / \langle u_{k}' , u_{k}' \rangle \} u_{k}' \) for $2 \leq r \leq q$.  The set $\fB_{0}'$ lies in $\bq$ - span ($\fB_{0}$) by the hypotheses on $\langle , \rangle$, and hence $\bq$ - span ($\fB_{0}'$) $ = \bq$ - span ($\fB_{0}$).  If $\fL_{0}' = \fB_{0}'~ \cup~ \fC_{0}$, then $\fL_{0}'$ is a basis of $\fN_{0}$ with rational structure constants by the discussion above.  Moreover, $\fL_{0}'$ is orthogonal by (3.4) and the orthogonality of U and $\fG_{0}$.

    2)  We prove rationality only for the subspaces U$_{\lambda}, \lambda \in \Lambda ^{+}$, since the proof for U$_0$ is similar.  Let $\lambda \in \Lambda^{+}$ be given, and let $\fH_0$ be a maximal abelian subalgebra of $\fG_0$.  Let H$_0 \in \fH_0$ be an element such that $\lambda$(H$_0$) is nonzero and the numbers $\{ \lambda(H_0)^{2} : \lambda \in \Lambda^{+} \}$ are all distinct.
The set of vectors H$_0 \in \fH_0$ that satisfy these conditions is a dense open subset of $\fH_0$, and hence we may impose the additional condition that H$_0 \in \bq$-span $\{\tilde{\tau}_{\alpha} : \alpha \in \Delta \}$ $\subseteq \bq$-span ($\fC_0$).  Note that $\lambda$(H$_0$)$^2 \in \bq$ since $\lambda$($\tilde{\tau}_{\alpha}) \in$ i $\bz$ for all $\alpha \in \Delta$.  Since H$_{0}^{2} \in$ End(U) leaves invariant $\bq$-span($\fB_{0}$) it follows that the matrix A$_0$ of H$_{0}^{2} $ relative to $\fB_{0}$ has entries in $\bq$.

    By 5) of (6.3) H$_0^2 = - \lambda(H_0)^2$ Id on U$_{\lambda}$.  If T$_{\lambda}$ : U $\rightarrow$ U is given by T$_{\lambda}$ = $H_0^2 + \lambda(H_0)^2$ Id, then by the conditions on H$_0$ it follows that U$_{\lambda}$ = Ker(T$_{\lambda}$).  Moreover, T$_{\lambda}$
leaves invariant U$_{\bq} = \bq$-span($\fB$).  The $\bq$-rank of T$_{\lambda}$ acting on U$_{\bq}$ is the same as the $\br$-rank of T$_{\lambda}$ acting on U ; both are the largest integer k such that the rational matrix A$_0 + \lambda(H_0)^2$ Id has a k x k minor with nonzero determinant.  Hence dim$_{\bq}U_{\bq ,\lambda}$, where $U_{\bq, \lambda}$ is the $\bq$-kernel of T$_{\lambda}$ acting on U$_{\bq}$, is the same as dim$_{\br}U_{\br, \lambda}$, where $U_{\br, \lambda}$ is the $\br$-kernel of T$_{\lambda}$ acting on U.  By the Gram-Schmidt process we can find an orthogonal $\bq$-basis $\fB' \subset U_{\bq}$ for $U_{\bq, \lambda}$.  The set $\fB'$ is linearly independent over $\br$ since it is orthogonal, and hence $\fB'$ is an $\br$-basis for U$_{\br, \lambda} =$ Ker(T$_{\lambda}$).
\end{proof}

\section{Admissible abstract weights}

    We introduce a notion of admissibility for an abstract weight that has already been relevant for (6.12) and will be used again in section 10 to construct rational, totally geodesic subalgebras of $\fN_{0} = U \oplus \fG_{0}$.

    Let $\fG$ be a complex semisimple Lie algebra.  Let $\fH$ be a Cartan subalgebra with roots $\Phi$, and let $\Delta = \{\alpha_{1},~ ...~ , \alpha_{n} \} \subset \Phi^{+}$ be a base of simple roots.  For $\lambda \in$ Hom ($\fH,\bc$) and $\alpha \in \Phi$ we introduce the notation $\langle \lambda,\alpha \rangle :~= \lambda(\tau_{\alpha})$.  Note that $\langle , \rangle$ is $\br$ - linear in the first entry but not in the second.

    Let $\tilde{\Lambda} = \{\lambda \in$ Hom($\fH,\bc$) :
$ \langle \lambda,\alpha \rangle \in \bz$ for all $\alpha \in \Phi$ \}, the lattice of $\mathit{abstract~ weights}$ determined by $\fH$.  We let $\{\omega_{1},~ ...~ , \omega_{n} \}$ denote the fundamental weights defined by the condition $\delta_{ij} = \langle \omega_{i}~,~\alpha_{j} \rangle$.

    Let $\mu \in \tilde{\Lambda}$ be an abstract weight of a complex semisimple Lie algebra $\fG$.   We say that  $\mu$ is $\mathit{admissible}$  if $2\mu + p~\alpha$ does not lie in $\Phi$ for any integer p.  Since $\tilde{\Lambda}$ and $\Phi$ are invariant under the Weyl group W  it follows immediately from the definition that the set of admissible weights is invariant under W.  Clearly the inadmissible weights are also invariant under W.  We shall classify the inadmissible abstract weights by picking out a special element in each  orbit of W.

    As we have seen already in (6.12) admissible abstract weights are relevant to describing the image of ad $\xi : \fN \rightarrow \fN$ for a noncentral element $\xi$ of $\fN$.  In the next section we shall see that admissible abstract weights are also relevant to describing totally geodesic subalgebras of $\fN$.

    Most abstract weights $\mu$ are admissible, and we list the exceptions below.  We consider separately the cases that $\fG$ is simple or a direct sum of at least two simple ideals.  We also distinguish between the cases that the integer p above is even or odd.  We use the standard notation and classification for complex simple Lie algebras.
\newline

$\mathit{Notation}$  We denote the highest long root by $\beta_{max}$ and the highest short root by  $\beta_{min}$.  A table of the roots  $\beta_{max}$ and  $\beta_{min}$ may be found in [Hu, p. 66] and also in the appendix.
\newline

$\mathit{Simple~ Lie~ algebras}$

    We describe the abstract weights $\mu$ for a simple complex Lie algebra $\fG$ that are not admissible.

\begin{proposition}  Let $\fG$ denote a complex simple Lie algebra.  Let $\mu \in \tilde{\Lambda}$ be an abstract weight such that  $2\mu = (2k+1)~\alpha + \beta$ for some elements $\alpha,\beta \in \Phi$ and some integer k.  Then all Weyl group orbits W($\mu$) are exhausted by the following possibilities :

A$_{n}$ \hspace{.5in} $\mu = (k+1)(\omega_{1} + \omega_{n})  = (k+1)\beta_{max} = (k+1)\beta_{min}$.

A$_{3}$ \hspace {.5in}   $\mu = \omega_{2} + k(\omega_{1} + \omega_{3}) =   \omega_{2} + k\beta_{max} =\omega_{2} + k\beta_{min}$.

B$_{n},  n \geq 3$

    \hspace{.7in}     $\mu = (k+1)\omega_{2} = (k+1)\beta_{max}$.

    \hspace{.7in}  $\mu = (k+1)\omega_{1} = (k+1)\beta_{min}$.

    \hspace{.7in}  $\mu = \omega_{1} + k\omega_{2} = \beta_{min} + k\beta_{max}$.

B$_{4}$ \hspace{.5in}   $\mu = \omega_{4} + k\omega_{2} = \omega_{4} + k\beta_{max}$.

B$_{3}$ \hspace{.5in}   $\mu = \omega_{3} + k\omega_{2} = \omega_{3} + k\beta_{max}$.

    \hspace{.7in}           $\mu = \omega_{3} + k\omega_{1} = \omega_{3} + k\beta_{min}$.

C$_{n}, n \geq 2$

    \hspace {.7in}  $\mu = 2(k+1)\omega_{1} = (k+1)\beta_{max}$.

    \hspace{.7in}   $\mu = \omega_{1} + k \omega_{2} = \frac{1}{2} \beta_{max} + k \beta_{min}$

    \hspace{.7in} $\mu = (k+1)\omega_{2} = (k+1)\beta_{min}$.

    \hspace{.7in} $\mu = \omega_{2} + 2k\omega_{1} = \beta_{min} + k\beta_{max}$.

D$_{n}, n \geq 4$

    \hspace {.7in}  $\mu = (k+1)\omega_{2} = (k+1)\beta_{max} = (k+1)\beta_{min}$.

    \hspace {.7in}  $\mu = \omega_{1} + k\omega_{2} = \omega_{1} + k\beta_{max} =  \omega_{1} + k\beta_{min}$.

D$_{4}$ \hspace{.5in}   $\mu = \omega_{3} + k\omega_{2} = \omega_{3} + k\beta_{max} =  \omega_{3} + k\beta_{min}$.

    \hspace{.7in}   $\mu = \omega_{4} + k\omega_{2} = \omega_{4} + k\beta_{max} =  \omega_{4} + k\beta_{min}$.

E$_{6}$ \hspace{.5in}   $\mu = (k+1)\omega_{2} = (k+1)\beta_{max} = (k+1)\beta_{min}$.

E$_{7}$ \hspace{.5in}   $\mu = (k+1)\omega_{1} = (k+1)\beta_{max} = (k+1)\beta_{min}$.

E$_{8}$ \hspace{.5in}   $\mu = (k+1)\omega_{8} = (k+1)\beta_{max} = (k+1)\beta_{min}$.

F$_{4}$ \hspace{.5in}       $\mu = (k+1)\omega_{1} = (k+1)\beta_{max}$.

    \hspace{.7in}   $\mu = (k+1)\omega_{4} = (k+1)\beta_{min}$.

    \hspace{.7in}   $\mu = \omega_{4} + k\omega_{1} = \beta_{min} + k\beta_{max}$.

G$_{2}$ \hspace{.5in}   $\mu = (k+1)\omega_{2} = (k+1)\beta_{max}$.

    \hspace{.7in} $\mu = (k+1)\omega_{1} = (k+1)\beta_{min}$.

    \hspace{.7in} $\mu = \omega_{1} +  k\omega_{2} = \beta_{min} + k\beta_{max}$.

\end{proposition}
\begin{proof}  See the Appendix.
\end{proof}

\begin{proposition}  Let $\fG$ denote a complex simple Lie algebra.  Let $\mu \in \tilde{\Lambda}$ be an abstract weight such that  $2\mu = (2k)~\alpha + \beta$ for some elements $\alpha,\beta \in \Phi$ and some integer k.  Then all Weyl group orbits W($\mu$) are exhausted by the following possibilities :

A$_{1}  \hspace{.5in}   \mu = 2k + 1$.

C$_{n}, n \geq 3$

    \hspace{.7in}   $\mu = \omega_{1} + k\alpha_{1}$.

    \hspace{.7in}   $\mu = \omega_{1} + k\alpha_{2}$.

    \hspace{.7in}   $\mu = \omega_{1} + k\alpha_{n}$.

    \hspace{.7in}   $\mu = (2k+1)\omega_{1}  = \omega_{1} + k\beta_{max}$.

C$_{2}  \hspace{.55in} \mu = \omega_{1} + k\alpha_{1}$.

    \hspace{.7in}   $\mu = \omega_{1} + k\alpha_{2}$.

    \hspace{.7in}   $\mu = (2k+1)\omega_{1}  = \omega_{1} + k\beta_{max}$.

\end{proposition}
\begin{proof}  See the Appendix.
\end{proof}

$\mathit{Semisimple~ Lie~ algebras}$

    We describe the abstract weights $\mu$ for a semisimple complex Lie algebra $\fG$ that are not admissible.  We first explain how to reduce consideration to the case of a simple complex Lie algebra.

    Let $\fG$ be a complex semisimple Lie algebra, and let $\fG = \fG_{1} \oplus~ ...~ \oplus~ \fG_{N}$ be its decomposition into a direct sum of simple ideals \{$\fG_i$\}.  Any Cartan subalgebra $\fH$ of $\fG$ may be written $\fH = \fH_{1} \oplus~ ...~ \oplus~ \fH_{N}$, where $\fH_i$ is a Cartan subalgebra of $\fG_i$ for 1 $\leq$ i $\leq$ N.  For convenience we let $\fH$* and $\fH_{i}$* denote Hom($\fH,\bc$) and Hom($\fH_{i},\bc$) , 1 $\leq$ i $\leq$ N. We may identify $\fH_{i}^{*}$ with a subspace of $\fH$* by setting $\lambda_{i} \equiv 0$ on $\fH_{j}$ for each $\lambda_{i} \in \fH_{i}^{*}$ and all
j $\neq$ i, 1 $\leq$ i,j $\leq$ N.  With this identification it is easy to see that that $\fH^{*} = \fH_{1}^{*} \oplus~ ... ~\fH_{N}^{*}$, (direct sum).

    We regard the Weyl group W$_{i} \subset \fH_{i}$* as a normal subgroup of W $\subset \fH$* by letting the elements of  W$_{i}$ act as the identity on $\fH_{j}$* for all j $\neq$ i.  It is easy to see that W = W$_{1}$ x ... x W$_{N}$, (direct product).

    We let $\tilde{\Lambda} \subset \fH$* and $\tilde{\Lambda_{i}} \subset \fH_{i}^{*}$ denote the abstract weights defined by $\fH$ and $\fH_{i}$; that is, $\tilde{\Lambda} = \{\lambda \in \fH^{*} :  \langle \lambda,\alpha \rangle \in \bz$ for all $\alpha \in \Phi$ \} with
a similar definition for $\tilde{\Lambda}_{i}$. It is easy to see that $\tilde{\Lambda} =
\tilde{\Lambda_{1}} \oplus~ ... \oplus ~\tilde{\Lambda_{N}}$.

    Let $\Phi \subset \fH$* and $\Phi_i \subset \fH_{i}^{*}$ denote the roots of $\fH$ and $\fH_i$ respectively.  The root space decompositions of each factor $\fG_{i}$ combine to yield a root space decomposition of $\fG$, and it follows that $\Phi = \Phi_{1} \cup~ ...~ \cup~ \Phi_{N}$.  Moreover, let $\Delta$ be a base of simple roots of $\Phi$.  Then  $\Delta = \Delta_{1} \cup~ ...~ \cup~ \Delta_{N}$, where each $\Delta_{i}$ is a base for the simple roots of $\Phi_{i}$.

    Let p$_i$ : $\fH$* $\rightarrow \fH_{i}^{*}$ denote the surjective homomorphism such that p$_i(\lambda)$ is the restriction of $\lambda$ to $\fH_i$ for 1 $\leq$ i $\leq$ N.  Note that if $\alpha_{j} \in \Phi_{j}$, then p$_i(\alpha_{j}) = \delta_{ij} \alpha_{j}$ for 1 $\leq$ i,j $\leq$ N.

    We are now ready to state the extensions of (9.1) and (9.2) to the case that $\fG$ is a complex semisimple Lie algebra.  It is evident from the next three results that if $\fG$ is semisimple but not simple, then for abstract weights  $\mu \in \tilde{\Lambda}$ it is much more difficult than in the case $\fG$ is simple to find an element $\alpha \in \Phi$ and an integer p such that $2\mu + p\alpha \in \Phi$.

 \begin{proposition}  Let $\fG$ be a complex semisimple Lie algebra with Cartan subalgebra $\fH$, roots $\Phi \subset \fH$* and  base $\Delta$ for $\Phi$ as above.  Let $\mu \in \tilde{\Lambda}$ be an abstract weight and suppose that $2\mu = (2k+1)\alpha + \beta$ for some elements $\alpha, \beta \in \Phi$ and some integer k.  Then one of the following cases occurs :

 \indent    1)  There exists an integer j with 1 $\leq$ j $\leq$ N such that $\mu \in \tilde{\Lambda_{j}}, \alpha \in \Phi_{j}$ and $\beta \in \Phi_{j}$.

 \indent    2)  There exist distinct integers i,j with 1 $\leq$ i,j $\leq$ N such that

  \indent \indent a)  $p_{k}(\mu) = 0$ for k $\neq$ i, k $\neq$ j.

 \indent \indent    b)  The root $\alpha$ lies in $\Phi_{i}$ and $\fG_{i} = C_{n}$ for some integer n $\geq 2$.  Moreover, there exists an element w$_{i}$ of W$_{i}$ such that w$_{i}($p$_{i}(\mu)) = (2k+1)\omega_{1}$ and w$_{i}($p$_{i}(\alpha)) = \beta_{max}$.

 \indent \indent    c) Either $\fG_{j} = A_{1}$  or $\fG_{j} = C_{n}$ for some integer n $\geq 2$.  In the first case p$_{j}(\mu) = 2k+1$ or $2k -1$.  In the second case there exists an element w$_{j}$ of W$_{j}$ such that  w$_{j}($p$_{j}(\mu)) = \omega_{1}$

 \end{proposition}

$\mathit{Remark}$  It is easy to check that each of the cases above actually occurs.  In case 1) we are reduced to the case that $\fG$ is simple, and we may apply (9.1) and (9.2).

\begin{proof}  Since $\Phi = \Phi_{1}~ \cup~ ... ~\cup \Phi_{N}$ there exist integers i,j such that $\alpha \in \Phi_{i},\beta \in \Phi_{j}$ and $2\mu - (2k+1)\alpha = \beta$.  If k $\neq$ i, k $\neq$ j, then $0 = p_{k}(\beta) = p_{k}(2\mu - (2k+1)\alpha) = 2p_{k}(\mu)$.  If i = j, then we obtain 1).  If i $\neq$ j, then by 2a) we may write $\mu = \mu_{i} + \mu_{j}$, where $\mu_{i} \in \tilde{\Lambda_{i}} $ and $\mu_{j} \in \tilde{\Lambda_{j}} $.  From the equation $\beta = 2\mu - (2k+1)\alpha = (2\mu_{i} - (2k+1)\alpha) + 2\mu_{j}$ and the fact that $\alpha \in \Phi_{i}$ we obtain the equations i) $2\mu_{i} - (2k+1)\alpha = 0$ and ii)  $2\mu_{j} = \beta$.

    By (i) and (9.2) with $\alpha = \beta$ we see that either $\fG_{i} = A_{1}$ or $\fG_{i} = C_{n}$ for some integer n $\geq 2$.  If $\fG_{i} = A_{1}$, then $\mu_{i}$ and $\alpha$ would be integers, but the equation (i) has no nonzero integer solutions.  Therefore $\fG_{i} = C_{n}$ for some integer n $\geq 2$. The equations i) and ii) together with (9.2) now yield 2b) and 2c).

\end{proof}

 \begin{proposition}  Let $\fG$ be a complex semisimple Lie algebra with Cartan subalgebra $\fH$, roots $\Phi \subset \fH$* and base $\Delta$ for $\Phi$ as above.  Let $\mu \in \tilde{\Lambda}$ be an abstract weight such that 2$\mu = (2k)\alpha + \beta$ for some elements $\alpha, \beta \in \Phi$ and some nonzero integer k.  Then one of the following cases occurs :

  \indent   1)  There exists an integer j with 1 $\leq$ j $\leq$ N such that $\mu \in \tilde{\Lambda_{j}}, \alpha \in \Phi_{j}$ and $\beta \in \Phi_{j}$.

 \indent    2)  There exist distinct integers i,j with 1 $\leq$ i,j $\leq$ N such that

  \indent \indent a)  $p_{k}(\mu) = 0$ for k $\neq$ i, k $\neq$ j.

 \indent \indent    b)  The root $\alpha$ lies in $\Phi_{i}$, and  p$_{i}(\mu) = k\alpha$.

 \indent \indent    c) Either $\fG_{j} = A_{1}$  or $\fG_{j} = C_{n}$ for some integer n $\geq 2$.  In the first case p$_{j}(\mu) = \pm 1$.  In the second case there exists an element w$_{j}$ of W$_{j}$ such that
 w$_{j}($p$_{j}(\mu)) = \omega_{1}$

  \end{proposition}

\begin{proof} Since  $\Phi = \Phi_{1} \cup~ ... ~\cup \Phi_{N}$ we may choose integers i,j such that $\alpha \in \Phi_{i}$ and $\beta \in \Phi_{j}$.  If i = j, then we are in case 1).  If i $\neq$ j, then 2a) follows as in the proof of the preceding proposition.  A trivial modification of that proof also shows that $\beta =
2p_{j}(\mu)$ and $0 = 2p_{i}(\mu) - (2k)\alpha$.  From (9.2) we obtain 2b) and 2c).
\end{proof}

\begin{proposition}   Let $\fG$ be a complex semisimple Lie algebra with Cartan subalgebra $\fH$, roots $\Phi \subset \fH$* and base $\Delta$ for $\Phi$ as above.  Let $\mu \in \tilde{\Lambda}$ be an abstract weight such that 2$\mu \in \Phi$.  Then  $\mu \in \tilde{\Lambda_{j}}$ for some integer j with 1 $\leq$ j $\leq$ N.  Either $\fG_{j} = A_{1}$  or $\fG_{j} = C_{n}$ for some integer n $\geq 2$.  In the first case $\mu = \pm 1$.  In the second case there exists an element w of W such that w$(\mu) = \omega_{1}$.
 \end{proposition}

 \begin{proof}  This follows immediately from the proof of c) in the previous proposition.
 \end{proof}

\section{Totally geodesic subalgebras}

    Let N$_{0}$ be a simply connected nilpotent Lie group with a left invariant metric $\ip$, and let $\fN_{0}$ denote the Lie algebra of N$_{0}$.  The Lie group exponential map exp : $\fN_{0} \rightarrow N_{0}$ is a diffeomorphism, and if $\fN_{0}$* is a subalgebra of $\fN_{0}$, then N$_{0}$* = exp($\fN_{0}$*) is a simply connected subgroup of N$_{0}$ with Lie algebra $\fN_{0}$*.  A subalgebra $\fN_{0}$* of N$_{0}$ is said to be $\mathit{totally~ geodesic}$ if the subgroup N$_{0}$* is a totally geodesic subgroup of N.   It is not difficult to prove the following :

\begin{lem}  A subalgebra $\fN_{0}$* of $\fN_{0}$ is totally geodesic $\iff$ ${\bigtriangledown}_{X}Y \in \fN_{0}$* whenever X,Y $\in \fN_{0}$*.
\end{lem}

    In the case that $\fN_{0} = U \oplus~ \fG_0$ is a 2-step nilpotent Lie algebra as defined in (2.3)  we now derive a sufficient condition for a subalgebra $\fN_{0}$* to be totally geodesic.

\begin{proposition}   A subalgebra $\fN_{0}$* of  $\fN_{0} = U \oplus~ \fG_0$ is totally geodesic if the following conditions are satisfied :

\indent 1)  $\fN_{0}$* = U* $\oplus~ \fG_{0}$*, where U* = $\fN_{0}$* $\cap~$ U and $\fG_{0}$* = $\fN_{0}$* $\cap~ \fG_{0}$.

\indent 2)  [U*, U*] $\subseteq \fG_0$*.

\indent 3)  Z(U*) $\subseteq$ U* for all Z $\in \fG_0$*.
\end{proposition}

\begin{proof}  Let $\fN_{0}$* satisfy 1), 2) and 3).  We write U* = U$_0$* $\oplus$ U$_1$*, where U$_0$* = $\{u \in $ U* : Z(u) = 0 for all Z $\in \fG_0$* \} and U$_1$* is the orthogonal complement in U* of U$_0$*.  It is not difficult to see that $\fZ_{0}$* = U$_{0}$* $\oplus~ \fG_{0}$* is the center of $\fN_{0}$*, and it follows that U$_1$* is the orthogonal complement of $\fZ_{0}$* in $\fN_{0}$*.

    From the formulas for the covariant derivative of a left invariant metric on a Lie group it is easy to see (e.g. (2.1) and (2.2) of [E3]) the following :

\indent a) $\bigtriangledown_{u} \xi = \bigtriangledown_{\xi}u = 0$ for all u $\in$ U$_0$ and all $\xi \in \fN_{0}$*.

\indent b) $\bigtriangledown_{X}Y = \frac{1}{2}~[X,Y]$ for all X,Y $\in$ U$_1$*.

\indent c) $\bigtriangledown_{Z}Z* = 0$ for all Z , Z* $\in \fG_0$*.

\indent d) $\bigtriangledown_{X}Z = \bigtriangledown_{Z}X = - \frac{1}{2}~Z(X)$ for all X $\in$ U$_1$* and all Z $\in \fG_0$*.

Since conditions 1), 2) and 3) hold for $\fN_{0}$* it follows from the lemma above and a) through d) that $\fN_{0}$* is a totally geodesic subalgebra.\end{proof}


    We now construct some examples of totally geodesic subalgebras.  We first construct totally geodesic subalgebras $\fN_{0}(\lambda)$ for weights $\lambda \in \Lambda$ such that 2$\lambda \notin \Phi$ and $\fN_{0}(\lambda,\beta)$ for pairs $(\lambda,\beta) \in \Lambda~$ x$~ \Phi$ such that  2$\lambda + $p$\beta \notin \Phi$ for any integer p.  These weights  have been classified in (9.1) and (9.2).

    The subalgebras $\fN_{0}(\lambda)$ and $\fN_{0}(\lambda,\beta)$ have centers of dimension 1 and 4 respectively.  The center of $\fN_{0}(\lambda,\beta)$ is a 4-dimensional subalgebra of $\fG_0$ that is Lie algebra isomorphic to the quaternions $\bh$ with the Lie algebra structure given by [x,y] = xy $-$ yx for x,y $\in \bh$.
\newline

$\mathit{Notation}$ We define or recall some notation and basic information that will be used in the next three results.

    For $\beta \in \Phi$ let $\fG_{0,\beta} = (\fG_{\beta} \oplus \fG_{-\beta}) \cap \fG_{0}$, a 2-dimensional subspace of $\fG_{0}$ by (6.5).

    Let $\fC_0$ denote the compact Chevalley basis of $\fG_{0}$ defined in the discussion preceding (3.3).  By (8.1) there exists a basis $\fB_{0}$ of U such that $\fC_{0}$ leaves invariant $\bq$-span($\fB_{0}$).  By 2) of (8.2) we may assume, without altering $\bq$-span($\fB_{0}$), that $\fB_{0}$ is a union of bases $\fB_{0}$* of U$_{0}$ and $\fB_{\mu}$ of U$_{\mu}, \mu \in \Lambda$. Let $\fN_{0, \bq}=\bq$-span($\fB_{0}~ \cup~ \fC_0$) denote the corresponding Chevalley rational structure for $\fN$.

    We let $\ip$ denote an inner product on $\fN =$ U $\oplus~ \fG_0$ satisfying the conditions of (2.4) such that $\langle \xi,\xi' \rangle \in \bq$ for all $\xi,\xi' \in \fB_{0}~  \cup~ \fC_0$.  The existence of such an inner product $\ip$ is part of the proof of 1) in (8.2).  We require further that $\ip = -$ B$_0$ on $\fG_0$, where B$_0$ denotes the Killing form on $\fG_0$. This allows us to use (6.8) and (6.9).

    Let $\tilde{\tau}_{\alpha} = i \tau_{\alpha}$ for $\alpha \in \Phi$.  Recall from the discussion in (6.8) $-$ (6.10) that for each $\lambda \in \Lambda$ there is a real weight vector $\tilde{H}_{\lambda}
 \in \fH_0 =$ i $\fH_{\br}$  defined by the condition $\langle \tilde{H}, \tilde{H}_{\lambda} \rangle  = -$ i $\lambda(\tilde{H})$ for all $\tilde{H} \in \fH_0$.  In particular, if $\ip = -$ B$_0$ on $\fG_0$, then $\tilde{H}_{\lambda}$ lies in $\bq$-span \{$\tilde{\tau}_{\alpha} : \alpha \in \Delta$ \} $\subset \bq$-span($\fC_0$).
 \newline

$\mathit{Rational~ totally~ geodesic~ subalgebras~ with~ 1-dimensional~ center}$

\begin{proposition}  For $\lambda \in \Lambda \subset$ Hom($\fH,\bc$) let $\fN_{0}(\lambda) = U_{\lambda} \oplus~ \br H_{\lambda}$.  If 2$\lambda \notin \Phi$  , then $\fN_{0}(\lambda)$ is a totally geodesic subalgebra of $\fN_{0}$ with a 1-dimensional center that is rational with respect to the Chevalley rational structure $\fN_{0, \bq}$ on $\fN_{0}$.
\end{proposition}

\begin{proof}  From 1) of (6.10) it follows that [$\fN_{0}(\lambda)$, $\fN_{0}(\lambda)$] $\subset \br \tilde{H}_{\lambda}$, and hence $\fN_{0}(\lambda)$ is a subalgebra of $\fN_{0}$.  By 1) of (6.3) $\tilde{H}_{\lambda}$(U$_{\lambda}$) $\subset$ U$_{\lambda}$, and it now follows from (10.1) that $\fN_{0}(\lambda)$ is a totally geodesic subalgebra of $\fN$.

    Clearly the center of $\fN_{0}(\lambda)$ is spanned by $\tilde{H}_{\lambda}$ and hence is 1-dimensional.  It remains only to prove that $\fN_{0}(\lambda)$ is rational with respect to the Chevalley rational structure $\fN_{0,\bq}$.  By hypothesis the basis $\fB_{\lambda}$ of U$_{\lambda}$ is a subset of $\fB_{0}$, and $\tilde{H}_{\lambda} \in \bq$-span ($\fC_0$) by (6.9).  Hence
$\fB_{\lambda} \cup$ \{$\tilde{H}_{\lambda}$ \} is a basis of $\fN_{0}(\lambda)$ that is contained in $\fN_{0, \bq}$. In particular, for any elements X$_{\lambda}$ , Y$_{\lambda}$ of $\fB_{\lambda}$ we have [X$_{\lambda}$, Y$_{\lambda}$] = q$_{\lambda} \tilde{H}_{\lambda}$ for some q$_{\lambda} \in \bq$.
\end{proof}

\begin{proposition}  Let $\lambda \in \Lambda$ and $0 \neq u \in$ U$_{\lambda}$.  Let $\fN_{0}(u,\lambda) = \br$-span \{u~,$\tilde{H}_{\lambda}$(u)~, $\tilde{H}_{\lambda}$ \}.  If $2\lambda \notin \Phi$ then $\fN_{0}(u,\lambda)$ is a 3-dimensional, totally geodesic subalgebra of $\fN_{0}$  with 1-dimensional center.  If u $\in \bq$-span($\fB_{\lambda}$) , then $\fN_{0}(u,\lambda)$ is a rational subalgebra of $\fN_{0}$ with respect to the Chevalley rational structure $\fN_{0, \bq}$ on $\fN_{0}$.
\end{proposition}

$\mathit{Remark}$  All 3-dimensional nilpotent Lie algebras with 1-dimensional center are isomorphic to the Heisenberg algebra.

\begin{proof}  Note that $\tilde{H}_{\lambda}$(u) $\in U_{\lambda}$ by 1) of (6.3).  From 1) of (6.10) it follows that
[$\fN_{0}(u,\lambda), \fN_{0}(u,\lambda)$] $\subseteq \br \tilde{H}_{\lambda}$, and hence $\fN_{0}(u,\lambda)$ is a subalgebra of $\fN_{0}$.  By  6) of (6.3) if $\lambda_{0} = i \lambda$, then $\lambda_{0} \in$ Hom($\fH_0,\br$) and $\tilde{H}_{\lambda}^{2} =
 -~ \lambda_{0}(\tilde{H}_{\lambda})^{2}$ Id on U$_{\lambda}$.  In particular $\tilde{H}_{\lambda}$ leaves invariant $\br$-span\{u, $\tilde{H}_{\lambda}$(u) \}, and it follows from (10.1) that $\fN_{0}(u,\lambda)$ is a totally geodesic subalgebra of $\fN_{0}$.  If u $\in \bq$-span($\fB_{\lambda}$) $=$ U$_{\lambda} \cap \bq$-span($\fB_{0}$), then $\tilde{H}_{\lambda}$(u) $\in \bq$-span($\fB_{\lambda}$) by (6.9). Hence $\fN_{0}(u,\lambda)$ is a rational subalgebra of $\fN$ since the basis \{u,$\tilde{H}_{\lambda}$(u),$\tilde{H}_{\lambda}$\} lies in $\fN_{0, \bq}=\bq$-span ($\fB_{0}~ \cup~ \fC_0$).
 \end{proof}

 $\mathit{Rational~ totally~ geodesic~ subalgebras~ of~ quaternionic~ type}$

\begin{proposition}  Let $\lambda \in \Lambda$ and $\beta \in \Phi^{+}$ be linearly independent.  Define U$'_{\lambda,\beta} =$ \(\sum_{k=-\infty}^{\infty}U_{\lambda+k\beta} \).  Let $\fN_{0}(\lambda,\beta) = U'_{\lambda,\beta} \oplus \fG_{\lambda,\beta}$, where $\fG_{\lambda,\beta} = \fG_{0,\beta} \oplus \br \tilde{H}_{\lambda} \oplus \br \tilde{\tau}_{\beta}$.  Then

\indent 1)  $\fG_{\lambda,\beta}$ is a Lie subalgebra of $\fG_0$ that is Lie algebra isomorphic to the quaternions $\bh$ with the Lie algebra structure given by [x,y] = xy $-$ yx for all x,y $\in \bh$.  Moreover, $\fG_{\lambda,\beta}$ is rational with respect to the rational structure $\bq$-span($\fC_{0}$) for $\fG_{0}$.

\indent 2)  U$'_{\lambda,\beta}$ is a $\fG_{\lambda,\beta}$-module.

\indent 3)  If 2$\lambda$ + p$\beta \notin \Phi$ for any integer p, then $\fN_{0}(\lambda,\beta)$ is a totally geodesic subalgebra of $\fN_{0}$ with 4-dimensional center $\fG_{\lambda,\beta}$ that is rational with respect to the Chevalley rational structure $\fN_{0, \bq}$ on $\fN$.
\end{proposition}

$\mathit{Remark}$  We note that if  2$\lambda$ + p$\beta \notin \Phi$ for any integer p, then $\lambda$ and $\beta$ are linearly independent.

\begin{proof}  1).  Let $\fG_{\lambda,\beta} = \fG_{0,\beta} \oplus \br \tilde{H}_{\lambda} \oplus \br \tilde{\tau}_{\beta}$, regarded now as a 4-dimensional subspace of $\fG_0$.  Recall from section 3 that $\fC_{0} =$ \{$\tilde{\tau}_{\alpha} : \alpha \in \Delta, A_{\beta}, B_{\beta} : \beta \in \Phi^{+}$\}, where A$_{\beta} =~$X$_{\beta}~ -$ X$_{-\beta}$, B$_{\beta} =$ iX$_{\beta}~ +$ iX$_{-\beta}$ and X$_{\beta}$, X$_{-\beta}$ are elements of $\fG_{\beta},\fG_{-\beta}$ such that [X$_{\beta}$, X$_{-\beta}$] $= \tau_{\beta}$.  Moreover, \{A$_{\beta}$, B$_{\beta}$\} is a basis of $\fG_{0,\beta}$ by the discussion at the beginning of section 7 and following (7.1).  We recall from (6.8) that $\tilde{H}_{\lambda} = i H_{\lambda}$, where $H_{\lambda}$ is the complex weight vector determined by $\lambda$ and the Killing form B on $\fG$.

    In terms of the bracket operation [ , ] in $\fG = \fG_0^{\bc}$ it is not difficult to obtain the following bracket relations from 4) of (4.6), (7.1) and the remarks following (7.1) :

\indent (*) \hspace{.01in}  [A$_{\beta}$, B$_{\beta}$] = 2$\tilde{\tau}_{\beta}$ ; \hspace{.565in} [$\tilde{\tau}_{\beta}$, A$_{\beta}$]
= 2B$_{\beta}$; \hspace{.76in} [$\tilde{\tau}_{\beta}$, B$_{\beta}$] = $-$ 2A$_{\beta}$

\hspace{.2in}       [$\tilde{H}_{\lambda}$, B$_{\beta}$] = i$\beta$($\tilde{H}_{\lambda}$) A$_{\beta}$; \hspace{.2in}
[$\tilde{H}_{\lambda}$, A$_{\beta}$] = $-$ i$\beta$($\tilde{H}_{\lambda}$) B$_{\beta}$ ; \hspace{.2in}[$\tilde{\tau}_{\beta};
\tilde{H}_{\lambda}$] = 0

    Note that i $\beta(\tilde{H}_{\lambda}$) $\in \bq$ since $\beta(\tilde{\tau}_{\alpha}) =$ i $\beta(\tau_{\alpha}) \in$ i $\bz$ and
$\tilde{H}_{\lambda} \in \bq$-span $\{\tilde{\tau}_{\alpha}: \alpha \in \Delta \}$ by (6.9).  Hence $\fG_{\lambda,\beta}$ is a 4-dimensional subalgebra of $\fG_0$ that is rational with respect to the rational structure $\bq$-span($\fC_{0}$) on $\fG_{0}$.

    We now show that $\fG_{\lambda,\beta}$ is Lie algebra isomorphic to the quaternions $\bh$ with the Lie algebra structure given by [x,y] = xy $-$ yx for all x,y $\in \bh$.  Let c $= - \frac{1}{2} \beta(H_{\lambda}) \in \bq$.  Then  $\xi =$ H$_{\lambda} + c \tau_{\beta} \in \fH_{\br} \subset \fH_0^{\bc}$ and $\beta(\xi) = 0$.  It follows from the choice of $\xi$ that ad $\xi$ annihilates X$_{\beta}$, X$_{-\beta}$ and the root vector H$_{\beta} \in \fG$.  If $\tilde{\xi}=$ i $\xi$, then by (6.1) and (6.8) $\tilde{\xi} = \tilde{H}_{\lambda} + c \tilde{\tau}_{\beta} \in \fG_{\beta,\lambda} \cap \fH_{0}$.  Moreover, ad $\tilde{\xi}$ annihilates A$_{\beta}$, B$_{\beta}$, $\tilde{H}_{\lambda}$ and $\tilde{\tau}_{\beta}$; that is, $\tilde{\xi}$ lies in the center of $\fG_{\lambda,\beta}$.  Now let $\varphi  : \fG_{\lambda,\beta} \rightarrow \bh$ be the linear isomorphism such that $\varphi(\tilde{\xi}) = 1,  \varphi$(A$_{\beta}$) = i, $\varphi$(B$_{\beta}$) = j , and  $\varphi(\tilde{\tau}_{\beta}) = k$.  From the bracket relations (*) above it follows that $\varphi$ is a Lie algebra isomorphism.

    2)   By (6.6) $\fG_{0,\beta}$(U$_{\lambda+k\beta}) \subset$ U$_{\lambda+(k+1)\beta}~ \oplus$~ U$_{\lambda+(k-1)\beta}$) for every integer k, and it follows that $\fG_{0,\beta}$(U$'_{\lambda,\beta}$) $\subset$ U$'_{\lambda,\beta}$.  By 1) of (6.3) $\tilde{H} (U_{\lambda+k\beta}) \subseteq U_{\lambda+k\beta}$ for every integer k and every  $\tilde{H} \in \fH_0$.  Hence $\tilde{H}_{\lambda}$ and $\tilde{\tau}_{\beta}$ leave U$'_{\lambda,\beta}$ invariant, and we conclude that $\fG_{\lambda,\beta}$(U$'_{\lambda,\beta}$) $\subseteq$ U$'_{\lambda,\beta}$.

    3)  If k and j are any integers, then by the hypothesis and (6.7) we obtain
 [U$_{\lambda+k\beta}$, U$_{\lambda+j\beta}$] $\subseteq \fG_{0,2\lambda+(k+j)\beta}~ \oplus~ \fG_{0,(k-j)\beta} = \fG_{0,(k-j)\beta}$. Note that $\fG_{0,(k-j)\beta} = \fG_{0,\beta}$ if $|k-j| = 1$.  If k = j and $\lambda + k\beta \in \Lambda$, then  [U$_{\lambda+k\beta}$, U$_{\lambda+j\beta}$] $ = \br \tilde{H}_{\lambda + k\beta}$ by (6.10).  In all other cases $\fG_{0,(k-j)\beta} =\{0 \}$.  By (6.8) we note that $\tilde{H}_{\beta}=$ i H$_{\beta}$ is a real multiple of $\tilde{\tau}_{\beta}= i~\tau_{\beta}$.  Moreover,  $\tilde{H}_{\lambda+k\beta}= \tilde{H}_{\lambda} + k \tilde{H}_{\beta}$ for all integers k.  Hence [U$'_{\lambda,\beta}$, U$'_{\lambda,\beta}$] $\subset \fG_{\lambda,\beta}$ and $\fN_{0}(\lambda,\beta)$ is a subalgebra of $\fN_{0}$.  Proposition 10.1 and 2) now imply that $\fN_{0}(\lambda,\beta)$ is a totally geodesic subalgebra of $\fN_{0}$ since $\fN_{0}(\lambda,\beta) \cap \fG_0 = \fG_{\lambda,\beta}$ and $\fN_{0}(\lambda,\beta) \cap U = U'_{\lambda,\beta}$.

    The subspaces U$_{\lambda+k\beta}, k \in \bz$, have bases $\fB_{\lambda+k\beta} \subset \fN_{\bq}=\bq$-span($\fB \cup \fC_0$)
whenever $\lambda + k\beta \in \Lambda$. Hence U$'_{\lambda,\beta}$ is a rational subspace of $\fN_{0}$, and by 1) $\fN_{0}(\lambda,\beta)$ is a rational subalgebra of $\fN_{0}$.
\end{proof}

\section{Appendix}

    We prove Propositions 9.1 and 9.2.  The proofs of both results are similar.

    It will be convenient to list  the highest long and short roots $\beta_{max}$ and $\beta_{min}$ for each complex simple Lie algebra, both as $\bz -$ linear combinations of the simple roots $\{\alpha_{i}\}$ and as $\bz -$ linear combinations of the fundamental weights $\{\omega_{i}\}$.  These lists can be determined from standard sources such as sections 12.2 and 13.2 of [Hu].
\newline

$\mathbf{\beta_{max}~ and~ \beta_{min}~ for~ the~ complex~ simple~ Lie~ algebras}$

A$_{n}, n \geq 2  \hspace{.5in} \beta _{max} = \beta_{min} = \alpha_{1} +~ ...~ \alpha_{n} = \omega_{1} + \omega_{n}$.

B$_{n}, n \geq 3 \hspace{.5in} \beta_{max} = \alpha_{1} + 2\alpha_{2} +~ ...~ +2\alpha_{n} = \omega_{2}$.

    $\hspace{1.05in} \beta_{min} = \alpha_{1} +~ ...~ + \alpha_{n} = \omega_{1}$.

C$_{n}, n \geq 2 \hspace{.5in} \beta_{max} = 2\alpha_{1} + 2\alpha_{2} +~ ...~ +2\alpha_{n-1} +\alpha_{n} = 2\omega_{1}$.

    $\hspace{1.05in} \beta_{min} = \alpha_{1} + 2\alpha_{2} +~ ...~ + 2\alpha_{n-1} + \alpha_{n}  = \omega_{2}$.

D$_{n}, n \geq 4 \hspace {.5in} \beta_{max} = \beta_{min} = \alpha_{1} + 2\alpha_{2} +~ ...~ + 2\alpha_{n-2} + \alpha_{n-1} + \alpha_{n} = \omega_{2}$.

E$_{6} \hspace{.9in} \beta_{max} = \beta_{min} = \alpha_{1} + 2\alpha_{2} + 2\alpha_{3} + 3\alpha_{4} + 2\alpha_{5} + \alpha_{6} = \omega_{2}$.

E$_{7} \hspace{.9in} \beta_{max} = \beta_{min} = 2\alpha_{1} + 2\alpha_{2} +  3\alpha_{3} + 4\alpha_{4} + 3\alpha_{5} + 2\alpha_{6} + \alpha_{7} = \omega_{1}$.

E$_{8} \hspace{.9in} \beta_{max} = \beta_{min} = 2\alpha_{1} + 3\alpha_{2} +  4\alpha_{3} + 6\alpha_{4} + 5\alpha_{5} + 4\alpha_{6} + 3\alpha_{7} + 2\alpha_{8} = \omega_{8}$.

F$_{4} \hspace{.9in} \beta_{max} = 2\alpha_{1} + 3\alpha_{2} + 4\alpha_{3} + 2\alpha_{4} = \omega_{1}$.

    $\hspace{1.05in} \beta_{min} = \alpha_{1} + 2\alpha_{2} + 3\alpha_{3} + 2\alpha_{4} = \omega_{4}$

G$_{2} \hspace{.9in} \beta_{max} = 3\alpha_{1} + 2\alpha_{2} = \omega_{2}$.

    $\hspace{1.05in} \beta_{min} = 2\alpha_{1} + \alpha_{2} = \omega_{1}$.

$\mathit{Proof~of~9.1}$

    Our problem is to find all solutions, up to the action of the Weyl group W, to the equation 2$\mu = (2k+1)\alpha + \beta$ for an abstract weight $\mu$ and elements $\alpha, \beta \in \Phi$.  We recall that $\fG$ is assumed to be a simple complex Lie algebra.
    The proof of (9.1) follows the statement of Lemmas (11.4) and (11.5), which are the central parts of the proof.  We begin with some preparatory results.

\begin{lemma}  If we replace $\mu$ by w($\mu$) for a suitable element w of W, then there exist elements $\alpha' , \beta' \in \Phi$ and a dominant weight $\mu' \in \tilde{\Lambda}$ such that $\mu = \mu' + k\alpha'$ and $2\mu' = \alpha' + \beta'$.
\end{lemma}

\begin{proof} Let $\mu' = \mu - k\alpha$.  Clearly $2\mu' = \alpha + \beta$ and $\mu = \mu' + k\alpha$.  Now choose an element w of W such that w($\mu'$) is a dominant weight.
\end{proof}

    From the result above we see that $\mu' \in \tilde{\Lambda}^{+}$, the set of dominant weights in $\tilde{\Lambda}, 2\mu'  \in \bz-span(\Delta)$ and $2\mu' \leq 2\beta_{max}$, where the final statement means that $2\beta_{max} - 2\mu'$ is a sum of elements from $\Delta$.

    The next result is the major tool used in the proof of (11.4).

\begin{lemma}   Let $\mu =$ \(\sum_{i=1}^{n} q_{i}\omega_{i}\)  $\in \tilde{\Lambda}^{+}$ be an element such that $2\mu \in \bz-span(\Delta)$ and $2\mu \leq 2\beta_{max}$.  Write $\mu =$ \(\sum_{i=1}^{n} q_{i}\omega_{i}\), where the integers q$_{i}$ are nonnegative.  Write $\beta_{max} =$ \(\sum_{j=1}^{n} \beta_{j}\alpha_{j}\), for suitable nonnegative integers $\beta_{j} = (\beta_{max})_{j}$. Let C$^{ij}$ denote the inverse of the Cartan matrix C$_{ij} = \langle \alpha_{i}~,~\alpha_{j} \rangle$. Then

\indent  \indent  1)  \(\sum_{i=1}^{n} q_{i} C^{ij}\) $ \leq (\beta_{max})_{j}$  for 1 $\leq$ j $\leq$ n.

\indent \indent 2)  $q_{i} C^{ij} \leq (\beta_{max})_{j}$ for all i,j.  If equality holds for some i,j, then $\mu = q_{i}\omega_{i}$.
\end{lemma}

\begin{proof}

\indent 1)  By definition of the Cartan matrix $\alpha_{i} =$
$\sum_{j=1}^{n} C_{ij}\omega_{j}$, and hence $\omega_{i} =$
$\sum_{j=1}^{n} C^{ij}\alpha_{j}$.  It follows that $\mu =$
$\sum_{i=1}^{n} q_{i}\omega_{i} = \sum_{j=1}^{n} (\sum_{i=1}^{n}
q_{i} C^{ij}) \alpha_{j}$. Since $2\mu \leq 2\beta_{max}$ we obtain
2$(\sum_{i=1}^{n} q_{i} C^{ij})$ $\leq 2(\beta_{max})_{j}$ for 1
$\leq$ j $\leq$ n.  This proves 1).

\indent 2)  A case by case inspection of the inverse Cartan matrices (C$^{ij}$) shows that C$^{ij} > 0$ for all i,j.  Hence  $q_{i} C^{ij} \leq$ \(\sum_{k=1}^{n} q_{k} C^{kj}\) $ \leq (\beta_{max})_{j}$ by 1).  If $q_{i} C^{ij} = (\beta_{max})_{j}$ for some i,j, then q$_{k} = 0$ for k $\neq$ i since C$^{kj} > 0$, which proves 2).
\end{proof}

    The following special case of (9.1) will be needed frequently in the proof of (9.1).

\begin{lemma} If $2\beta_{max} = \alpha + \beta$ for $\alpha, \beta \in \Phi$, then $\alpha = \beta = \beta_{max}$.
\end{lemma}

\begin{proof}  We choose integers $a_{i}, b_{i}$ and $(\beta_{max})_{i}$ such that $\beta_{max} =$ \(\sum_{i=1}^{n}(\beta_{max})_{i} \alpha_{i}\), $\alpha =$ \(\sum_{i=1}^{n} a_{i} \alpha_{i}\), and $\beta =$ \(\sum_{i=1}^{n} b_{i} \alpha_{i}\).  Note that $a_{i} \leq (\beta_{max})_{i}$ and $b_{i} \leq (\beta_{max})_{i}$ for all i since $\alpha \leq \beta_{max}$ and $\beta \leq \beta_{max}$.  If \(\sum_{i=1}^{n}2(\beta_{max})_{i} \alpha_{i}\) $= 2\beta_{max} = \alpha + \beta =$ \(\sum_{i=1}^{n} (a_{i} + b_{i}) \alpha_{i}\), then it follows that $a_{i} = b_{i} = (\beta_{max})_{i}$ for all i.  Hence $\alpha = \beta = \beta_{max}$.
\end{proof}

    The proof of (9.1) will follow from the next two results :

\begin{lemma}  Let $\mu \in \tilde{\Lambda}^{+}$ be an element such that $2\mu \in \bz-span(\Delta)$ and $2\mu \leq 2\beta_{max}$.  Then $\mu$ belongs to the following list :

A$_{n}, n \geq 2$ \hspace{.1in} $\mu = (\omega_{1} + \omega_{n})  = \beta_{max} = \beta_{min}$.

A$_{3}$ \hspace {.5in}   $\mu = \omega_{2}$.

B$_{n},n \geq 3$

    \hspace{.7in}     $\mu = \omega_{2} = \beta_{max}$.

    \hspace{.7in}   $\mu = \omega_{1} = \beta_{min}$.

B$_{4}$ \hspace{.5in}   $\mu = \omega_{4}$.

B$_{3}$ \hspace{.5in}   $\mu = \omega_{3}$.

C$_{n}, n \geq 2$

    \hspace {.7in}  $\mu = \omega_{1} = \frac{1}{2} \beta_{max}$.

    \hspace{.7in} $\mu = 2\omega_{1} = \beta_{max}$.

    \hspace{.7in} $\mu = \omega_{2}  = \beta_{min}$.

D$_{n}, n \geq 4$

    \hspace {.7in}  $\mu = \omega_{1}$.

    \hspace {.7in}  $\mu = \omega_{2} = \beta_{max} = \beta_{min}$.

D$_{4}$ \hspace{.5in}   $\mu = \omega_{3}$.

    \hspace{.7in}   $\mu = \omega_{4}$.

E$_{6}$ \hspace{.5in}   $\mu = \omega_{2} = \beta_{max} = \beta_{min}$.

E$_{7}$ \hspace{.5in}   $\mu = \omega_{1} = \beta_{max} = \beta_{min}$.

E$_{8}$ \hspace{.5in}   $\mu = \omega_{8} = \beta_{max} = \beta_{min}$.

F$_{4}$ \hspace{.5in}       $\mu = \omega_{1} = \beta_{max}$.

    \hspace{.7in}   $\mu = \omega_{4} = \beta_{min}$.

G$_{2}$ \hspace{.5in}   $\mu = \omega_{2} = \beta_{max}$.

    \hspace{.7in} $\mu = \omega_{1} = \beta_{min}$.

\end{lemma}

    If $\mu = \omega_{i}$ is a fundamental weight such that (*) $2\omega_{i} = \alpha + \beta$ for $\alpha, \beta \in \Phi$, then clearly $2\omega_{i} \in \bz-span(\Delta)$ and $2\mu \leq 2\beta_{max}$, the hypotheses of (11.4).  Note that the elements $\alpha, \beta \in \Phi$ that satisfy (*) are invariant under the action of W$_{i} = \{w \in W : w(\omega_{i}) = \omega_{i}\}$.   However, the next result shows that there are at most two such W$_{i}$ orbits for any of the complex simple Lie algebras.

\begin{lemma}  Let $\omega_{i}, 1 \leq i \leq n$ be a fundamental dominant weight such that $2\omega_{i} = \alpha + \beta$ for $\alpha, \beta \in \Phi$.  Then there exist elements w$_{1}$, w$_{2} \in$ W such that

\indent 1)  w$_{1}(\omega_{i}) =$ w$_{2}(\omega_{i}) = \omega_{i}$.

\indent 2)  w$_{1}(\alpha) = \beta_{max}$ or $\beta_{min}$ and w$_{1}(\beta) \in \Phi^{+}$.

\indent 3)  w$_{2}(\beta) = \beta_{max}$ or $\beta_{min}$ and w$_{2}(\alpha) \in \Phi^{+}$.
\newline

    Moreover, if $2\omega_{i} = \alpha + \beta$, where $\alpha$ or $\beta$ is $\beta_{max}$ or $\beta_{min}$, then the following is a complete list of solutions.

\indent A$_{3}  \hspace{.5in} \mu = \omega_{2}, 2\mu = \beta_{max} + \alpha_{2} = \beta_{min} + \alpha_{2}$.

\indent B$_{n}, n \geq 3$

\indent  $\hspace{.7in} \mu = \omega_{2} = \beta_{max}, 2\mu = \beta_{max} + \beta_{max}$.

\indent $\hspace{.7in} \mu = \omega_{1} = \beta_{min}, 2\mu = \beta_{min} + \beta_{min} = \beta_{max} + \alpha_{1}$.

\indent B$_{4}  \hspace{.5in} \mu = \omega_{4}, 2\mu = \beta_{max} + \alpha_{3} + 2\alpha_{4} = \beta_{max} + \sigma_{\alpha_{4}}(\alpha_{3})$.

\indent B$_{3}  \hspace{.5in} \mu = \omega_{3}, 2\mu = \beta_{max} + \alpha_{3}  = \beta_{min} + \alpha_{2} + 2\alpha_{3} = \beta_{min} + \sigma_{\alpha_{3}}(\alpha_{2})$.

\indent C$_{n}, n \geq 2$

\indent $\hspace{.7in} \mu = \omega_{1} = \frac{1}{2} \beta_{max}, 2\mu = \beta_{min} + \alpha_{1}$.

\indent $\hspace{.7in} \mu = 2\omega_{1} = \beta_{max}, 2\mu = \beta_{max} + \beta_{max}$.

\indent $\hspace{.7in} \mu = \omega_{2} = \beta_{min}, 2\mu = \beta_{min} + \beta_{min} = \beta_{max} + 2\alpha_{2} + ~...~ + 2\alpha_{n-1} + \alpha_{n} = $

    $\hspace{.7in} \beta_{max} + \sigma_{\alpha_{1}}(\beta_{max})$.

\indent D$_{n}, n \geq 4$

\indent $\hspace{.7in} \mu = \omega_{1}, 2\mu = \beta_{max} + \alpha_{1} = \beta_{min} + \alpha_{1}$.

\indent $\hspace{.7in} \mu = \omega_{2} = \beta_{max} = \beta_{min}, 2\mu = \beta_{max} + \beta_{max}$.

\indent D$_{4}  \hspace{.5in} \mu = \omega_{3}, 2\mu = \beta_{max} + \alpha_{3}$

\indent $\hspace{.7in} \mu = \omega_{4}, 2\mu = \beta_{max} + \alpha_{4}$.

\indent E$_{6}  \hspace{.5in} \mu = \omega_{2} = \beta_{max} = \beta_{min}, 2\mu = \beta_{max} + \beta_{max}$.

\indent E$_{7}  \hspace{.5in} \mu = \omega_{1} = \beta_{max} = \beta_{min}, 2\mu = \beta_{max} + \beta_{max}$.

\indent E$_{8}  \hspace{.5in} \mu = \omega_{8} = \beta_{max} = \beta_{min}, 2\mu = \beta_{max} + \beta_{max}$.

\indent F$_{4}  \hspace{.5in} \mu = \omega_{1} = \beta_{max}, 2\mu = \beta_{max} + \beta_{max}$.

\indent $\hspace{.7in} \mu = \omega_{4} = \beta_{min}, 2\mu = \beta_{min} + \beta_{min} = \beta_{max} + \alpha_{2} + 2\alpha_{3} + 2\alpha_{4} =$

    $\hspace{.7in}  \beta_{max} + \sigma_{\alpha_{4}}  \sigma_{\alpha_{3}} (\alpha_{2})$.

\indent G$_{2}  \hspace{.5in} \mu = \omega_{2} = \beta_{max}, 2\mu = \beta_{max} + \beta_{max}$.

\indent $\hspace{.7in} \mu = \omega_{1} = \beta_{min}, 2\mu = \beta_{min} + \beta_{min} = \beta_{max} + \alpha_{1}$.

\end{lemma}

$\mathit{Proof~of~(9.1)}$  Before proving (11.4) and (11.5) we complete the proof of (9.1).  We consider abstract weights $\mu = \mu' + k\alpha$ as in (11.1), and we note that $\mu'$ lies in the list of (11.4).  If $\mu' = \beta_{max}$ for A$_{n}, n \geq 2$ or C$_{n}, n \geq 2$, then the assertions of (9.1) for these cases follow from (11.3).  In all other cases of (11.4) we note that $\mu'$ is a fundamental weight $\omega_{i}$ for $1 \leq i \leq n$.  In these cases we use (11.5) to replace $\alpha'$ in (11.1) by $\alpha' = \beta_{max}$ or $\alpha' = \beta_{min}$.  The assertions of (9.1) now follow.

    We now prove (11.4) and (11.5).
\newline

$\mathit{Proof~of~Lemma~11.4}$  Let C$^{ij}$ denote the inverse of the Cartan matrix C$_{ij} = \langle \alpha_{i}~,~\alpha_{j} \rangle$.  Recall that $\omega_{i} =$ \(\sum_{j=1}^{n} C^{ij} \alpha_{j} \), or equivalently, that $C^{ij}$ is the $\alpha_{j}$ coefficient of $\omega_{i}$.  For a table of the inverse Cartan matrices see, for example, section 13.2 of [Hu].

    We use (11.2) to prove (11.4) in each case of the classification of complex simple Lie algebras.

A$_{n}, n \geq 4$.

    We omit the proof in the cases n = 2,3, which is similar to the proof for n $\geq 4$.  We recall

\indent 1) $\omega_{i} = (\frac{n-i+1}{n+1})\{\alpha_{1} + 2\alpha_{2} + ... + i\alpha_{i}\} +
(\frac{i}{n+1})\{(n-i)\alpha_{i+1} + (n-i-1)\alpha_{i+2} + ... + \alpha_{n}\}$.

\indent 2)  $\beta_{max} = \beta_{min} = \alpha_{1} +~ ...~ \alpha_{n} = \omega_{1} + \omega_{n}$.  In particular, $(\beta_{max})_{j} = 1$ for all j.
\newline

$\mathbf{Sublemma}$  Let $\mu =$ \(\sum_{i=1}^{n} q_{i}\omega_{i}\)  $\in \tilde{\Lambda}^{+}$ be an element such that $2\mu \in \bz-span(\Delta)$ and $2\mu \leq 2\beta_{max}$.  Suppose that C$^{ii} \leq 1$ for some integer i with $2 \leq$ i $\leq n-1$.  Then n = 3 and i = 2.  In this case C$^{22} = 1$.
\newline

$\mathit{Proof~of~the~Sublemma}$  By hypothesis $1 \geq C^{ii} = \frac{i(n-i+1)}{n+1}$, which implies that $n+1 \geq in - i^2 + i$, or equivalently, $1 \geq (i-1)n - i(i-1) = (i-1)(n-i)$.  Since $i-1$ and $n-i$ are positive integers we conclude that $i-1 = n - i = 1$, which completes the proof. $\hspace{.2in} \square$
\newline

$\mathbf{Corollary}$  Let n $\geq 4$, and let $\mu =$ \(\sum_{i=1}^{n} q_{i}\omega_{i}\), where q$_{i} \in \bz^{+}, 2\mu \in \bz-span(\Delta)$ and $2\mu \leq 2\beta_{max}$.  Then $\mu = \omega_{1} + \omega_{n} = \beta_{max} = \beta_{min}$.
\newline

$\mathit{Proof~of~the~Corollary}$  By the lemma above we have \(\sum_{i=1}^{n} q_{i} C^{ij}\) $ \leq (\beta_{max})_{j} = 1$ for all j.  If q$_{i} \neq 0$ for some i with $2 \leq i \leq n-1$, then $1 \geq q_{i}C^{ii} \geq C^{ii} > 2$ by the sublemma.  Since this is impossible we conclude that $\mu = q_{1}\omega_{1} + q_{n}\omega_{n}$ for some integers $q_{1}, q_{n} \in \bz^{+}$.

    Next, suppose that both $q_{1}$ and $q_{n}$ are positive integers.  From the lemma above we obtain $1 = (\beta_{max})_{j} \geq$  \(\sum_{i=1}^{n} q_{i} C^{ij}\) $= q_{1}C^{1j} + q_{n}C^{nj} \geq C^{1j} + C^{nj} = \frac{n-j+1}{n+1} + \frac{j}{n+1} = 1$.  It follows that $q_{1} = q_{n} = 1$ and $\mu = \omega_{1} + \omega_{n}$.

    It remains only to rule out the cases $\mu = q_{1}\omega_{1}$ and $\mu = q_{n}\omega_{n}$ for integers $q_{1} \geq 1$ and  $q_{n} \geq 1$.  Suppose first that $\mu = q_{1}\omega_{1}$ for some integer $q_{1} \geq 1$.  Since $1 \geq q_{1}C^{11} = \frac{q_{1}n}{n+1}$ we must have $q_{1} = 1$, which implies that $\mu = \omega_{1}$.  The $\alpha_{2}$ coefficient of $2\mu = 2\omega_{1}$ is $\frac{2n-2}{n+1}$.  However, $\frac{2n-2}{n+1}$ is not an integer for n $\geq 4$, which contradicts the hypothesis that $2\mu \in \bz-span(\Delta)$.  This rules out the possibility $\mu = \omega_{1}$.

    Similarly, if  $\mu = q_{n}\omega_{n}$ for some integer $q_{n} \geq 1$, then $1 \geq q_{n}C^{nn} = q_{n} \frac{n}{n+1}$.  This implies that $q_{n} = 1$ and $\mu = \omega_{n}$.  However, $2\mu = 2\omega_{n}$ has $\alpha_{n}$ coefficient $\frac{2n}{n+1}$, which is not an integer for n $\geq 4$.  This contradicts the hypothesis that $2\mu \in \bz-span(\Delta)$ and rules out the possibility $\mu = \omega_{n}$.  This completes the proof of the Corollary. $\hspace{.5in} \square$
\newline

B$_{n}, n \geq 3$

    We recall that $\beta_{max} = \alpha_{1} + 2\alpha_{2} + ~...~ + 2\alpha_{n} = \omega_{2}$ and $\beta_{min} = \alpha_{1} + \alpha_{2} + ~...~ + \alpha_{n} = \omega_{1}$.  In particular we always have the solutions $\mu = \omega_{1} = \beta_{min}$ and  $\mu = \omega_{2} = \beta_{max}$.
\newline

$\mathbf{Sublemma~1}$  Let $\mu =$ \(\sum_{i=1}^{n} q_{i}\omega_{i}\)  $\in \tilde{\Lambda}^{+}$ be an element such that $2\mu \in \bz-span(\Delta)$ and $2\mu \leq 2\beta_{max}$.  Then either $\mu = \omega_{i}$ for 1 $\leq$ i $\leq$ n or $\mu = 2\omega_{n}$.

$\mathit{Proof}$   We note that C$^{i1} = 1$ for 1 $\leq$ i $\leq$ n $-1$ and C$^{n1} = \frac{1}{2}$.  Hence $1 = (\beta_{max})_{1} \geq$ \(\sum_{i=1}^{n} q_{i}C^{i1}\) $= (q_{1} + ~...~ + q_{n-1} + \frac{1}{2}q_{n})$.  Exactly one q$_{i}$ is nonzero, and the sublemma follows immediately.
\newline

$\mathbf{Sublemma~2}$  Let $\mu$ be as above.  If n $\geq 5$, then either $\mu = \omega_{1}$ or $\mu = \omega_{2}$.

$\mathit{Proof}$  If $3 \leq i \leq n-1$, then C$^{ii} = i \geq 3$.  Hence $2 = (\beta_{max})_{i} \geq q_{i}C^{ii} \geq 3q_{i}$, which implies that $q_{i} = 0$.  By Sublemma 1 only $\mu = \omega_{1},  \mu = \omega_{2},  \mu = \omega_{n}$ and $\mu = 2 \omega_{n}$ are possible.  However, if $\mu = \omega_{n}$, then $2\mu = 2\omega_{n} = \alpha_{1} + 2\alpha_{2} +~ ...~ + n\alpha_{n}$.  Since $(\beta_{max})_{n} = 2$ and n $\geq 5$ we obtain a contradiction to the condition $2\mu \leq 2\beta_{max}$.  Similarly, if $\mu = 2\omega_{n}$, then we obtain a contradiction to the condition $2\mu \leq 2\beta_{max}$.
\newline

$\mathbf{Sublemma~3}$ Let $\mu$ be as above.  If n $= 4$, then $\mu = \omega_{1}, \mu = \omega_{2}$ or $\mu = \omega_{4}$.

$\mathit{Proof}$  By sublemma 1 we need only consider $\mu = \omega_{3} = \alpha_{1} + 2\alpha_{2} + 3\alpha_{3} + 3\alpha_{4}, \omega_{4} = \frac{1}{2} (\alpha_{1} + 2\alpha_{2} + 3\alpha_{3} + 4\alpha_{4})$ and $\mu = 2\omega_{4}$.  If $\mu = \omega_{3}$, then the condition $2\mu \leq 2\beta_{max}$ fails.  If $\mu = \omega_{4}$, then the condition $2\mu = 2\omega_{4} \leq 2\beta_{max}$ holds, but if $\mu = 2\omega_{4}$, then the condition $2\mu = 4\omega_{4} \leq 2\beta_{max}$ fails.
\newline

$\mathbf{Sublemma~4}$   Let $\mu$ be as above.  If n $= 3$, then $\mu = \omega_{1}, \mu = \omega_{2}$ or $\mu = \omega_{3}$.

$\mathit{Proof}$  By sublemma 1 we need only consider $\omega_{3} = \frac{1}{2} (\alpha_{1} + 2\alpha_{2} + 3\alpha_{3})$ and $\mu = 2\omega_{3}$.  If $\mu = \omega_{3}$, then the condition $2\mu \leq 2\beta_{max}$ holds, but if $\mu = 2\omega_{3}$, then the condition $2\mu \leq 2\beta_{max}$ fails.
\newline

C$_{n}, n \geq 2$

    We recall that $\beta_{max} = 2\alpha_{1} + ~...~ + 2\alpha_{n-1} + \alpha_{n} = 2\omega_{1}$ and $\beta_{min} = \alpha_{1} + 2\alpha_{2} ~...~ + 2\alpha_{n-1} + \alpha_{n} = \omega_{2}$.  We must show that if $\mu$ satisfies the hypotheses of (11.4), then $\mu = \omega_{1} = \frac{1}{2}\beta_{max}, \mu = 2\omega_{1} = \beta_{max}$ or $\mu = \omega_{2} = \beta_{min}$.

    We note that C$^{in} = \frac{1}{2}i$.  By (11.2) we have $1 = (\beta_{max})_{n} \geq q_{i}C^{in} = \frac{1}{2}iq_{i}$ for $1 \leq$ i $\leq$ n.  Hence if q$_{i} \neq 0$, then i = 1 or 2.  We conclude that $\mu = q_{1}\omega_{1} + q_{2}\omega_{2}$ for $q_{1},~q_{2} \in \bz^{+}$.  Moreover, from (11.2) we also obtain $\frac{1}{2}q_{1} + q_{2} = q_{1}C^{1n} + q_{2}C^{2n} \leq 1$.  If $q_{2} \neq 0$, then $q_{2} = 1$ and $q_{1} = 0$; that is, $\mu = \omega_{2}$.   If $q_{1} \neq 0$, then $q_{2} = 0$ and $q_{1} = 1$ or $2$; that is, $\mu = \omega_{1} = \frac{1}{2}\beta_{max}$ or $\mu = 2\omega_{1} = \beta_{max}$.  In both cases the condition $2\mu \leq 2\beta_{max}$ is satisfied.
\newline

D$_{n}, n \geq 4$

    We recall that $\beta_{max} = \beta_{min} = \alpha_{1} + 2\alpha_{2} + ~...~ + 2\alpha_{n-2} + \alpha_{n-1} + \alpha_{n} = \omega_{2}$.
\newline

$\mathbf{Sublemma~1}$  Let $\mu =$ \(\sum_{i=1}^{n} q_{i}\omega_{i}\)  $\in \tilde{\Lambda}^{+}$ be an element such that $2\mu \in \bz-span(\Delta)$ and $2\mu \leq 2\beta_{max}$.  Then $\mu$ satisfies one of the following :

\indent 1)  $\mu = \omega_{i}$ for some integer i with 1 $\leq$ i $\leq$ n.

\indent 2)  $\mu = \omega_{n-1} + \omega_{n}$.

\indent 3)  $\mu = 2\omega_{n-1}$.

\indent 4)  $\mu = 2\omega_{n}$.

$\mathit{Proof}$ Note that C$^{i1} = 1$ if 1 $\leq$ i $\leq$ n $- 2$ and C$^{(n-1)1} = $ C$^{n1} = \frac{1}{2}$.  By (11.2) we obtain $1 = (\beta_{max})_{1} \geq$  \(\sum_{i=1}^{n} q_{i} C^{i1}\) $ = (q_{1} + ~...~ q_{n-2}) + \frac{1}{2}(q_{n-1} + q_{n})$.  Hence if $q_{i} \neq 0$ for 1 $\leq$ i $\leq$ n $-2$, then $q_{k} = 0$ for k $\neq$ i and we conclude that $\mu = \omega_{i}$.  If $q_{i} \neq 0$ for i $=$ n $-1$ or n, then $q_{k} = 0$ for 1 $\leq$ k $\leq$ n $-2$ and $q_{n-1} + q_{n} \leq 2$.  This proves the sublemma.

    For n $\geq$ 4 it is routine to check that $\mu = \omega_{1}$ and $\mu = \omega_{2}$ are solutions.  We show next that the converse is true for n $\geq$ 5.
\newline

 $\mathbf{Sublemma~2}$  Let $\mu$ be as above.  If n $\geq$ 5, then the only solutions are $\mu = \omega_{1}$ and $\mu = \omega_{2}$.

 $\mathit{Proof}$  From (11.2) we obtain

 \indent    a)  $1 = (\beta_{max})_{n} \geq$ \(\sum_{i=1}^{n} q_{i} C^{in}\) $= \frac{1}{2} \{q_{1} + 2q_{2} + ~...~ (n-2)q_{n-2}\} + \frac{1}{4}(n-2)q_{n-1} + \frac{1}{4}nq_{n}$.

 \indent    b)  $1 = (\beta_{max})_{n-1} \geq$ \(\sum_{i=1}^{n} q_{i} C^{i(n-1)}\) $= \frac{1}{2} \{q_{1} + 2q_{2} + ~...~ (n-2)q_{n-2}\} + \frac{1}{4}nq_{n-1} + \frac{1}{4}(n-2)q_{n}$.

    Since n $\geq$ 5 it follows from a) and b) that $q_{n-1} = q_{n} = 0$.  Hence from  a) and b) we obtain $1 \geq \frac{1}{2} \{q_{1} + 2q_{2} + ~...~ (n-2)q_{n-2}\}$.  It follows that $q_{i} = 0$ for 3 $\leq$ i $\leq$ n$-2$, and we obtain $\mu = q_{1}\omega_{1} + q_{2}\omega_{2}$, where $1 \geq \frac{1}{2}(q_{1} + 2q_{2})$.

    If $q_{1} \neq 0$, then $q_{2} = 0$ and $q_{1} = 1$ or $2$.  The case $q_{1} = 2$ is ruled out by Sublemma 1, so only $\mu = \omega_{1}$ is possible if  $q_{1} \neq 0$.

    If $q_{2} \neq 0$, then $q_{2} = 1$ and $q_{1} = 0$, which shows that $\mu = \omega_{2}$.  This completes the proof of Sublemma 2.
\newline

    It remains only to consider the case n $= 4$.  If $q_{4} \neq 0$, then from a) of Sublemma 2 we obtain $1 \geq \frac{1}{2} \{q_{1} + 2q_{2}\} + \frac{1}{2}q_{3} + q_{4} \geq q_{4} \geq 1$.  We conclude that $q_{4} = 1$ and $q_{k} = 0$ for k $\neq$ 4; that is, $\mu = \omega_{4}$.  If $q_{3} \neq 0$, then from b) we obtain $1 \geq \frac{1}{2} \{q_{1} + 2q_{2}\} + q_{3} + \frac{1}{2}q_{4} \geq q_{3} \geq 1$.  We conclude that $q_{3} = 1$ and $q_{k} = 0$ for k $\neq$ 3; that is, $\mu = \omega_{3}$.  If $\mu = \omega_{3}$ or $\mu = \omega_{4}$, then the conditions $\mu \in \bz$ - span($\Delta$) and $2 \mu \leq 2 \beta_{max}$ are satisfied.
\newline

E$_{6}$

    We recall that $\beta_{max} = \beta_{min} = \alpha_{1} + 2\alpha_{2} + 2\alpha_{3} + 3\alpha_{4} + 2\alpha_{5} + \alpha_{6} = \omega_{2}$.  Let $\mu =$ \(\sum_{i=1}^{6} q_{i}\omega_{i}\)  $\in \tilde{\Lambda}^{+}$ be an element such that $2\mu \in \bz-span(\Delta)$ and $2\mu \leq 2\beta_{max}$.  By (11.2) we have $1 = (\beta_{max})_{1} \geq$ \(\sum_{i=1}^{6} q_{i} C^{i1}\).  Note that $C^{21} = 1, C^{61} = \frac{2}{3}$ and $C^{i1} > 1$ if i $\neq 2$ and i $\neq 6$.  Hence q$_{i} = 0$ if i $\neq$ 2, i $\neq$ 6, and we obtain $1 \geq q_{2}C^{21} + q_{6}C^{61} = q_{2} + \frac{2}{3}q_{6}$.  If q$_{2} \neq 0$, then q$_{2} = 1$ and q$_{6} = 0$; that is, $\mu = \omega_{2} = \beta_{max}$.  This is a solution.  If q$_{6} \neq 0$, then q$_{6} = 1$ and q$_{2} = 0$; that is, $\mu = \omega_{6}$.  However, $2\mu \notin \bz-span(\Delta)$ so $\mu = \omega_{6}$ is not a solution.

    We conclude that $\mu = \omega_{2} = \beta_{max} = \beta_{min}$ is the only solution.
\newline

E$_{7}$

    We recall that $\beta_{max} = \beta_{min} = 2\alpha_{1} + 2\alpha_{2} + 3\alpha_{3} + 4\alpha_{4} + 3\alpha_{5} + 2\alpha_{6} + \alpha_{7} = \omega_{1}$.   Let $\mu =$ \(\sum_{i=1}^{7} q_{i}\omega_{i}\)  $\in \tilde{\Lambda}^{+}$ be an element such that $2\mu \in \bz-span(\Delta)$ and $2\mu \leq 2\beta_{max}$.  From (11.2) we have $1 = (\beta_{max})_{7} \geq$   \(\sum_{i=1}^{7} q_{i} C^{i7}\).  We note that C$^{i7} > 1$ if i $\geq 2$.  Hence q$_{i} = 0$ for i $\geq 2$, which shows that $\mu = q_{1}\omega_{1}$, where $1 \geq q_{1}C^{17} = q_{1}$.  Hence q$_{1} = 1$ and $\mu = \omega_{1} = \beta_{max} = \beta_{min}$.

        We conclude that $\mu = \omega_{1} = \beta_{max} = \beta_{min}$ is the only solution.
\newline

E$_{8}$

    We recall that $\beta_{max} = \beta_{min} = 2\alpha_{1} + 3\alpha_{2} + 4\alpha_{3} + 6\alpha_{4} + 5\alpha_{5} + 4\alpha_{6} + 3\alpha_{7} + 2\alpha_{8} = \omega_{8}$.  Let $\mu =$ \(\sum_{i=1}^{8} q_{i}\omega_{i}\)  $\in \tilde{\Lambda}^{+}$ be an element such that $2\mu \in \bz-span(\Delta)$ and $2\mu \leq 2\beta_{max}$.  From (11.2) we have $2 = (\beta_{max})_{1} \geq$ \(\sum_{i=1}^{8} q_{i} C^{i1}\).  We note that C$^{i1} > 2$ if i $\neq 8$ and C$^{81} = 2$.  It follows that q$_{i} = 0$ if i $\neq 8$, and we obtain $2 \geq q_{8}C^{81} = 2q_{8}$.  Hence q$_{8} = 1$ and $\mu = \omega_{8}$.

    We conclude that $\mu = \omega_{8} = \beta_{max} = \beta_{min}$ is the only solution.
\newline

F$_{4}$

    We recall that $\beta_{max} = 2\alpha_{1} + 3\alpha_{2} + 4\alpha_{3} + 2\alpha_{4} = \omega_{1}$ and $\beta_{min} = \alpha_{1} + 2\alpha_{2} + 3\alpha_{3} + 2\alpha_{4} = \omega_{4}$.   Let $\mu =$ \(\sum_{i=1}^{4} q_{i}\omega_{i}\)  $\in \tilde{\Lambda}^{+}$ be an element such that $2\mu \in \bz-span(\Delta)$ and $2\mu \leq 2\beta_{max}$.  From (11.2) we have $2 = (\beta_{max})_{4} \geq$ \(\sum_{i=1}^{4} q_{i} C^{i4}\).  We note that $C^{14} = 2, C^{24} = 4, C^{34} = 3$ and $C^{44} = 2$.  hence $q_{2} = q_{3} = 0$, and we obtain $2 \geq q_{1}C^{14} + q_{4}C^{44} = 2(q_{1} + q_{4})$.  Hence either $q_{1} = 1, q_{4} = 0$ and $\mu = \omega_{1} = \beta_{max}$ or $ q_{4} = 1, q_{1} = 0$ and $\mu = \omega_{4} = \beta_{min}$.

    We conclude that $\mu = \omega_{1} = \beta_{max}$ and  $\mu = \omega_{4} = \beta_{min}$ are the only solutions.
\newline

G$_{2}$

    We recall that $\beta_{max} = 3\alpha_{1} + 2\alpha_{2} = \omega_{2}$ and $\beta_{min} = 2\alpha_{1} + \alpha_{2} = \omega_{1}$.   Let $\mu =$ \(\sum_{i=1}^{2} q_{i}\omega_{i}\)  $\in \tilde{\Lambda}^{+}$ be an element such that $2\mu \in \bz-span(\Delta)$ and $2\mu \leq 2\beta_{max}$.  From (11.2) we have $2 = (\beta_{max})_{2} \geq q_{1}C^{12} + q_{2}C^{22} = q_{1} + 2 q_{2}$.  If $q_{1} \neq 0$, then $q_{2} = 0$ and $q_{1} = 1$ or $2$; that is, $\mu = \omega_{1}$ or $2\omega_{1}$.  If $\mu = \omega_{1} = \beta_{min}$, then $2\mu = 4\alpha_{1} + 2\alpha_{2}$ satisfies the condition $2\mu \leq 2\beta_{max}$, but $\mu = 2\omega_{1}$ does not.  If $q_{2} \neq 0$, then $q_{1} = 0$ and $q_{2} = 1$.  In this case $\mu = \omega_{2} = \beta_{max}$, which does satisfy the condition $2\mu \leq 2\beta_{max}$.

    We conclude that $\mu = \omega_{1} = \beta_{min}$ and $\mu = \omega_{2} = \beta_{max}$ are the only solutions.  This completes the proof of Lemma 11.4. $\square$
\newline

$\mathit{Proof~of~Lemma~11.5}$  Let $\omega_{i}$, 1 $\leq$ i $\leq$ n, be a dominant weight such that $2\omega_{i} = \alpha + \beta$ for $\alpha, \beta \in \Phi$.
\newline

$\mathbf{Sublemma~1}$  Either $\alpha \in \Phi^{+}$ or $\beta \in \Phi^{+}$.

$\mathit{Proof}$ We know that $\omega_{i} =$ \(\sum_{j=1}^{n} C^{ij} \alpha_{j}\), where (C$^{ij}$) denotes the inverse Cartan matrix.  Recall that C$^{ij} > 0$ for all i,j.  Write $\alpha =$ \(\sum_{j=1}^{n} a_{j} \alpha_{j}\) and $\beta =$ \(\sum_{j=1}^{n} b_{j} \alpha_{j}\) for suitable integers a$_{j}$,b$_{j}$.  Since $2\omega_{i} = \alpha + \beta$ it follows that $C^{ij} = a_{j} + b_{j}$ for 1 $\leq$ j $\leq$ n.  Hence for each j either a$_{j} > 0$ or b$_{j} > 0$.  Since the integers \{a$_{j}$\}  all have the same sign, as do the integers \{b$_{j}$\}, it follows that either $\alpha \in \Phi^{+}$ or $\beta \in \Phi^{+}$.
\newline

    The next fact is well known but not easy to find in the literature.
\newline

$\mathbf{Sublemma~ 2}$  Suppose that there exist roots $\alpha, \beta \in \Phi$ such that $ < \alpha, \beta > = 3$.  Then $\fG = G_{2}$.

$\mathit{Proof}$  Define $\Phi_{\alpha \beta} = \Phi~ \cap~ \bz$ - span ($\alpha, \beta$).  It is easy to check that $\Phi_{\alpha \beta}$ satisfies the axioms for a root system, and clearly $\Phi_{\alpha \beta}$ has rank 2.  The hypothesis and a short argument show that the complex simple Lie algebra with root system $\Phi_{\alpha \beta}$ is $G_{2}$ (cf. section 9.4 of [Hu]).  It follows that $\fG = G_{2}$ since the Dynkin diagram classification shows that $G_{2}$ is the only root system of a complex simple Lie algebra that contains $G_{2}$ as a sub root system. $\square$
\newline

$\mathbf{Sublemma~3}$ Let $2\omega_{i} = \alpha + \beta$, where $\alpha \in \Phi^{+}$ and $\beta \in \Phi$.  Then there exists w $\in$ W such that w$(\omega_{i}) = \omega_{i}$ and w($\alpha$) $= \beta_{max}$ or $\beta_{min}$.

$\mathit{Proof}$  It is easy to prove this result by hand in the case that $\fG = G_{2}$, and we omit the details.  In the remainder of the proof we assume that $\fG \neq G_{2}$.  In particular we may assume that $\linebreak | < \alpha, \beta > | \leq 2$ for all $\alpha, \beta \in \Phi$ by Sublemma 2.

For an element $\alpha =$ \(\sum_{j=1}^{n} a_{j} \alpha_{j}\) $\in \Phi^{+}$ define $\ell(\alpha) =$ \(\sum_{j=1}^{n} a_{j}\) $ > 0$.  Let W$_{i} = \{w \in W : w(\omega_{i}) = \omega_{i}\}$.  It follows that $2\omega_{i} = w(\alpha) + w(\beta)$ for all w $\in$ W$_{i}$.  Choose $\alpha \in \Phi^{+}$ and $\beta \in \Phi$ such that $2\omega_{i} = \alpha + \beta$ and $\ell(\alpha) \geq \ell(\alpha'))$ if  $2\omega_{i} = \alpha' + \beta'$, where $\alpha' \in \Phi^{+}~ \cap~ W_{i}(\alpha)$ and $\beta' \in \Phi$.  It suffices to show that $\alpha$ is a dominant weight; that is, $ < \alpha, \alpha_{i} >~  \geq 0$ for 1 $\leq$ i $\leq$ n.  It will then follow that $\alpha = \beta_{max}$ or $\beta_{min}$ since these are the only two roots that are also dominant weights.  We use the fact that the W orbit of any abstract weight contains exactly one dominant weight (cf. section 13.2 of [Hu]).

    If j $\neq$ i, then $\sigma_{\alpha_{j}} \in$ W$_{i}$ and we obtain $2\omega_{i} = \sigma_{\alpha_{j}}(\alpha) + \sigma_{\alpha_{j}}(\beta)$.  Now $\sigma_{\alpha_{j}}(\alpha) = \alpha - < \alpha, \alpha_{j} > \alpha_{j}$.  If $ < \alpha, \alpha_{j} >$ is negative, then $ - < \alpha, \alpha_{j} >~ \geq 1,\sigma_{\alpha_{j}}(\alpha) \in \Phi^{+}~ \cap~ W_{i}(\alpha)$ and $\ell(\sigma_{\alpha_{j}}(\alpha)) = \ell(\alpha) - < \alpha, \alpha_{j} >~ \geq \ell(\alpha) + 1$.  This contradicts the maximality of $\ell(\alpha)$ among the solutions to  $2\omega_{i} = \alpha' + \beta'$, where $\alpha' \in \Phi^{+}~ \cap~ W_{i}(\alpha)$ and $\beta' \in \Phi$.  We conclude

\indent \indent (1) $\hspace{.5in} < \alpha, \alpha_{j} >~~ \geq 0$ if j $\neq$ i.

    Next, we assume that $ < \alpha, \alpha_{i} >~ <~ 0$ and obtain a contradiction.  If $ < \alpha, \alpha_{i} >~ <~ 0$, then  $ < \alpha, \alpha_{i} >~ \leq -1$, which implies that $ < \beta, \alpha_{i} > = 2 <\omega_{i}, \alpha_{i} > - < \alpha, \alpha_{i} >~ \geq 3$.  It is well known that $ < \beta, \sigma >~ \leq 3$ for all $\beta, \sigma \in \Phi$ (cf. section 9.4 of [Hu]).  Hence  $ < \beta, \alpha_{1} > = 3$.  However, by Sublemma 2 this contradicts our assumption at the beginning of the proof that $\fG \neq G_{2}$.  Therefore  $ < \alpha, \alpha_{1} > \geq 0$ and we conclude

\indent \indent (2) $\hspace{,5in}  < \alpha, \alpha_{i} >~~ \geq 0$.

    From (1) and (2) we conclude that $\alpha$ is a dominant root, which proves that $\alpha = \beta_{max}$ or $\beta_{min}. \square$
\newline

$\mathbf{Sublemma~4}$  Let $2\omega_{i} = \alpha + \beta$, where $\alpha = \beta_{max}$ or $\beta_{min}$ and $\beta \in \Phi$.  Then $\alpha$ and $\beta$ belong to the list in (11.5).  In particular $\beta \in \Phi^{+}$ in all cases.

$\mathit{Proof}$  If  $2\omega_{i} = \alpha + \beta$, then $\mu = \omega_{i}$ satisfies the conditions $2\mu \in \bz-span(\Delta)$ and $2\mu \leq 2\beta_{max}$.  Hence $\omega_{i}$ belongs to the list of (11.4).  If $\omega_{i}$ belongs to the list of (11.4), then straightforward computations now show that roots $\beta \in \Phi$ of the form $\beta = 2\omega_{i} - \beta_{max}$ or  $\beta = 2\omega_{i} - \beta_{min}$ must be precisely those listed in (11.5).  An examination of the list in (11.5) shows that $\beta \in \Phi^{+}$ in all cases.$\square$
\newline

    We now complete the proof of (11.5).  If $2\omega_{i} = \alpha + \beta$ for $\alpha, \beta \in \Phi$, then by sublemmas 1, 3 and 4 there exists an element w of W such that $w(\omega_{i}) = \omega_{i}, \alpha' = w(\alpha) \in \Phi^{+}$ and $\beta' = w(\alpha) \in \Phi^{+}$.  Hence $2\omega_{i} = \alpha' + \beta'$, where $\alpha'$ and $\beta'$ both lie in $\Phi^{+}$.  By sublemmas 3 and 4 there exists an element  $w_{1}^{*}$ of W such that $w_{1}^{*}(\omega_{i}) = \omega_{i},~ w_{1}^{*}(\alpha') = \beta_{max}$ or $\beta_{min}$ and  $w_{1}^{*}(\beta') \in \Phi^{+}$.  If $w_{1} = w_{1}^{*}w$, then $w_{1}$ satisfies assertions 1) and 2) of (11.5).  Similarly, there exists an element $w_{2}$ of W that satisfies assertions 1) and 3) of (11.5).  The list of (11.5) now follows from sublemma 4.  This completes the proof of Lemma 11.5.  $\square$
\newline

$\mathit{Proof~ of~ 9.2}$

    We first consider the special case of (9.2) where k $= 0$.

\begin{lemma}  Let $\fG$ be a complex simple Lie algebra.  Let $\mu \in \tilde{\Lambda}$ be an abstract weight such that $2\mu \in \Phi$.  Then all Weyl group orbits W($\mu$) are exhausted by the following possibilities :

A$_{1} \hspace{.7in} \mu = 1 = \frac{1}{2} \beta_{max} \in \bz$.

C$_{n}, n \geq 2 \hspace{.3in} \mu = \omega_{1} = \frac{1}{2} \beta_{max}$.

If $\mu$ is a dominant weight, then $\mu$ is one of the elements listed above.
\end{lemma}

\begin{proof}  From the hypothesis it is clear that $2\mu \leq \beta_{max} \leq 2\beta_{max}$.  Hence we need only consider those abstract weights $\mu$ that are listed in (11.4).  We may exclude the abstract weights $\mu = \beta_{max}$ or $\mu = \beta_{min}$ since $2\sigma \notin \Phi$ for any $\sigma \in \Phi$ by one of the axiomatic properties of root systems.  By inspection of the list in (11.4) this eliminates the Lie algebras E$_{6}$, E$_{7}$, E$_{8}$, F$_{4}$, G$_{2}$.  Of the remaining examples of $\mu$ in (11.4) all fail the condition $2\mu \leq \beta_{max}$ except the two examples listed above.

    The orbit W($\mu$) contains exactly one dominant weight (cf. section 13.2 of [Hu]).  The final assertion of the lemma now follows since the weights listed in the lemma are dominant.
\end{proof}

    The proof of (9.2) is an immediate consequence of the result just proved and the next three results.

\begin{lemma}  Let $\fG$ be a complex simple Lie algebra.  Let $\mu \in \tilde{\Lambda}$ be an abstract weight such that $2\mu = 2k\alpha + \beta$ for elements $\alpha, \beta \in \Phi$.  If we replace $\mu$ by w$(\mu)$ for a suitable element w of W, then there exist an element $\alpha' \in \Phi$ and a dominant abstract weight  $\mu' \in \tilde{\Lambda}$ such that $\mu = \mu' + k\alpha'$, where $2\mu' \in \Phi$.
\end{lemma}

\begin{proof}  Let $\mu' = \mu - k\alpha$.  Clearly $2\mu' = \beta \in \Phi$.  Now choose an element w of W such that w($\mu'$) is a dominant weight.
\end{proof}

    For the next result we recall that $\Delta = \{\alpha_{1}, \alpha_{2}, ~...,~ \alpha_{n}\}$ denotes the set of positive simple roots for a complex semisimple Lie algebra $\fG$.

\begin{lemma}  Let $\fG = C_{n}$ for n $\geq 3$.  Let W$_{1}$ be the subgroup of W generated by $\{\sigma_{\alpha_{j}} : 2 \leq j \leq n \}$.  Then $\Phi^{+} - \beta_{max} \subset W_{1}(\alpha_{1}) \cup W_{1}(\alpha_{2}) \cup W_{1}(\alpha_{n})$ (disjoint union).
\end{lemma}

\begin{proof}  We note that $\Phi^{+} - \beta_{max}$ is invariant under W$_{1}$ since W$_{1}$ fixes $\omega_{1}$ and $\beta_{max} = 2\omega_{1}$.   We show first that the orbits $W_{1}(\alpha_{1}), W_{1}(\alpha_{2})$ and $W_{1}(\alpha_{n})$ are disjoint.  The orbits $W_{1}(\alpha_{2})$ and $W_{1}(\alpha_{n})$ are disjoint since $\alpha_{2}$ is a short root and $\alpha_{n}$ is a long root.  Note that the $\alpha_{1}$ coefficient is zero in the elements of $W_{1}(\alpha_{2})$ and $W_{1}(\alpha_{n})$ while the  $\alpha_{1}$ coefficient is one in the elements of $W_{1}(\alpha_{1})$.  Hence the orbit $W_{1}(\alpha_{1})$ is disjoint from the orbits $W_{1}(\alpha_{2})$ and $W_{1}(\alpha_{n})$.

    The remaining assertions of the lemma follow from

\indent 1)  $\{\alpha_{2}, ~...,~ \alpha_{n-1}\} \subset W_{1}(\alpha_{2})$.

\indent 2)  $\Phi^{+} - \beta_{max} \subset W_{1}(\Delta)$.

    The assertion 1) follows since $\Delta_{1} = \{\alpha_{2}, ~...,~ \alpha_{n} \}$ is a base for a root system isomorphic to C$_{n-1}$ and W$_{1}$ is isomorphic to the Weyl group for C$_{n-1}$.  The set $\{\alpha_{2}, ~...,~ \alpha_{n-1} \}$ consists of the short roots in $\Delta_{1}$, and the short roots consist of a single W$_{1}$ orbit.

    It remains to prove 2).  For $\beta =$ \(\sum_{i=1}^{n} b_{i} \alpha_{i} \) $\in \Phi^{+}$ we define $\ell(\beta) =$ \(\sum_{i=1}^{n} b_{i} \in \bz^{+}$.

    We suppose that 2) is false and we obtain a contradiction.  Let $\alpha \in \Phi^{+} - \beta_{max}$ be an element such that $W_{1}(\alpha) \cap \Delta$ is empty.  Choose $\alpha' \in W_{1}(\alpha) \cap \Phi^{+}$ such that $\ell(\alpha') \leq \ell(w(\alpha'))$ for all w $\in W_{1}$ with $w(\alpha') \in \Phi^{+}$.  Since $\alpha' \notin \Delta$ it follows that $\alpha' \neq \alpha_{j}$ for 1 $\leq$ j $\leq$ n.  Hence $\sigma_{\alpha_{j}}(\alpha') \in \Phi^{+}$ for 2 $\leq$ j $\leq$ n by Lemma B in section 10.2 of [Hu].

    By the choice of $\alpha'$ we have $\ell(\alpha') \leq \ell(\sigma_{\alpha_{j}}(\alpha'))$ for 2 $\leq$ j $\leq$ n since $\sigma_{\alpha{_j}} \in W_{1}$.  It follows that $\langle \alpha', \alpha_{j} \rangle \leq 0$ for 2 $\leq$ j $\leq$ n since $\sigma_{\alpha_{j}}(\alpha') = \alpha' - \langle \alpha', \alpha_{j} \rangle \alpha_{j}$.  If we write $\alpha' =$\(\sum_{i=1}^{n} a_{i} \alpha_{i} \), where $a_{i} \in \bz^{+}$, then we obtain

\indent 1)  $\hspace{.2in} a)  -a_{n-1} + 2a_{n} = \langle \alpha', \alpha_{n} \rangle \leq 0$.

\indent $\hspace{.35in}      b)  -a_{n-2} + 2a_{n-1} - 2a_{n} = \langle \alpha', \alpha_{n-1} \rangle \leq 0$.

\indent $\hspace{.35in} c)-a_{i-1} + 2a_{i} - a_{i+1} =   \langle \alpha', \alpha_{i} \rangle \leq 0 \hspace{.2in}$ for 2 $\leq$ i $\leq$ n $-2$.

    Hence we obtain

\indent   2) $\hspace{.2in}     a) 2a_{n} \leq a_{n-1}$.

\indent $\hspace{.35in}     b) 2a_{n-1} \leq a_{n-2} + 2a_{n}$.

\indent $\hspace{.35in}     c) 2a_{i} \leq a_{i-1} + a_{i+1} \hspace{.2in}$ for 2 $\leq$ i $\leq$ n $-2$.

    Since $\alpha' \leq \beta_{max}$ we have

\indent 3) $\hspace{.2in}   0 \leq a_{j} \leq 2$ for $1 \leq j \leq n-1$ and $0 \leq a_{n} \leq 1$.

    We consider separately the cases $a_{n} = 1$ and $a_{n} = 0$.  In the first case we show that $\alpha' = \beta_{max}$, and in the second we show that $\alpha' = \alpha_{1}$.  However, by assumption $\alpha' \in \Phi^{+} - \beta_{max}$ and $\alpha' \notin \Delta$.  This contradiction will imply that $\Phi^{+} - \beta_{max} \subset W_{1}(\Delta)$.

    $\mathbf{Case~ 1}~ a_{n} = 1$

    From 2a) and 3) we have $2 \leq a_{n-1} \leq 2$, which implies that $a_{n-1} = 2$.  From 2b), 3)  and 4) we obtain $4 = 2a_{n-1} \leq a_{n-2} + 2a_{n} = a_{n-2} + 2 \leq 4$, which implies that $a_{n-2} =2$.  Continuing in this fashion yields $a_{i} = 2$ for $1 \leq i \leq n-1$.  Therefore $\alpha = \beta_{max}$.

    $\mathbf{Case~2}~ a_{n} = 0$

    If $a_{n-1} \neq 0$, then from 2b) and 3) we have $2 \leq 2 a_{n-1} \leq a_{n-2} + 2a_{n} = a_{n-2} \leq 2$.  Hence $a_{n-1} = 1$ and $a_{n-2} = 2$.  From 2c) we obtain $4 = 2a_{n-2} \leq a_{n-3} + a_{n-1} = a_{n-3} + 1$.  This implies that $a_{n-3} \geq 3$, which contradicts 3).  We conclude that $a_{n-1} = 0$.

    Continuing in this fashion we obtain $a_{i} = 0$ for $3 \leq i \leq n$.  If $a_{2} \neq 0$, then from 2c) we obtain $2 \leq 2a_{2} \leq a_{1} + a_{3} = a_{1} \leq 2$.  This implies that $a_{2} = 1$ and $a_{1} = 2$.  Hence $\alpha = 2\alpha_{1} + \alpha_{2}$.  However, $\sigma_{\alpha_{1}}(2\alpha_{1} + \alpha_{2}) = \alpha_{2} - \alpha_{1}$ is not a root since the coefficients of $\alpha_{1}$ and $\alpha_{2}$ have opposite signs.  Hence $2\alpha_{1} + \alpha_{2}$ is not a root.  This contradiction proves that $a_{2} = 0$, and hence $\alpha = a_{1} \alpha_{1}$ with $0 \leq a_{1} \leq 2$.  Since $2\alpha_{1} \notin \Phi$ it follows that $\alpha = \alpha_{1}$    \end{proof}

    An obvious modification of the proof of the preceding result yields the next result, which completes the proof of (9.2).

\begin{lemma}  If $\fG = C_{2}$, then $\Phi^{+} - \beta_{max} \subset W_{1}(\alpha_{1}) \cup W_{1}(\alpha_{2})$.
\end{lemma}

$\mathit{References}$

[B-tD]  T. Br$\ddot{o}$cker and T. tom Dieck, $\mathit{Representations~ of~ Compact~ Lie~ Groups}$, Springer, 1985, New York.

[D]  R. DeCoste,  Ph.D. dissertation, Univ. N. Carolina Chapel Hill, 2004.

[CG]  L. Corwin and F. Greenleaf, $\mathit{Representations~ of~ Nilpotent~ Lie~ Groups~ and~ their~ Applications}$, Cambridge University
Press, 1990, Cambridge.

[Eb1]  P. Eberlein, " The moduli space of 2-step nilpotent Lie algebras of type (p,q) ", Contemporary Math. 332 (2003), 37-72.

[Eb2] -------, " Riemannian submersions and lattices in 2-step nilpotent Lie groups ", Commun. Analysis and Geom. vol.11 (3), (2003), 441-488.

[Eb3] -------, " Geometry of 2-step nilpotent Lie groups with a left invariant metric ", Ann. Sci. $\acute{E}$cole Norm. Sup. 27, (1994), 611-660.

[El], A. G. Elashvili, " Stationary subalgebras of points of the common state for irreducible linear Lie groups ", Func. Analysis and Appl. 6 (1972), 139-148.

[FH]  W. Fulton and J. Harris, $\mathit{Representation~ Theory,~ A~ First~ Course}$, GTM $\#$ 129, Springer, 1991, New York.

[GM]    C. Gordon and Y. Mao, " Geodesic conjugacies in two step nilmanifolds ", Mich. Math. J. 45 (1998), 451-481.

[GW]    C. Gordon and E. Wilson, " Isospectral deformations of compact solvmanifolds ", J. Diff. Geom. 19 (1984), 241-256.

[He]  S. Helgason, $\mathit{Differential~ Geometry,~ Lie~ Groups~ and~ Symmetric~ Spaces}$, Academic Press, 1978, New York.

[HK]  K. Hoffman and R. Kunze, $\mathit{Linear~ Algebra}$, Prentice Hall, 1961, Englewood Cliffs.

[Hu]  J. Humphreys, $\mathit{Introduction~ to~ Lie~ Algebras~ and~ Representation~ Theory}$, Springer, 1972, New York.

[KL]  F. Knop and P. Littelmann, " Der grad erzeugender Funktionen von Invariantentringen ", Math Zeit. 96 (1987), 211-229.

[L]  J. Lauret, " Homogeneous nilmanifolds attached to representations of compact Lie groups", Manuscripta Math. 99 (1999), 287-309.

[Ma]  A. I. Mal'cev, " On a class of homogeneous spaces ", Amer. Math. Soc. Translations, vol. 39, 1951.

[Mo]  D. Morris (= D.Witte) " Real representations of semisimple Lie algebras have $\bq$-forms ", $\mathit{Proc.~ International~ Conf.~ on~
Algebraic~ Groups~ and~ Arithmetic~}$, TIFR, December 2001, to appear.

[R1] M. S. Raghunathan, " Arithmetic lattices in semisimple Lie groups ", Proc. Indian Acad. Sci., 91(2), (1984), 133-138.

[R2]  --------, $\mathit{Discrete~ subgroups~ of~ Lie~ groups}$, Ergebnisse der Mathematik, vol.68, Springer, 1972.

\end{document}